\documentclass[preprint,3p,12pt]{elsarticle}
\usepackage{mathrsfs}
\usepackage{amsmath}
\usepackage{stmaryrd}
\usepackage{bbding}
\usepackage{dcolumn}
\usepackage{graphicx}
\usepackage{amsfonts}
\usepackage{amssymb}
\usepackage{psfrag}
\usepackage{wrapfig}
\usepackage{subfigure}
\usepackage{makeidx}
\usepackage{bm}
\usepackage{epsf}
\usepackage{epsfig}
\usepackage{setspace}
\usepackage{graphicx}
\usepackage{epstopdf}
\usepackage{psfrag}
\usepackage{subfigure}
\usepackage{color}

\begin{document}

\title{Two-Stage Fourth-order Gas-kinetic Scheme for Three-dimensional Euler and Navier-Stokes Solutions}

\author[iapcm]{Liang Pan}
\ead{panliangjlu@sina.com}

\author[HKUST1,HKUST2]{Kun Xu\corref{cor1}}
\ead{makxu@ust.hk}

\address[iapcm]{Institute of Applied Physics and Computational Mathematics, Beijing, 100088, China}
\address[HKUST1]{Department of Mathematics, Department of Mechanical and Aerospace Engineering,
Hong Kong University of Science and Technology, Clear Water Bay, Kowloon, Hong Kong}
\address[HKUST2]{HKUST Shenzhen Research institute, Shenzhen 518057, China}
\cortext[cor1]{Corresponding author}

\begin{abstract}
For the one-stage third-order gas-kinetic scheme (GKS), success applications have been achieved for the three-dimensional compressible
flow computations \cite{GKS-high-1}.
The high-order accuracy of the scheme is obtained
directly by integrating a multidimensional time-accurate gas distribution function over the cell interface within a time step
without implementing Gaussian quadrature points and Runge-Kutta time-stepping technique.
However, for the further increasing the order of the scheme, such as the fourth-order one,
the formulation becomes very complicated for the multidimensional flow.
Recently, a two-stage fourth-order GKS with high efficiency has been constructed
for two-dimensional inviscid and viscous flow computations \cite{GRP-high, GKS-high},
and the scheme uses the time accurate flux function and its time derivatives.
In this paper, a fourth-order GKS is developed for the three-dimensional flows under the two-stage framework.
Based on the three-dimensional WENO reconstruction and flux evaluation at Gaussian quadrature points on a cell interface,
the high-order accuracy in space is achieved first. Then, the two-stage time stepping method provides the high accuracy in time.
In comparison with the formal third-order GKS \cite{GKS-high-1},
the current fourth-order method not only improves the accuracy of the scheme, but also reduces the complexity of the gas-kinetic solver greatly. More importantly, the fourth-order GKS has the same robustness as the second-order shock capturing scheme.
This scheme is applied to both inviscid and viscous, and low and high speed flow computations.
Numerical results validate the outstanding reliability and applicability of the scheme for three-dimensional flows, such as
turbulent one.
With the count of the degree of freedom in simulating a flow field, the current two-stage fourth-order multidimensional GKS becomes
one of the most reliable and efficient higher-order schemes among all shock capturing higher-order schemes for computational fluid dynamics.

\end{abstract}
\maketitle

\section{Introduction}
Over the past half century, the computational fluid dynamics has been developed into a powerful tool for solving the fluid flow problems in industrial applications.  Currently, there are a gigantic amount of numerical methods in literature.
In comparison with well-developed second-order schemes, the higher-order methods can provide more accurate solutions,
but they are less robust and more complicated. For the high-order schemes, there are generally three parts,
i.e. spatial discretization, temporal discretization, and flux solvers. For the spatial discretization, many high-order methods have been developed, including the discontinuous Galerkin (DG) \cite{DG-1, DG-2,DG-3}, essential non-oscillatory (ENO) \cite{ENO-1,ENO-2}, weighted essential non-oscillatory (WENO) \cite{WENO,WENO-JS,WENO-Z}, and many others.
For most of those methods, the exact and approximate Riemann solvers \cite{Riemann-appro} are used for flux evaluation.
Due to the first-order evolution of Riemann solvers, the Runge-Kutta method is implemented to achieve higher order temporal accuracy \cite{TVD-rk}, in which $n$-stage is needed for $n$-th order accuracy. Instead of Riemann solver, many schemes have been developed based on the time-dependent flux function with high-order evolution, such as the generalized Riemann problem solvers \cite{GRP1,GRP2,GRP3,GRP-high-2}, gas-kinetic scheme \cite{GKS-Xu1,GKS-Xu2}, and AEDR methods \cite{Riemann-appro,ADER}. High-order temporal accuracy can be achieved in a one-stage framework
with the time integration of the time-dependent flux function.
Recently, in order to increase the accuracy and efficiency of these schemes
a two-stage fourth-order method has been developed for the time-dependent flux solvers \cite{GRP-high,GKS-high}, where both the flux and temporal derivative of flux function are used in the construction of higher-order schemes.
By combining the multi-stage multi-derivative technique, a family of higher-order schemes has been constructed as well \cite{MSMD-GKS}.

In the past decades, the gas-kinetic scheme (GKS) based on the Bhatnagar-Gross-Krook (BGK) model \cite{BGK-1,BGK-2,BGK-3} has been developed systematically for the compressible flow computations \cite{GKS-Xu1,GKS-Xu2,GKS-Xu3}.
Different from the traditional CFD methods based on the macroscopic governing equations, the main advantages of the gas-kinetic scheme are the followings. (i) The inviscid and viscous coupling in the flux evolution \cite{GKS-Xu1, GKS-Xu2}; (ii) Multi-dimensionality with the inclusion of both normal and tangential gradients of flow variables in the flux function across a cell interface \cite{GKS-high-4}; (iii) Compact stencils can be constructed with the use of the time accurate cell interface flow variables at the next time level \cite{GKS-high-4}; (iv) Extension to the whole flow regime from the rarefied to the continuum one \cite{UGKS-1,UGKS-2}. Recently, with the high-order initial reconstruction, the third-order gas-kinetic schemes have been constructed \cite{GKS-high-2,GKS-high-3,GKS-high-4}, in which the flux evaluation is based on the moments of spatial-temporal dependent gas distribution function.
High-order accuracy can be achieved in a one-stage scheme without Gaussian point integration for spatial accuracy
and Runge-Kutta technique for temporal accuracy.
However, with the one-stage gas evolution model, the formulation of gas-kinetic scheme can become very complicated for the further development,
such as the fourth-order method \cite{GKS-high-5}, especially for multidimensional computations.
The two-stage fourth-order temporal discretization for time-accurate flux solvers in \cite{GRP-high} provides a reliable framework to further develop the GKS into fourth-order and even higher accuracy with the implementation of the traditional second-order or third-order flux functions \cite{GKS-high, GKS-high-6,MSMD-GKS}. Most importantly, this scheme is robust, and works perfectly from the subsonic to the hypersonic viscous heat conducting flows.
The robustness of the GKS in comparison with Riemann solver based CFD methods is solely due to the
differences in the dynamical evolution model of the flux function. The GKS flux follows the dynamics from the particle free transport, to including  collisional effect, and to the NS distribution function with the count of intensive particle collisions as the ratio of the time step $\Delta t$ over the particle collision time $\tau$ becoming large. The real physics used in the flux depends on the local $\Delta t /\tau$ \cite{xu-liu}.
In real NS computations for the compressible high Mach number flow, the ratio of $\Delta t /\tau$ is not too large as people think of for the
validation of the NS modeling \cite{liu-kn}.
However, for the Riemann solver based CFD methods, at the beginning of the step it is already assumed that there are infinity number of particle collision to generate
distinguishable waves in the Riemann solution, and the collision needs to be reduced for the NS solutions. Theoretically, there is no such a physical process for the Riemann solver-based Godunov-type schemes for the NS equations.
For the second-order schemes, it is hard to distinguish the dynamical differences from the GKS and Riemann solver.
However,  for the higher-order schemes it seems that a reliable physical evolution model becomes more
important due to the absence of large numerical dissipation in the second-order schemes, and the delicate flow structures captured in higher-order
schemes depend on the quality of the solvers greatly \cite{MSMD-GKS}.

In this paper, with the two-stage fourth-order discretization, a multidimensional fourth-order gas-kinetic scheme is constructed for simulating  three-dimensional flows. High-order accuracy in space is achieved by the three-dimensional WENO method \cite{WENO,WENO-JS,WENO-Z} and Gaussian quadrature points for the numerical fluxes. In comparison with the formal three-dimensional scheme \cite{GKS-high-1}, the current fourth-order scheme reduces the complexity of the gas-kinetic flux solver greatly, and improves the robustness of scheme.  Many numerical tests, including both inviscid and viscous, and low and high speed flow computations, will be used to validate the current fourth-order
method.  Numerical results show that the current scheme has the same reliability and applicability as the well-developed second-order scheme, but is much more accurate and effective to capture the complicated flow structures. With the count of the degrees of freedom for the description of a flow field, the current scheme provides state-of-art solutions from a higher-order scheme from the incompressible to the hypersonic flow simulations.

This paper is organized as follows. In Section 2, BGK equation and finite volume scheme are briefly reviewed. The general formulation for the two-stage temporal discretization is introduced in Section 3, and the procedure of spatial reconstruction is given in Section 4. Section 5 includes numerical examples to validate the current algorithm. The last section is the conclusion.

\section{BGK equation and finite volume scheme}
The three-dimensional BGK equation \cite{BGK-1,BGK-2,BGK-3} can be written as
\begin{equation}\label{bgk}
f_t+uf_x+vf_y+wf_z=\frac{g-f}{\tau},
\end{equation}
where $f$ is the gas distribution function, $g$ is the corresponding equilibrium state, and $\tau$ is the collision time. The collision term satisfies the compatibility condition
\begin{equation}\label{compatibility}
\int \frac{g-f}{\tau}\psi \text{d}\Xi=0,
\end{equation}
where $\psi=(1,u,v,w,\displaystyle \frac{1}{2}(u^2+v^2+w^2+\xi^2))$, the internal variables $\xi^2$ is equal to $\xi^2=\xi_1^2+...+\xi_K^2$, $\text{d}\Xi=\text{d}u\text{d}vd\text{d}w\text{d}\xi^1...\text{d}\xi^{K}$, $K$ is the degrees of freedom, i.e. $K=(5-3\gamma)/(\gamma-1)$ for three-dimensional flows and $\gamma$ is the specific heat ratio. In the continuum region, the gas distribution function can be expanded as
\begin{align*}
f=g-\tau D_{\textbf{u}}g+\tau D_{\textbf{u}}(\tau D_{\textbf{u}})g-\tau D_{\textbf{u}}[\tau D_{\textbf{u}}(\tau D_{\textbf{u}})g]+...,
\end{align*}
where $D_{\textbf{u}}=\displaystyle\frac{\partial}{\partial t}+\textbf{u}\cdot \nabla$. Based on the Chapman-Enskog expansion, the macroscopic equations can be derived \cite{GKS-Xu1,GKS-Xu2}. For the Euler equations, the zeroth-order truncation is taken, i.e. $f=g$. For the Navier-Stokes equations, the first-order truncation is
\begin{align*}
f=g-\tau (ug_x+vg_y+wg_z+g_t).
\end{align*}
With the higher order truncations, the Burnett and super-Burnett equations can be derived.

Taking moments of Eq.\eqref{bgk} and integrating over the control volume $V_{ijk}=\overline{x}_i\times\overline{y}_j\times \overline{z}_k$ with $\overline{x}_i=[x_i-\Delta x/2,x_i+\Delta x/2], \overline{y}_j=[y_j-\Delta y/2,y_j+\Delta y/2], \overline{z}_k=[z_k-\Delta z/2,z_k+\Delta z/2]$,  the
semi-discretized form of finite volume scheme can be written as
\begin{align}\label{finite}
&\frac{\text{d}Q_{ijk}}{\text{d}t}=\mathcal{L}(Q_{ijk})=\frac{1}{|V_{ijk}|}\Big[
\int_{\overline{y}_j\times\overline{z}_k}(F_{i-1/2,j,k}-F_{i+1/2,j,k})\text{d}y\text{d}z\nonumber\\
&+\int_{\overline{x}_i\times\overline{z}_k}(G_{i,j-1/2,k}-G_{i,j+1/2,k})\text{d}x\text{d}z+\int_{\overline{x}_i\times\overline{y}_j}(H_{i,j,k-1/2}-H_{i,j,k+1/2})\text{d}x
\text{d}y\Big],
\end{align}
where $Q=(\rho,\rho U,\rho V,\rho W,\rho E)$ are the conservative flow variables, $Q_{ijk}$ is the cell averaged value over control volume $V_{ijk}$ and $|V_{ijk}|=\Delta x\Delta y\Delta z$.
For the three-dimensional computation, the Gaussian quadrature for the numerical fluxes is used to achieve the accuracy in space, and the numerical fluxes in $x$-direction is given as an example
\begin{align}\label{gauss}
\int_{\overline{y}_j\times\overline{z}_k}F_{i+1/2,j,k}\text{d}y\text{d}z=\Delta y\Delta z\sum_{m,n=1}^M\omega_{mn} F(\textbf{x}_{i+1/2,j_m,k_n},t),
\end{align}
where $\textbf{x}_{i+1/2,,m,n}=(x_{i+1/2},y_{j_m},z_{k_n})$, $(y_{j_m},z_{k_n}), m, n= 1,...,M$ are the Gauss quadrature points for $\overline{y}_j\times\overline{z}_k$ and $\omega_{mn}$ are corresponding weights, and $M=2$ is used in this paper. Based on the spatial reconstruction, which will be presented in the following section, the reconstructed point value and the spatial derivatives at each Gauss quadrature points can be obtained and the numerical fluxes $F(\textbf{x}_{i+1/2,j_m,k_n},t)$ can be provided by the flow solvers. Similarly, the numerical fluxes  in the $y$ and $z$-directions can be obtained as well
\begin{align*}
\int_{\overline{x}_i\times\overline{z}_k}G_{i,j+1/2,k}\text{d}x\text{d}z=\Delta x\Delta z\sum_{m,n=1}^M\omega_{mn} G(\textbf{x}_{i_m,j+1/2,k_n},t),\\
\int_{\overline{x}_i\times\overline{y}_j}H_{i,j,k+1/2}\text{d}x\text{d}y=\Delta x\Delta y\sum_{m,n=1}^M\omega_{mn} H(\textbf{x}_{i_m,j_n,k+1/2},t).
\end{align*}
In the gas-kinetic scheme, the numerical fluxes at the Gauss quadrature point can be obtained by taking moments of the gas distribution function
\begin{align}\label{quadrature-flux}
F(\textbf{x}_{i+1/2,j_m,k_n},t)=\int\psi u f(\textbf{x}_{i+1/2,j_m,k_n},t,\textbf{u},\xi)\text{d}u\text{d}v\text{d}w\text{d}\xi,
\end{align}
where $f(\textbf{x}_{i+1/2,j_m,k_n},t,\textbf{u},\xi)$ is provided by the integral solution of BGK equation Eq.\eqref{bgk} at the cell interface
\begin{align}\label{integral}
f(\textbf{x}_{i+1/2,j_m,k_n},t,\textbf{u},\xi)=&\frac{1}{\tau}\int_0^tg(\textbf{x}',t',\textbf{u},\xi)e^{-(t-t')/\tau}dt'+e^{-t/\tau}f_0(-\textbf{u}t,\xi),
\end{align}
and $\textbf{x}_{i+1/2,j_m,k_n}=\textbf{0}$ is the location of cell interface, $\textbf{u}=(u,v,w)$ is the particle velocity, $x_{i+1/2}=x'+u(t-t'), y_{j_m}=y'+v(t-t'), z_{k_n}=z'+w(t-t')$ is the trajectory of particles.   For the second-order scheme, the second-order gas-kinetic solver for the three-dimensional flows \cite{GKS-Xu2,GKS-Xu1} can be expressed as
\begin{align}\label{flux}
f(\textbf{x}_{i+1/2,m,n},t,\textbf{u},\xi)=&(1-e^{-t/\tau})g_0+((t+\tau)e^{-t/\tau}-\tau)(\overline{a}_1u+\overline{a}_2v+\overline{a}_3w)g_0\nonumber\\
+&(t-\tau+\tau e^{-t/\tau}){\bar{A}} g_0\nonumber\\
+&e^{-t/\tau}g_r[1-(\tau+t)(a_{1r}u+a_{2r}v+a_{3r}w)-\tau A_r)]H(u)\nonumber\\
+&e^{-t/\tau}g_l[1-(\tau+t)(a_{1l}u+a_{2l}v+a_{3l}w)-\tau A_l)](1-H(u)).
\end{align}
where the coefficients in Eq.\eqref{flux} can be determined by the spatial derivatives of macroscopic flow variables and the compatibility condition as follows
\begin{align*}
\displaystyle\langle a_1\rangle =\frac{\partial W }{\partial x},
\langle a_2\rangle =\frac{\partial W }{\partial y}, \langle
a_3\rangle =\frac{\partial W }{\partial z},
\langle A+a_1u+a_2v+a_3w\rangle=0,
\end{align*}
where the superscripts or subscripts of these coefficients $a_1, a_2, a_3, A$ are omitted for simplicity, more details about the
determination of coefficient can be found in \cite{GKS-Xu2}. With the second-order gas-kinetic solver Eq.\eqref{flux}, the second-order accuracy in time can be achieved by one step integration.
In the one-stage gas evolution model, the third-order gas-kinetic solver has been developed as well \cite{GKS-high-1,GKS-high-2,GKS-high-3,GKS-high-4}.

\section{Fourth-order temporal discretization}
In the further development of higher-order scheme, the formulation of one-stage gas-kinetic solver can become very complicated, such as the fourth-order method \cite{GKS-high-5}, especially for multidimensional computations.
Recently, a two-stage fourth-order time-accurate discretization was developed for Lax-Wendroff flow solvers, particularly applied for hyperbolic equations with the generalized Riemann problem (GRP) solver \cite{GRP-high} and gas-kinetic scheme \cite{GKS-high}. Such method provides a reliable framework to develop a
three-dimensional fourth-order gas-kinetic scheme with a second-order flux function  Eq.\eqref{flux}. Consider the following time-dependent equation
\begin{align*}
\frac{\partial \textbf{w}}{\partial t}=\mathcal {L}(\textbf{w}),
\end{align*}
with the initial condition at $t_n$, i.e.,
\begin{align*}
\textbf{w}(t=t_n)=\textbf{w}^n,
\end{align*}
where $\mathcal {L}$ is an operator for spatial derivative of flux.
The time derivatives are obtained using the Cauchy-Kovalevskaya method,
\begin{align*}
\frac{\partial \textbf{w}^n}{\partial t}=\mathcal{L}(\textbf{w}^n),~ \ \ \
\frac{\partial }{\partial t}\mathcal
{L}(\textbf{w}^n)=\frac{\partial }{\partial \textbf{w}}\mathcal
{L}(\textbf{w}^n)\mathcal {L}(\textbf{w}^n).
\end{align*}
Introducing an intermediate state at $t^*=t_n+\Delta t/2$,
\begin{align}\label{step1}
\textbf{w}^*=\textbf{w}^n+\frac{1}{2}\Delta t\mathcal
{L}(\textbf{w}^n)+\frac{1}{8}\Delta t^2\frac{\partial}{\partial
t}\mathcal{L}(\textbf{w}^n),
\end{align}
the corresponding time derivatives are obtained as well for the intermediate stage state,
\begin{align*}
\frac{\partial \textbf{w}^*}{\partial t}=\mathcal{L}(\textbf{w}^*),~
\frac{\partial }{\partial t}\mathcal
{L}(\textbf{w}^*)=\frac{\partial }{\partial \textbf{w}}\mathcal
{L}(\textbf{w}^*) \mathcal {L}(\textbf{w}^*).
\end{align*}
Then, the state $\textbf{w}$ can be updated with the following formula,
\begin{align}\label{step2}
\textbf{w}^{n+1}=\textbf{w}^n+\Delta t\mathcal
{L}(\textbf{w}^n)+\frac{1}{6}\Delta t^2\big(\frac{\partial}{\partial
t}\mathcal{L}(\textbf{w}^n)+2\frac{\partial}{\partial
t}\mathcal{L}(\textbf{w}^*)\big).
\end{align}
It can be proved that for hyperbolic equations  the above time
stepping method  Eq.\eqref{step1} and Eq.\eqref{step2}  provides a fourth-order time
accurate solution for $\textbf{w}(t)$ at $t=t_n +\Delta t$. The
details of the analysis can be found in \cite{GRP-high}.

The semi-discretized form of finite volume scheme
Eq.\eqref{finite} can be applied for the above two-stage method
. For the gas-kinetic scheme, the gas evolution is a relaxation process from kinetic to hydrodynamic scale through the exponential function,
and the corresponding flux is a complicated function of time.
In order to obtain the time derivatives of the
flux function at $t_n$ and $t_*=t_n + \Delta t/2$ with the correct physics,
the flux function should be approximated as a linear function of time within a time interval.
According to the numerical fluxes at the Gauss quadrature points Eq.\eqref{quadrature-flux},  the following notation is introduced
\begin{align*}
\mathbb{F}_{i+1/2,j,k}(Q^n,\delta)=\int_{t_n}^{t_n+\delta}\textbf{F}_{i+1/2,j,k,}(Q^n,t)\text{d}t
=\sum_{m,n=1}^M\omega_{mn}\int_{t_n}^{t_n+\delta}F(\textbf{x}_{i+1/2,j_m,k_n},t)\text{d}t.
\end{align*}
In the time interval $[t_n, t_n+\Delta t]$, the flux is expanded as
the following linear form
\begin{align*}
\textbf{F}_{i+1/2,j,k}(Q^n,t)=\textbf{F}_{i+1/2,j,k}^n+ \partial_t \textbf{F}_{i+1/2,j,k}^n(t-t_n).
\end{align*}
The coefficients $\textbf{F}_{i+1/2,j,k}(Q^n,t_n)$ and $\partial_t\textbf{F}_{i+1/2,j,k}(Q^n,t_n)$ can be
determined as follows,
\begin{align*}
\textbf{F}_{i+1/2,j,k}(Q^n,t_n)\Delta t&+\frac{1}{2}\partial_t
\textbf{F}_{i+1/2,j,k}(Q^n,t_n)\Delta t^2 =\mathbb{F}_{i+1/2,j,k}(Q^n,\Delta t) , \\
\frac{1}{2}\textbf{F}_{i+1/2,j,k}(Q^n,t_n)\Delta t&+\frac{1}{8}\partial_t
\textbf{F}_{i+1/2,j,k}(Q^n,t_n)\Delta t^2 =\mathbb{F}_{i+1/2,j,k}(Q^n,\Delta t/2).
\end{align*}
By solving the linear system, we have
\begin{align*}
\textbf{F}_{i+1/2,j,k}(Q^n,t_n)&=(4\mathbb{F}_{i+1/2,j,k}(Q^n,\Delta t/2)-\mathbb{F}_{i+1/2,j,k}(Q^n,\Delta t))/\Delta t,\nonumber\\
\partial_t \textbf{F}_{i+1/2,j,k}(Q^n,t_n)&=4(\mathbb{F}_{i+1/2,j,k}(Q^n,\Delta t)-2\mathbb{F}_{i+1/2,j,k}(Q^n,\Delta t/2))/\Delta
t^2.
\end{align*}
Similarly, $\displaystyle \textbf{F}_{i+1/2,j,k}(Q^*,t_*), \partial_t
\textbf{F}_{i+1/2,j,k}(Q^*,t_*)$ for the intermediate state $t^*=t_n+\Delta t/2$ can be constructed. The corresponding fluxes in the $y$- and $z$-direction can be
obtained as well. With these notations, the numerical scheme for the three-dimensional flows can be constructed, and the detail procedure of two-stage algorithm can be found in \cite{GKS-high}.

\begin{figure}[!h]
\centering
\includegraphics[height=0.4\textwidth]{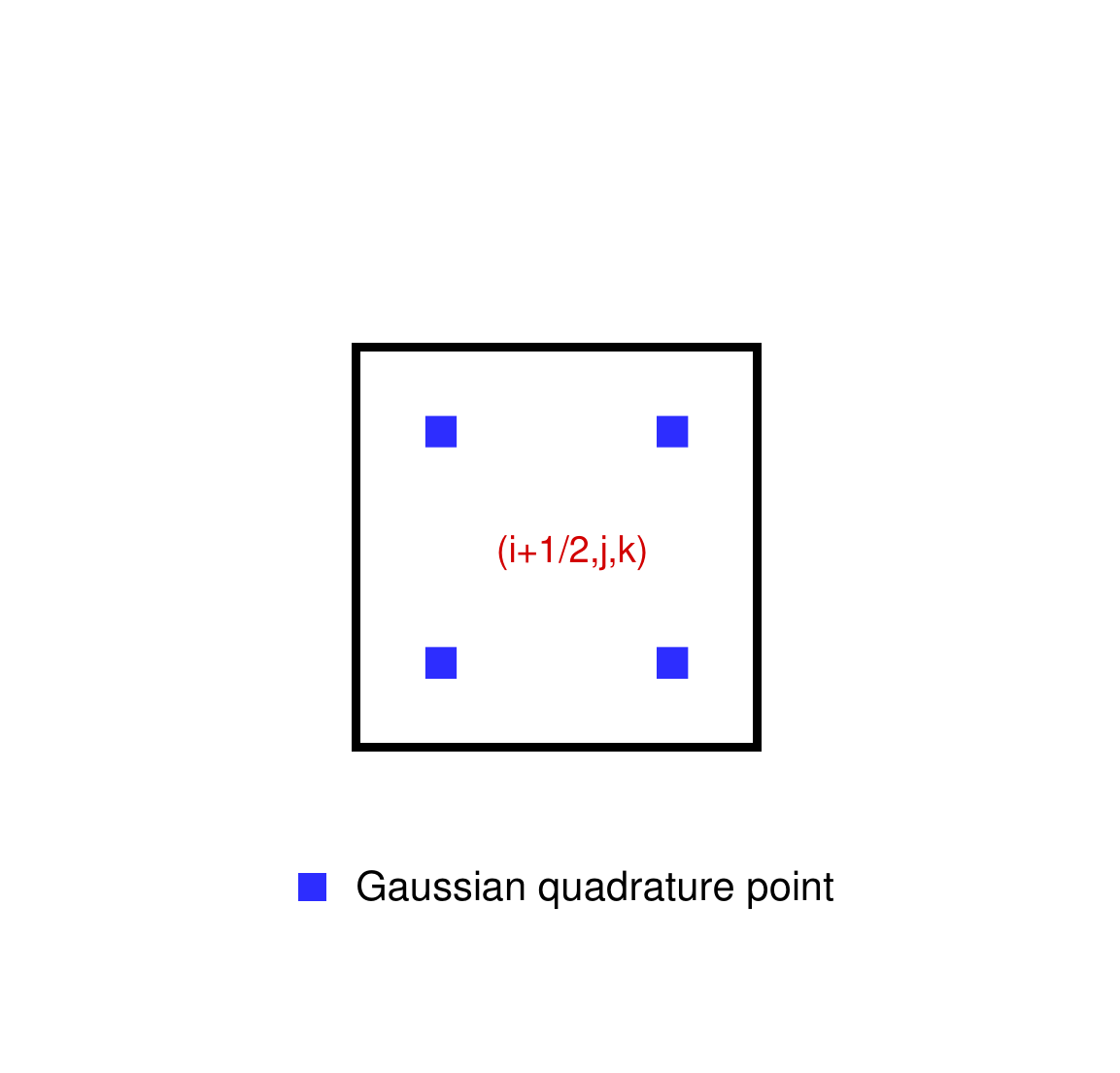}
\includegraphics[height=0.4\textwidth]{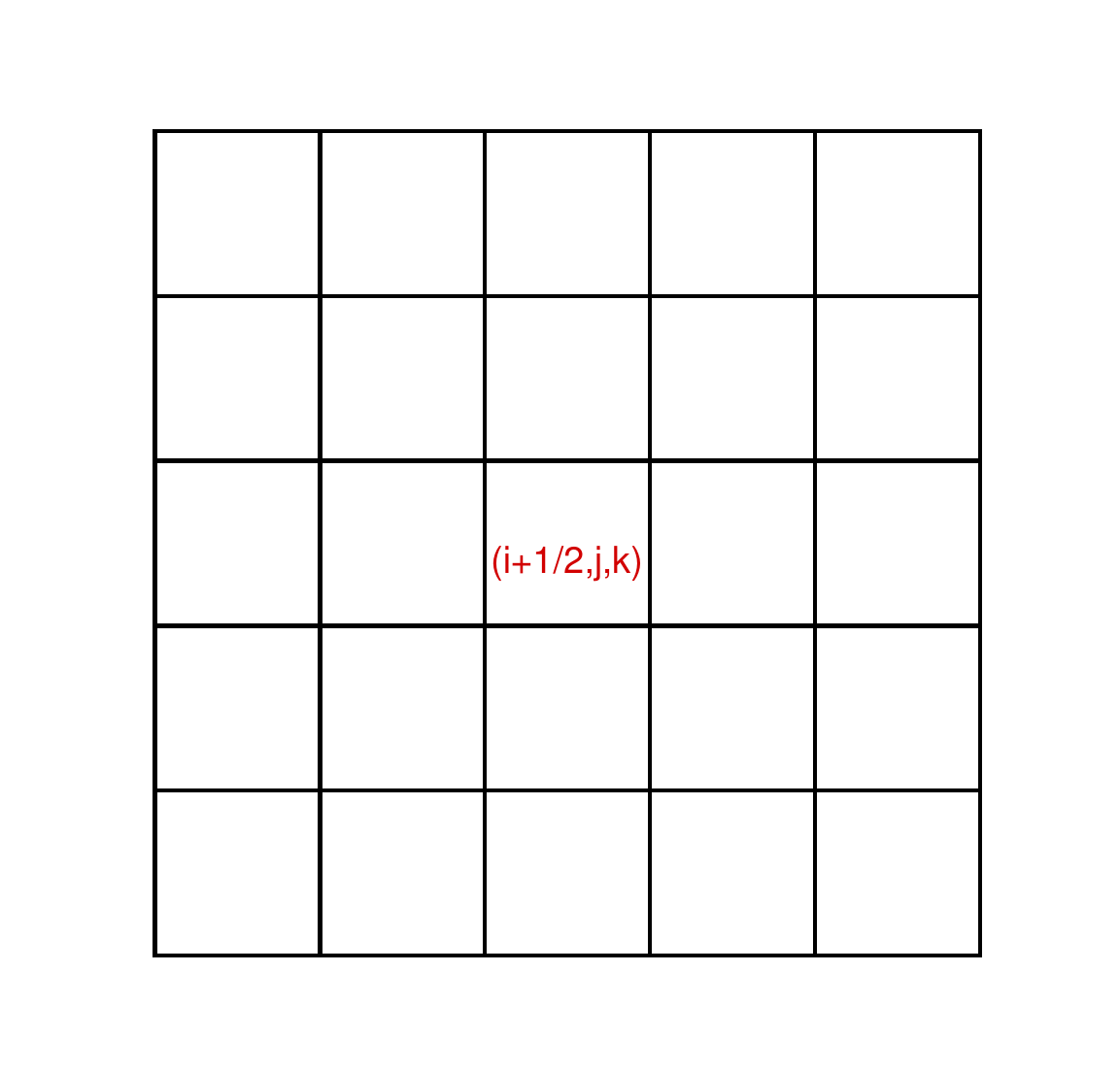}
\caption{\label{interplate}Schematics of spatial reconstruction at the cell interface $(i+1/2,j,k)$.}
\end{figure}

\section{Spatial reconstruction}
The above time evolution solution is based on the high-order initial reconstruction for macroscopic flow variables and  WENO reconstruction \cite{WENO, WENO-JS, WENO-Z} is adopted for the spatial reconstruction.
For the three dimensional computation, the reconstruction procedure for the cell interface $x_{i+1/2,j,k}$ is given as an example, and the stencil for reconstruction is given in Fig.\ref{interplate}.
The point value $Q_l, Q_r$ and $Q_0$ and first-order derivatives at the Gauss quadrature points $(x_{i+1/2},y_{j_m},z_{k_n})$, $m, n= 1,...,2$ need to be constructed. The detailed procedure is given as follows
\begin{enumerate}
\item According to one dimensional reconstruction, the cell averaged reconstructed value $(Q_{l})_{j-\ell_1,k-\ell_2}, (Q_{r})_{j-\ell_1,k-\ell_2},(Q_{0})_{j-\ell_1,k-\ell_2}, \ell_1,\ell_2=-2,...,2$ and cell averaged spatial derivatives $(\partial_xQ_{l})_{j-\ell_1,k-\ell_2}, (\partial_xQ_{r})_{j-\ell_1,k-\ell_2},(\partial_xQ_{0})_{j-\ell_1,k-\ell_2}$ can be constructed for the stencil shown in Fig.\ref{interplate}.
\item With the one-dimensional WENO reconstruction in the horizontal direction, the averaged value $(Q_{l})_{j_m,k-\ell_2}, (Q_{r})_{j_m,k-\ell_2},(Q_{0})_{j_m,k-\ell_2}$ and the averaged spatial derivatives  $(\partial_xQ_{l})_{j_m,k-\ell_2}, (\partial_xQ_{r})_{j_m,k-\ell_2},(\partial_xQ_{0})_{j_m,k-\ell_2}$ and $(\partial_yQ_{l})_{j_m,k-\ell_2}, (\partial_yQ_{r})_{j_m,k-\ell_2},(\partial_yQ_{0})_{j_m,k-\ell_2}$ over the interval $[z_{k-\ell_2}-\Delta z/2,z_{k-\ell_2}+\Delta z/2]$ with $y=y_{j_m}$ can be given.
\item With one-dimensional WENO reconstruction in the vertical direction,  the point value $(Q_{l})_{j_m,k_n}, (Q_{r})_{j_m,k_n},(Q_{0})_{j_m,k_n}$ and spatial derivatives $(\partial_xQ_{l})_{j_m,k_n}, (\partial_xQ_{r})_{j_m,k_n},(\partial_xQ_{0})_{j_m,k_n}$,  $(\partial_yQ_{l})_{j_m,k_n}, (\partial_yQ_{r})_{j_m,k_n},(\partial_yQ_{0})_{j_m,k_n}$
and $(\partial_zQ_{l})_{j_m,k_n}, (\partial_zQ_{r})_{j_m,k_n},(\partial_zQ_{0})_{j_m,k_n}$ can be fully given at the Gaussian quadrature points $(x_{i+1/2},y_{j_m},z_{k_n})$.
\end{enumerate}
In the computation, without special statement, the fifth-order WENO-JS reconstruction \cite{WENO-JS} is adopted for the flow with discontinuities and the linear scheme is used for the smooth flows to reduce the dissipation.

\section{Numerical experiments}
In this section, numerical tests for both inviscid and viscous flows will be presented to validate our numerical scheme. For the inviscid flow, the collision time $\tau$ takes
\begin{align*}
\tau=\epsilon \Delta t+C\displaystyle|\frac{p_l-p_r}{p_l+p_r}|\Delta t,
\end{align*}
where $p_l$ and $p_r$ denote the pressures on the left and right sides of the cell interface. In the computation, $\varepsilon=0.01$ and $C=1$. For the viscous flow, we have
\begin{align*}
\tau=\frac{\mu}{p}+\displaystyle|\frac{p_l-p_r}{p_l+p_r}|\Delta t,
\end{align*}
where $\mu$ is the viscous coefficient and $p$ is the pressure at the cell interface, and it will reduce to $\tau=\mu/p$ in the smooth flow regions. $\Delta t$ is the time step which is determined according to the CFL number, which takes $0.4$ in the computation. The reason for including artificial dissipation through the additional term in the particle collision time is to enlarge the kinetic scale physics in the discontinuous region for the construction of a numerical shock structure through the particle free transport and inadequate particle collision in order to keep the non-equilibrium property.

\begin{figure}[!h]
\centering
\includegraphics[width=0.475\textwidth]{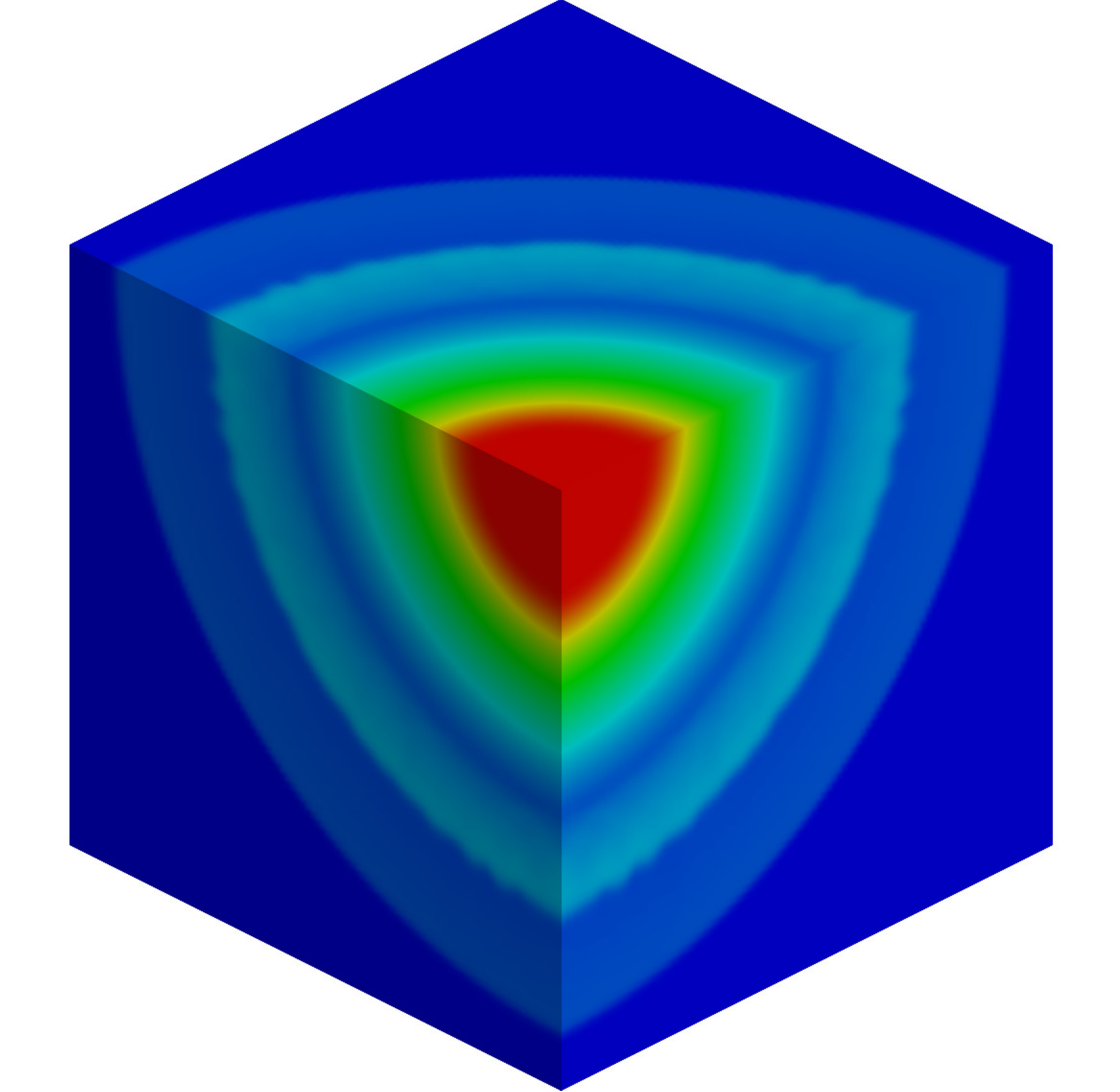}
\includegraphics[width=0.475\textwidth]{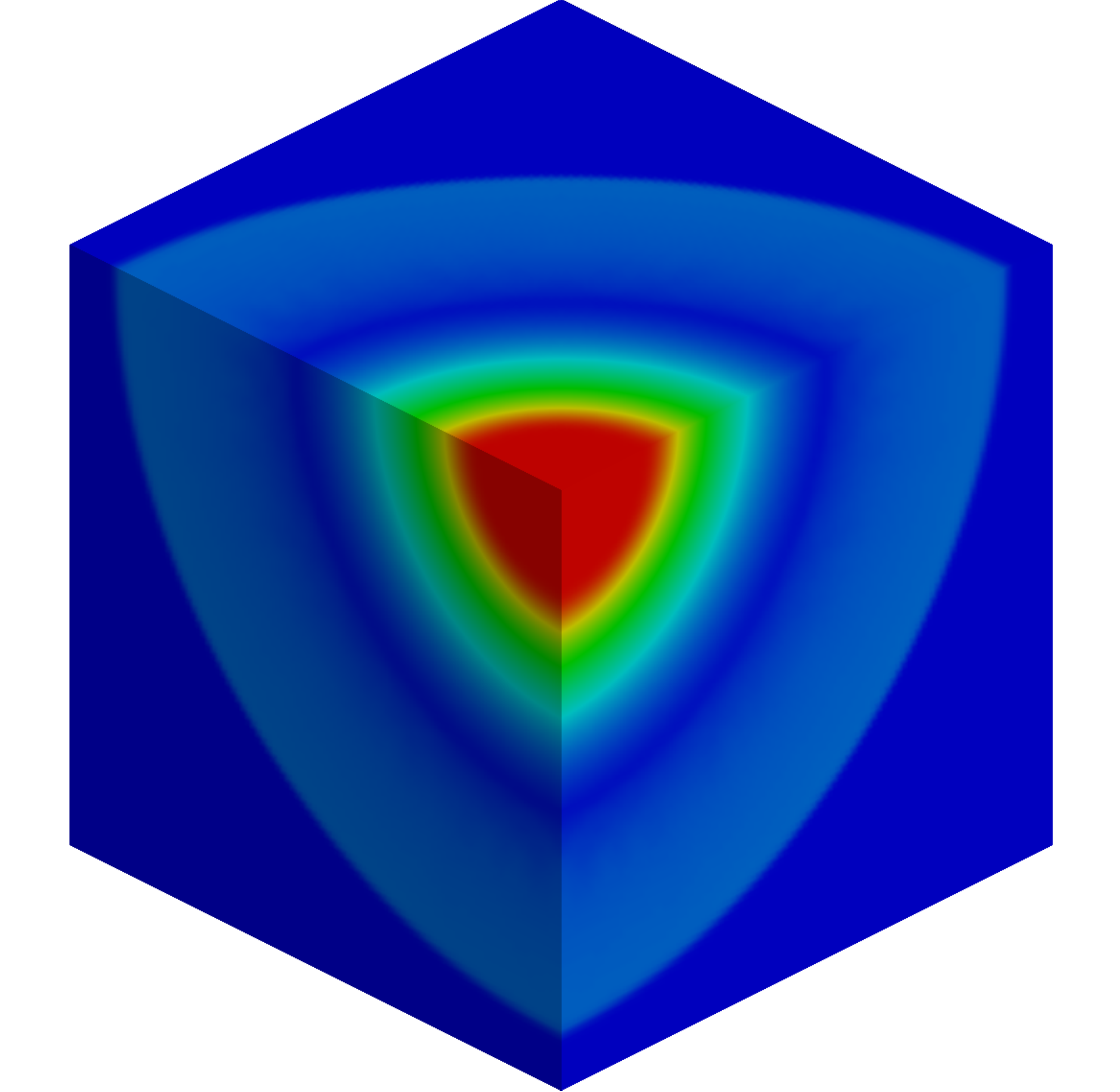}
\caption{\label{3d-riemann-1}Three-dimensional Sod problem: the three-dimensional density and pressure distributions with $\Delta x=\Delta y=\Delta z=1/100$.}
\includegraphics[width=0.475\textwidth]{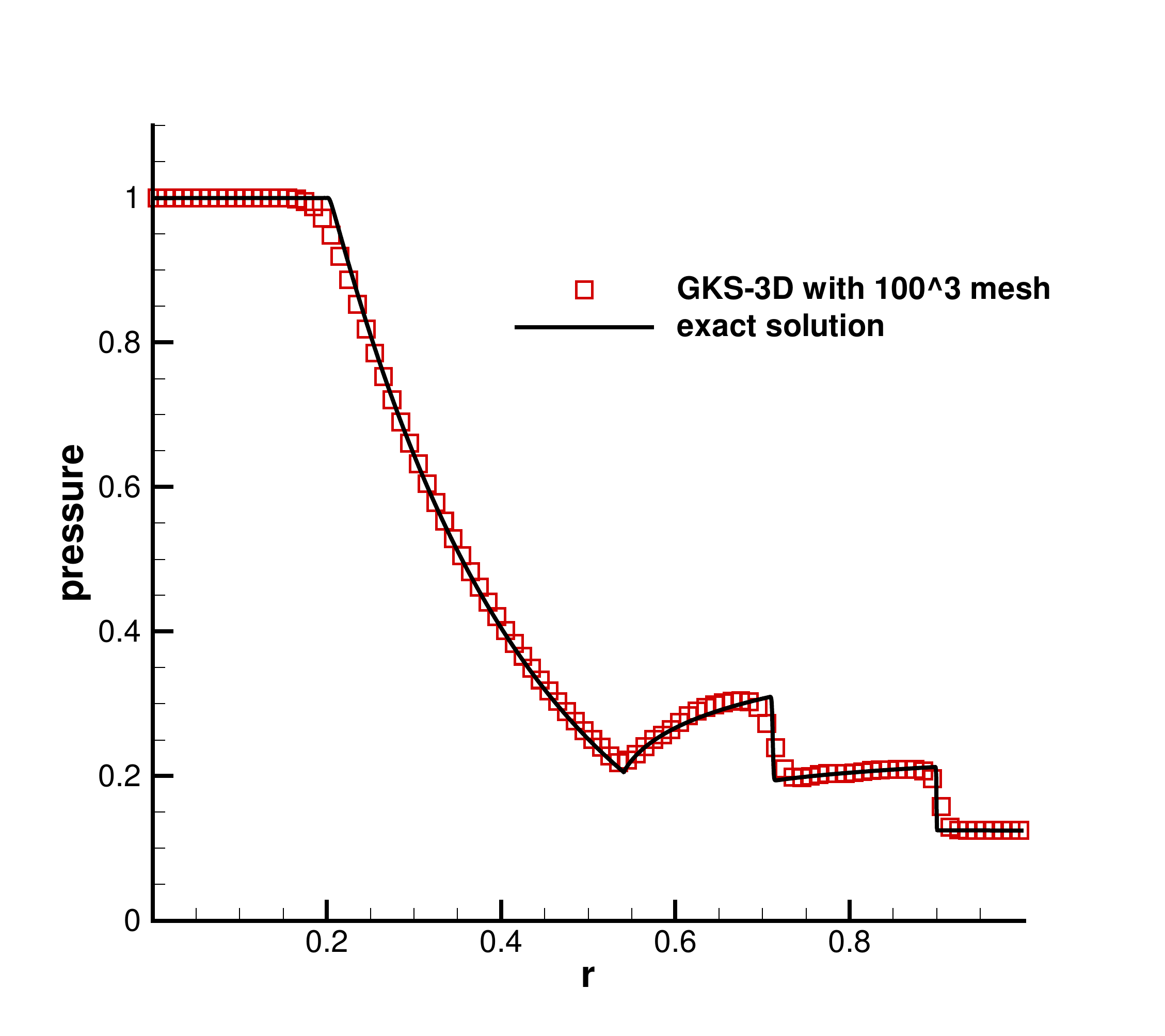}
\includegraphics[width=0.475\textwidth]{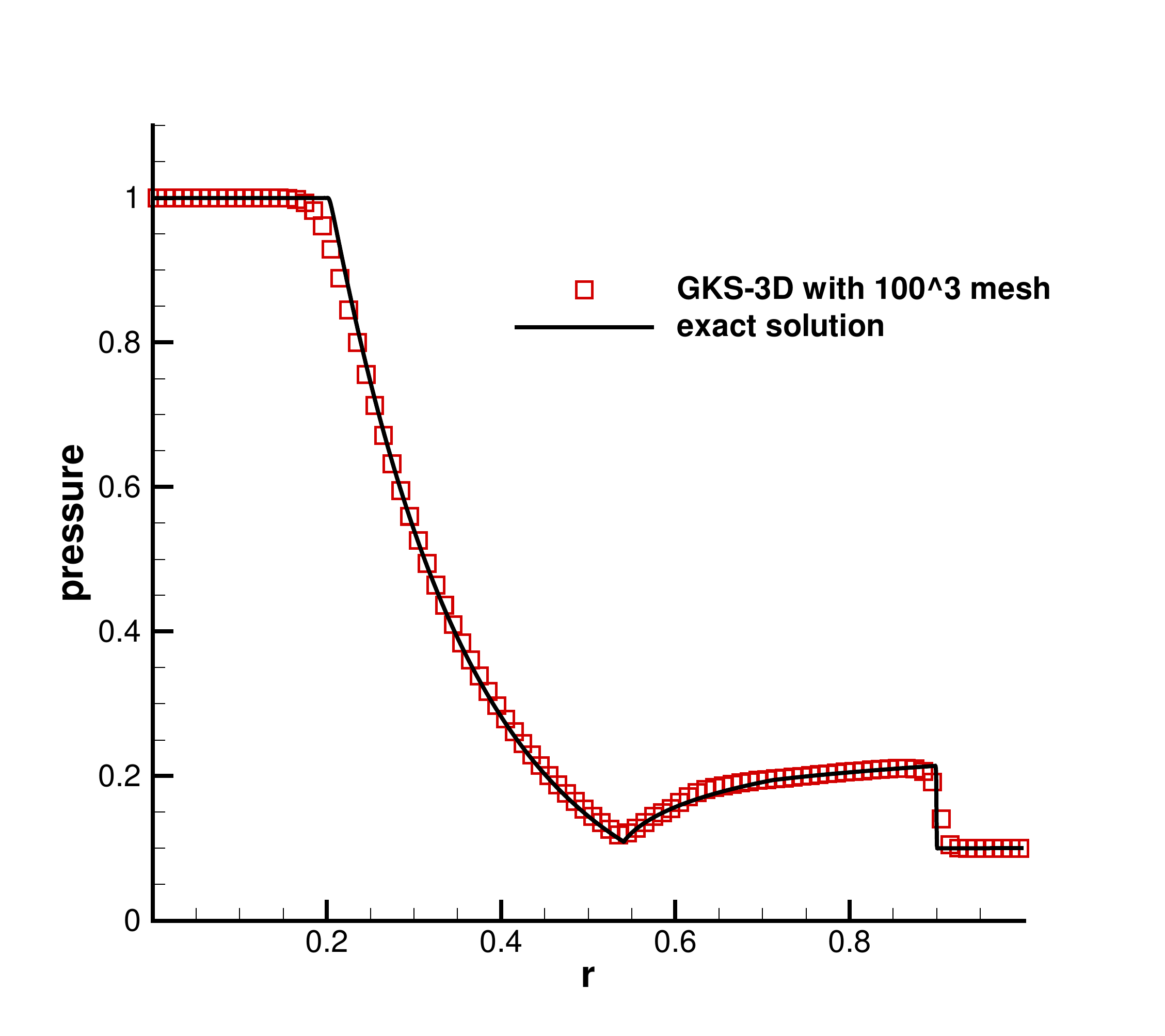}
\caption{\label{3d-riemann-2}Three-dimensional Sod problem: the density and pressure distributions with $\Delta x=\Delta y=\Delta z=1/100$ along the line $y=z=0$.}
\end{figure}

For the smooth flow, the WENO reconstruction can be used directly on the conservative flow variables. For the flow with strong discontinuity, the characteristic variables can be used in the reconstruction. Based on $A_{i+1/2,j}=(\partial F/\partial Q)_{Q=Q^*}$, where $Q$ are the conservative variables, $F(Q)$ are the corresponding fluxes, and $Q^*=(Q_{i,j,k}+Q_{i+1,j,k})/2$, the cell averaged and point conservative values can be projected into the characteristic field by $\omega=RQ$, where $R$ is the matrix corresponding to right eigenvectors of $A$. The reconstruction scheme is applied on the characteristic variables $\omega$. With the reconstructed polynomials for characteristic variables, the conservative flow variables can be recovered by the inverse projection.

In the following, we present many numerical examples. Based on the same WENO-JS reconstruction, the results presented below
have very high accuracy in comparison with  other solutions from same reconstruction.

\subsection{Three-dimensional Sod problem}
This problem is a fully three-dimensional extension of the Sod problem. The computational domain is $(x,y,z)\in[0, 1]\times[0, 1]\times[0, 1]$, and a sphere of radius $R =0.5$ separates two different constant states
\begin{equation*}
(\rho, U, V, p)=
\begin{cases}
(1, 0, 0, 1),  \ \ \ \ &  0<R<0.5,\\
(0.125, 0, 0, 0.1),    & 0.5<R<1.
\end{cases}
\end{equation*}
The exact solution of this problem can be given by the one-dimensional system with geometric source terms
\begin{align*}
\frac{\partial Q}{\partial t}+\frac{\partial F(Q)}{\partial
r}=S(Q),
\end{align*}
with
\begin{align*}
Q=\begin{pmatrix}
   \rho \\
   \rho U \\
   \rho E \\
\end{pmatrix},
F(Q)=\begin{pmatrix}
   \rho U\\
   \rho U^2+p \\
   U(\rho E+p) \\
 \end{pmatrix},
S(Q)=-\frac{d-1}{r}\begin{pmatrix}
   \rho U\\
   \rho U^2  \\
   U(\rho E+p) \\
 \end{pmatrix},
\end{align*}
where the radial direction is denoted by $r$, while $U$ is the radial velocity, $d$ is the number of space dimensions and $\gamma=1.4$. The two-stage fourth-order gas-kinetic scheme is used to solve this equation and the numerical results with $10000$ cells are given as the reference solutions. In the computation, the uniform mesh with $\Delta x=\Delta y=\Delta z=1/100$ is used. The symmetric boundary condition is imposed on the plane with $x=0$, $y=0$ and $z=0$, and the non-reflection boundary condition is imposed on the plane with $x=1$, $y=1$ and $z=1$. The three-dimensional density and pressure distributions are given in Fig.\ref{3d-riemann-1}, and the density and pressure distribution along the line $y=z=0$ are given in Fig.\ref{3d-riemann-2}. The numerical results agree well with the reference solutions.

\begin{figure}[!h]
\centering
\includegraphics[height=0.3\textwidth]{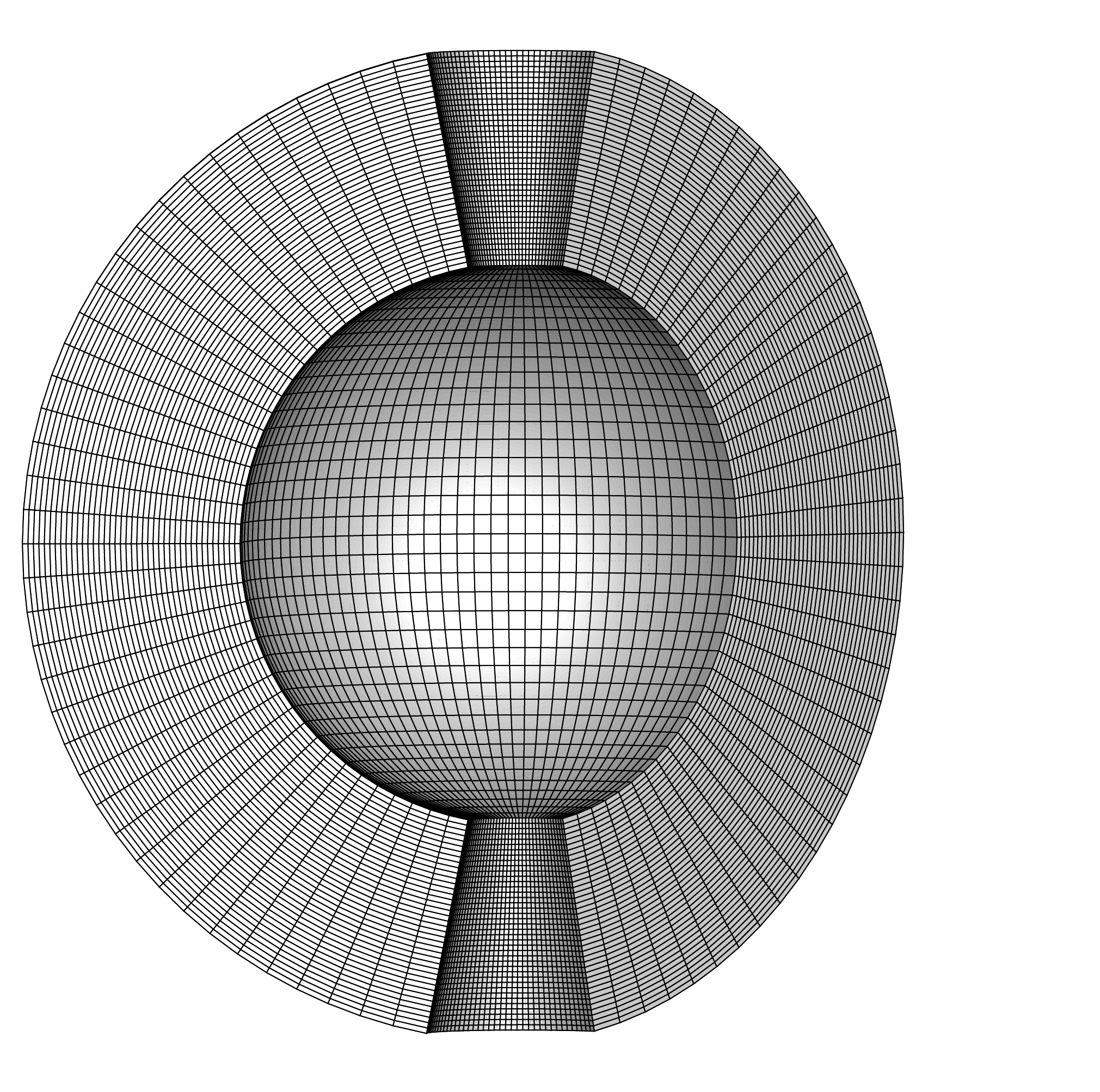}{a}
\includegraphics[height=0.3\textwidth]{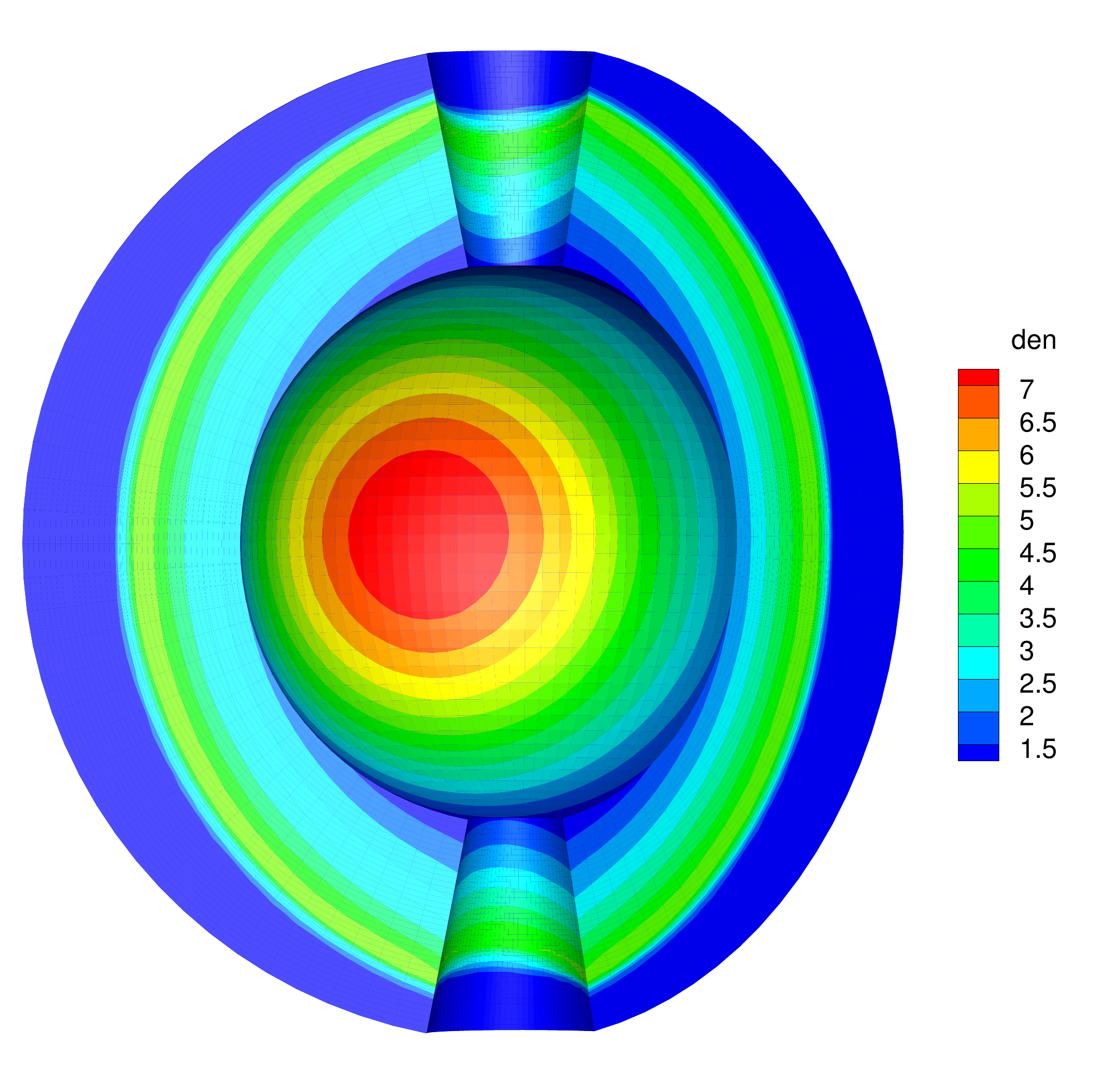}{b}
\includegraphics[height=0.3\textwidth]{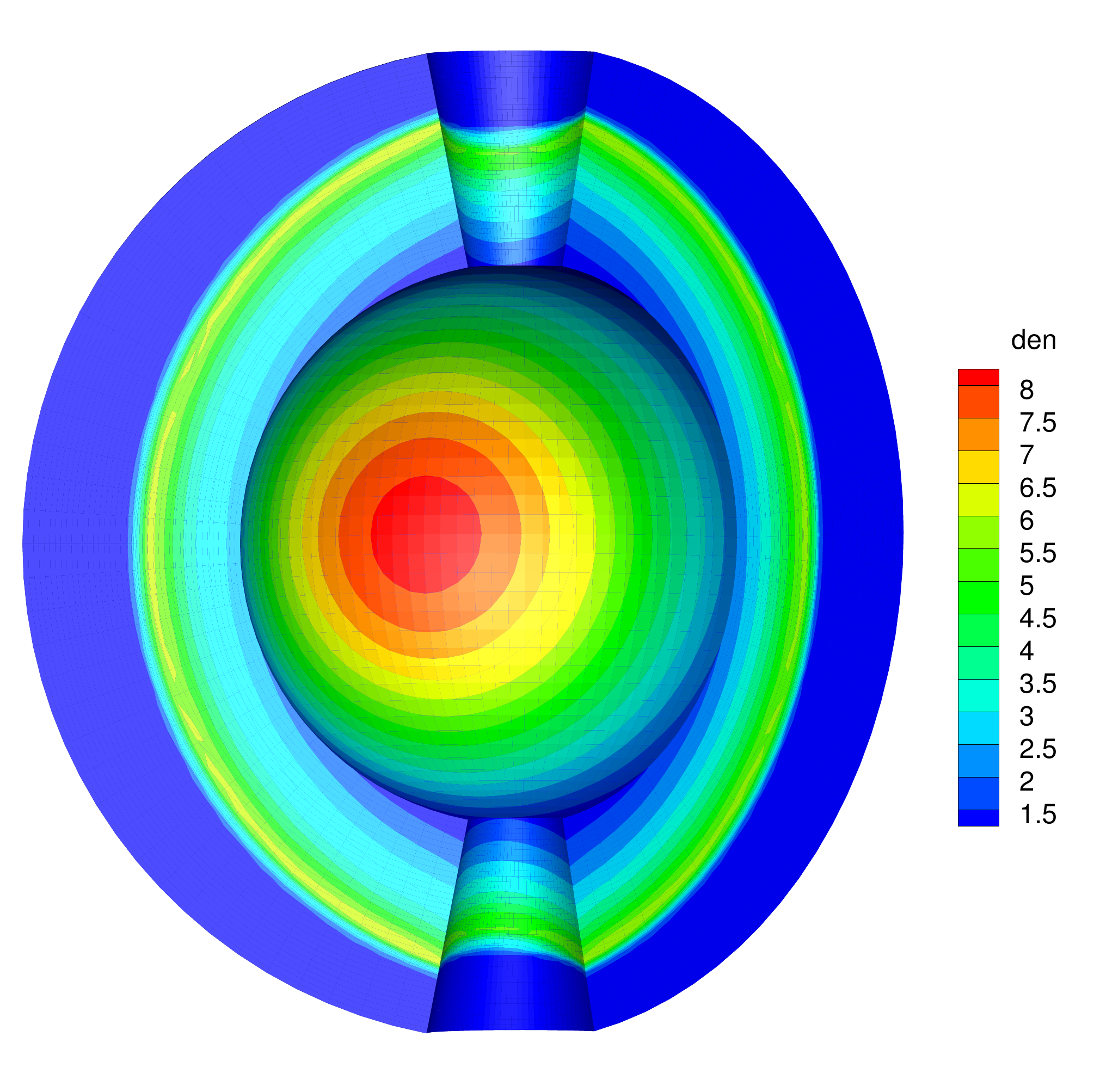}{c}
\caption{\label{sphere-0}Flow impinging on sphere: the computational mesh (a), density distribution with Mach number $Ma=5$ (b) and $8$ (c).}
\end{figure}

\begin{figure}[!h]
\centering
\includegraphics[height=0.375\textwidth]{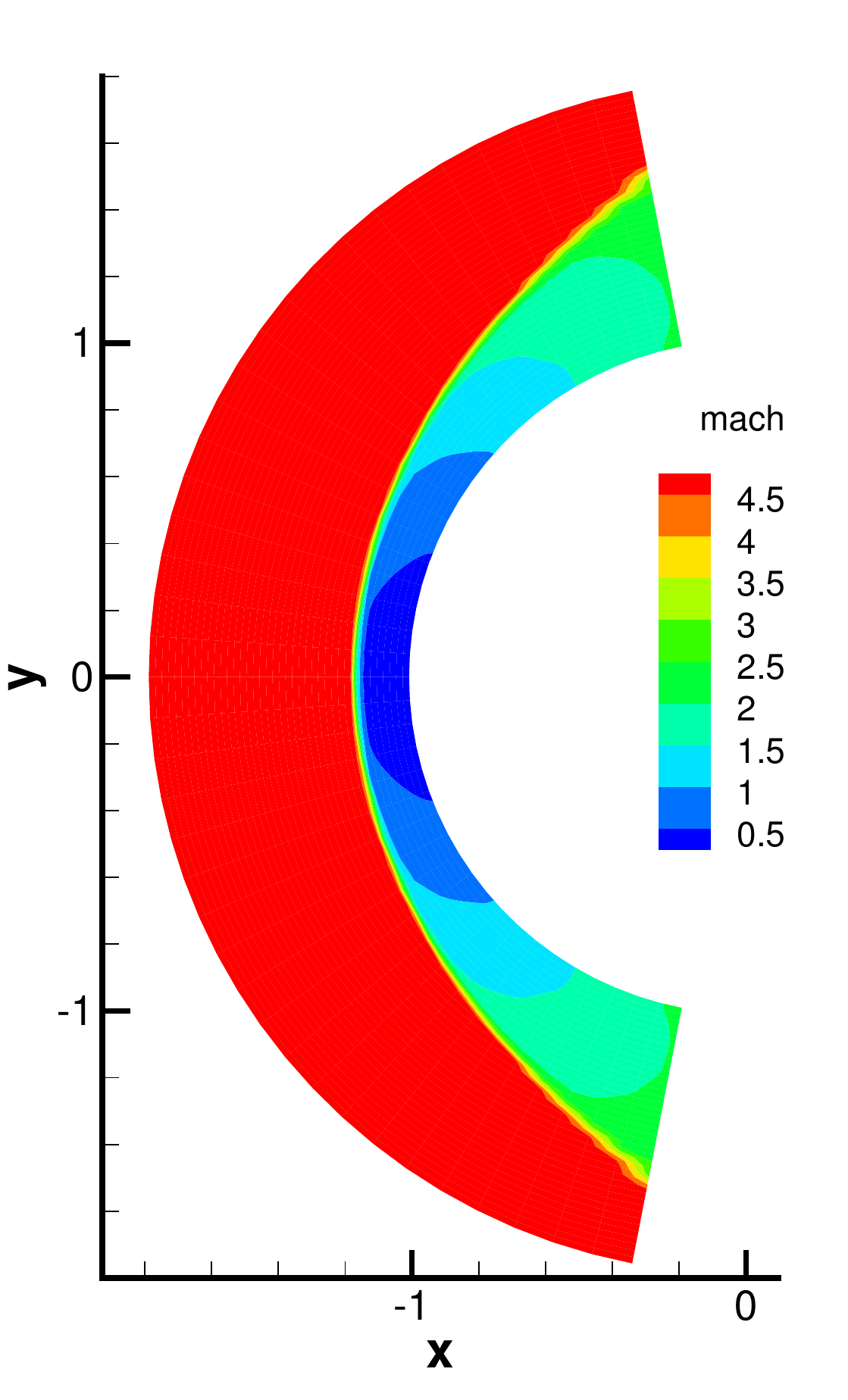}{a}
\includegraphics[height=0.375\textwidth]{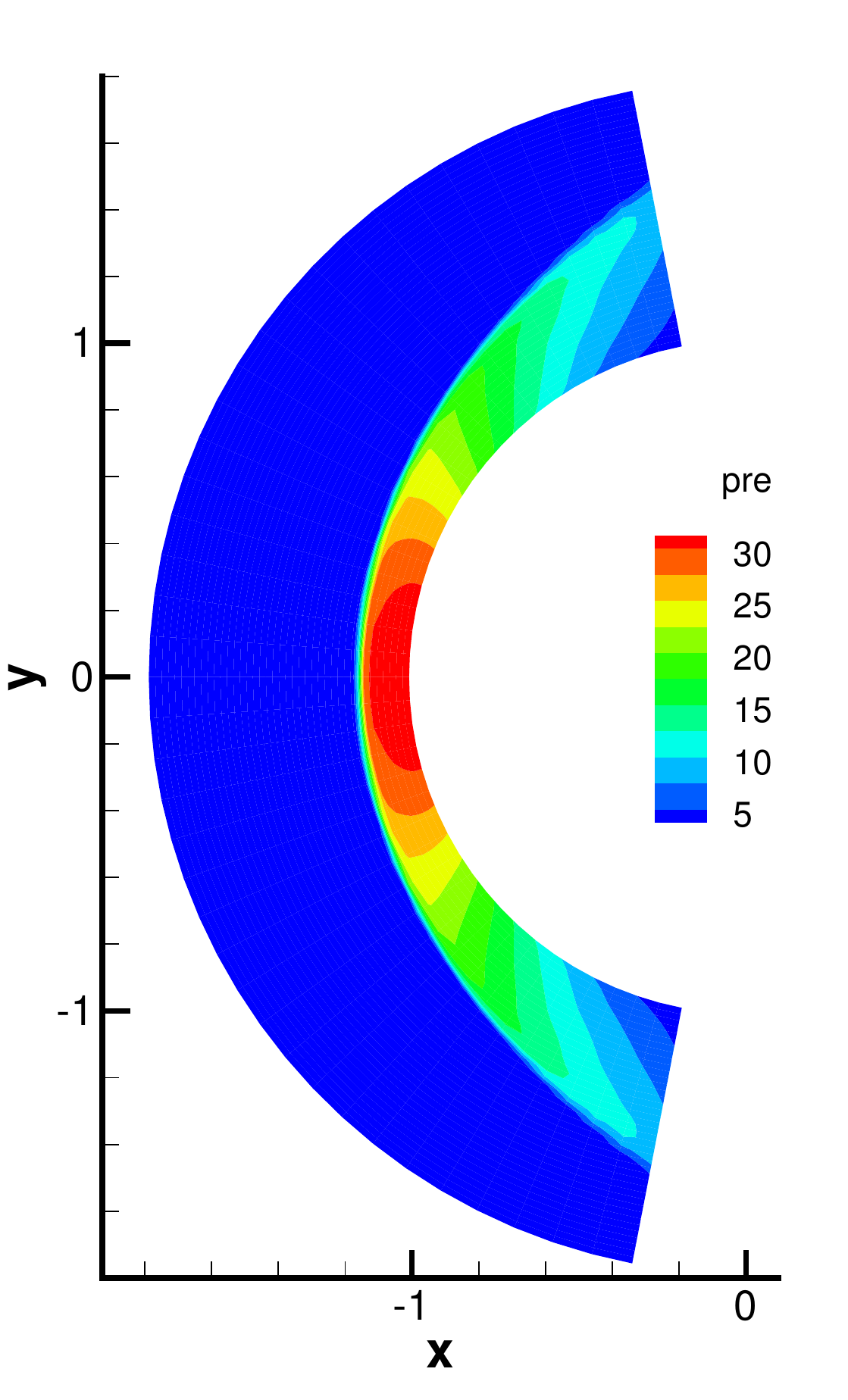}{b}
\includegraphics[height=0.375\textwidth]{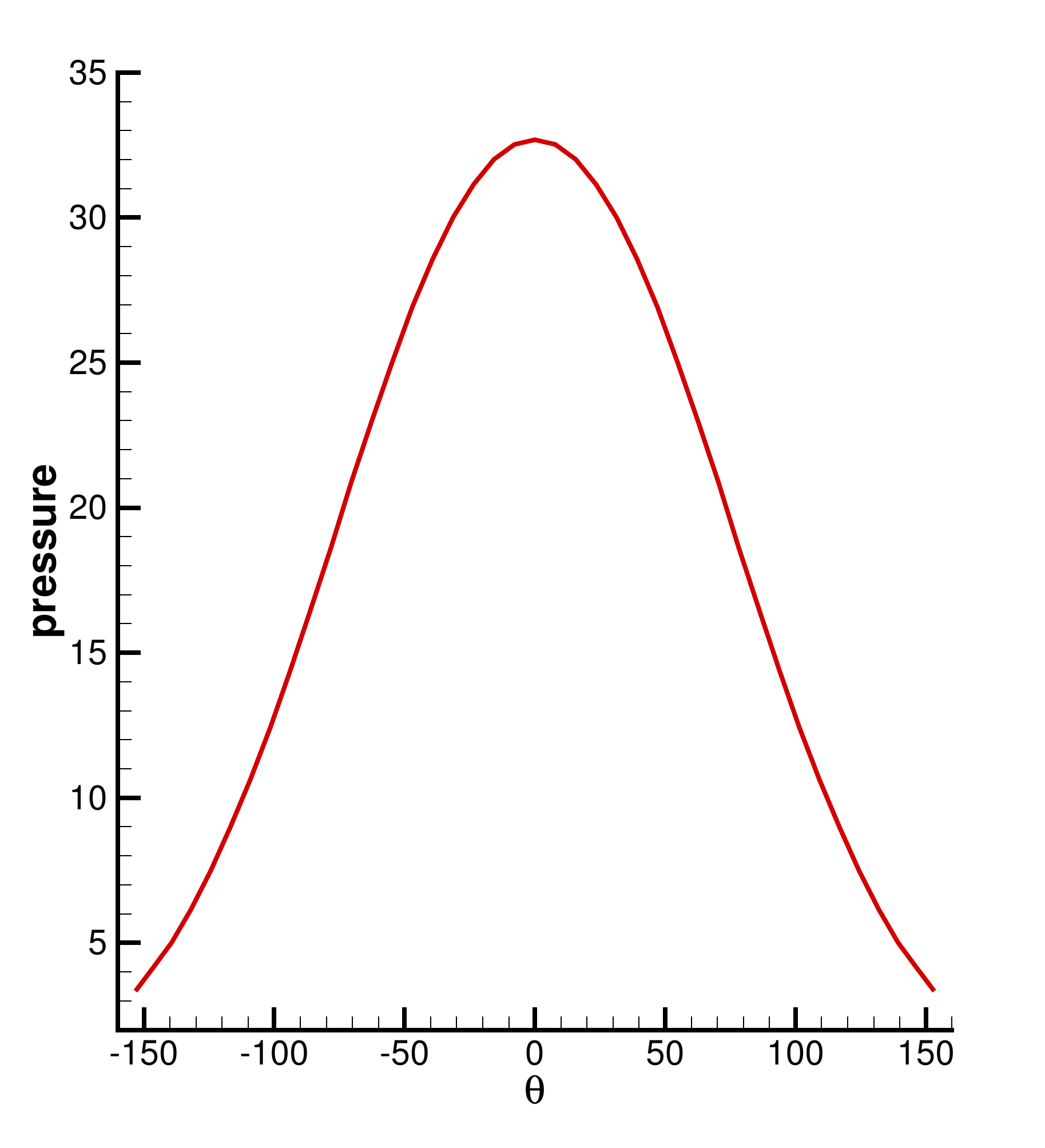}{c}
\caption{\label{sphere-1}Flow impinging on sphere: the Mach number (a), pressure (b) distributions at the plane with $\phi=0$, and pressure profile along the surface of sphere with $\phi=0$ (c) for Mach number $Ma=5$.}
\centering
\includegraphics[height=0.375\textwidth]{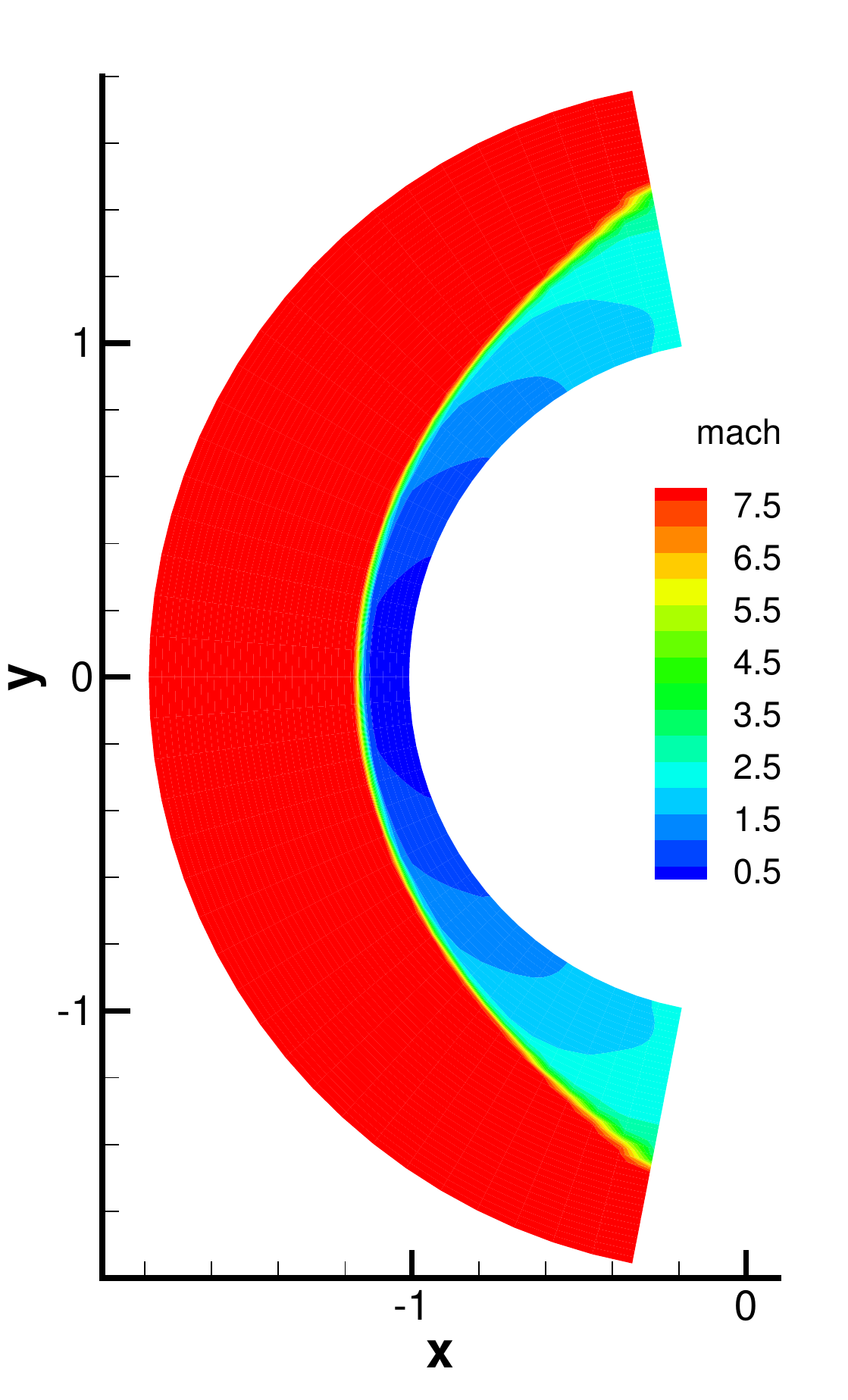}{a}
\includegraphics[height=0.375\textwidth]{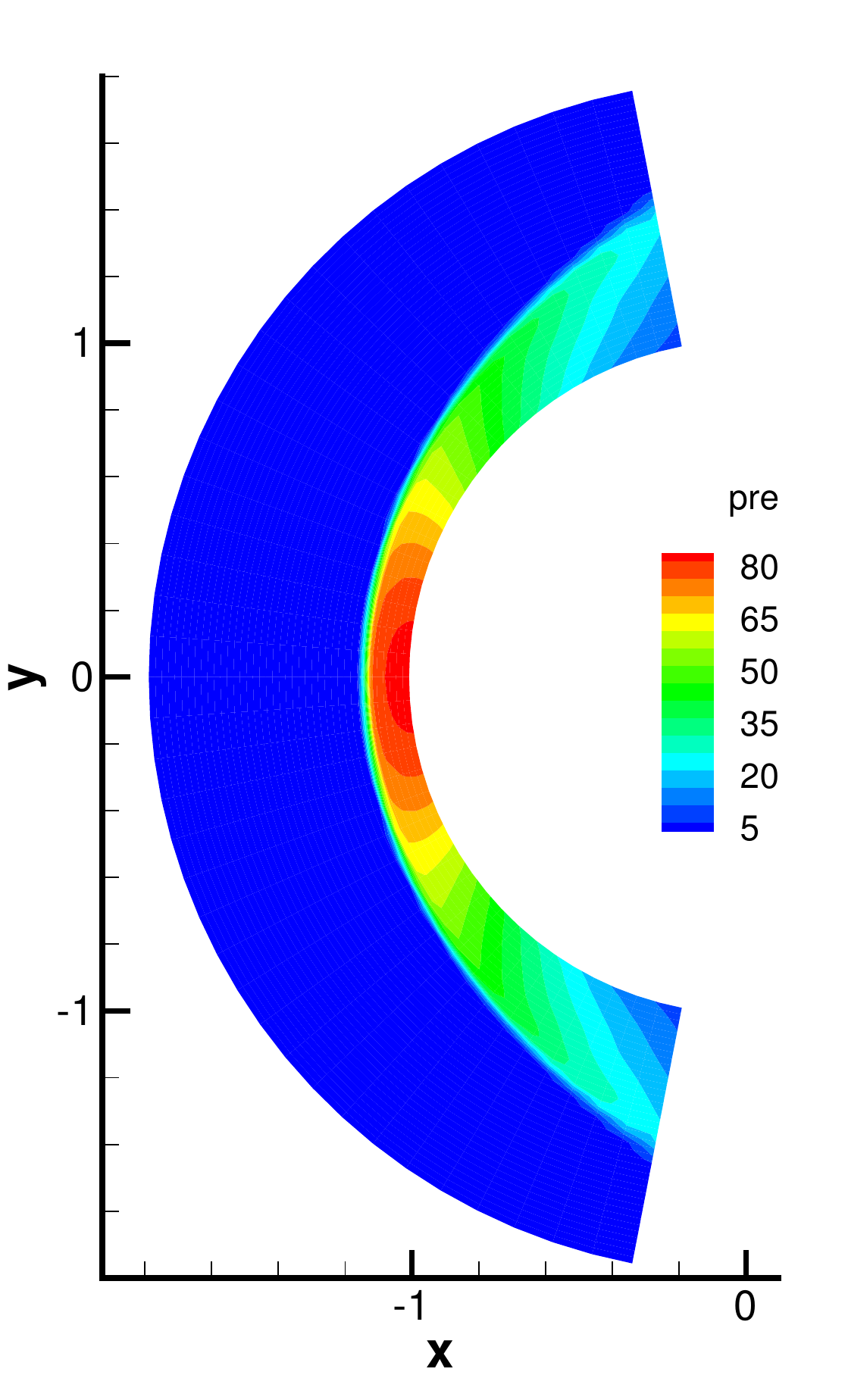}{b}
\includegraphics[height=0.375\textwidth]{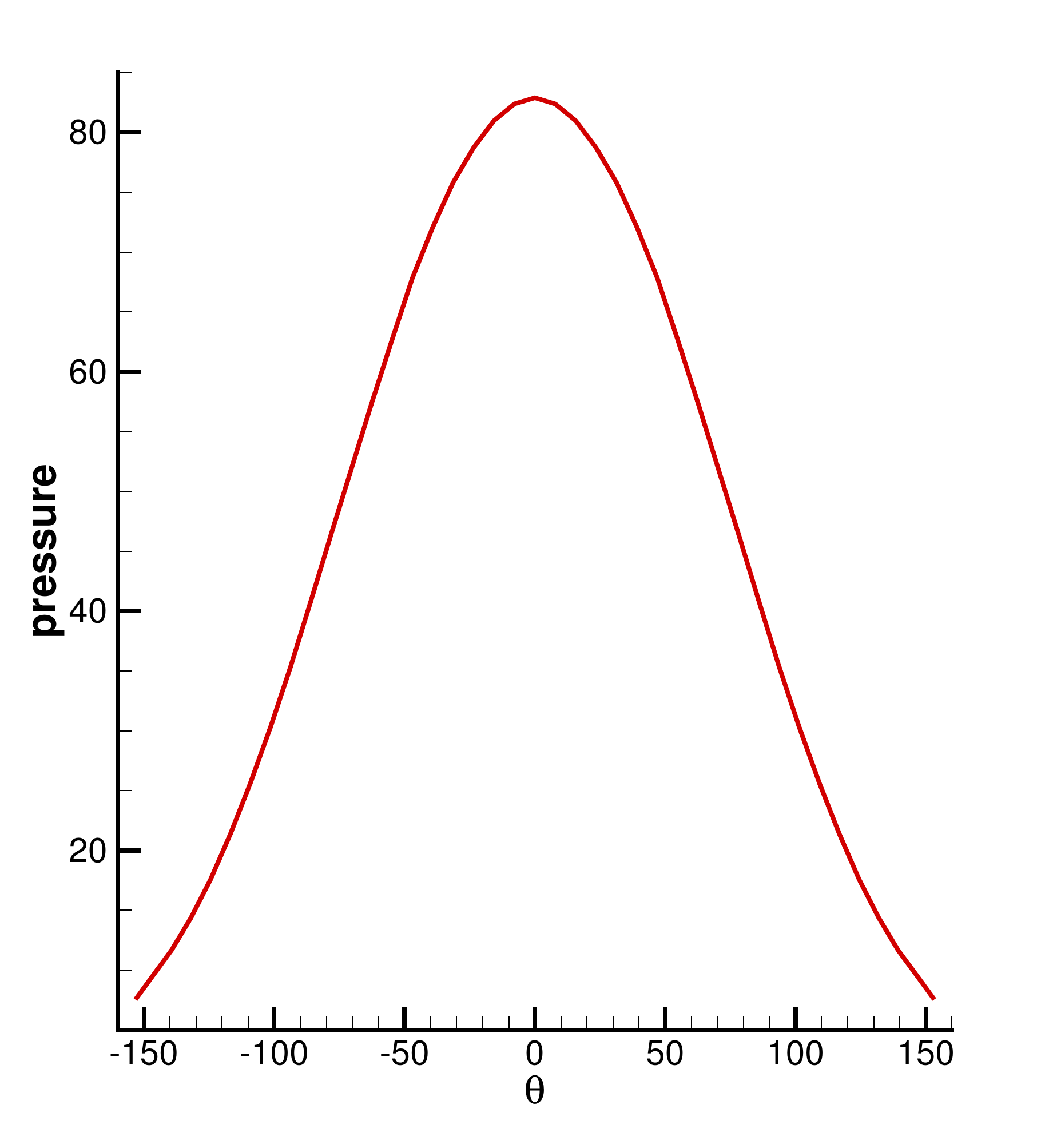}{c}
\caption{\label{sphere-2}Flow impinging on sphere: the Mach number (a), pressure (b) distributions at the plane with $\phi=0$, and pressure profile along the surface of sphere with $\phi=0$ (c) for Mach number $Ma=8$.}
\end{figure}

\subsection{Flow impinging on sphere}
\begin{figure}[!h]
\centering
\includegraphics[height=0.4\textwidth]{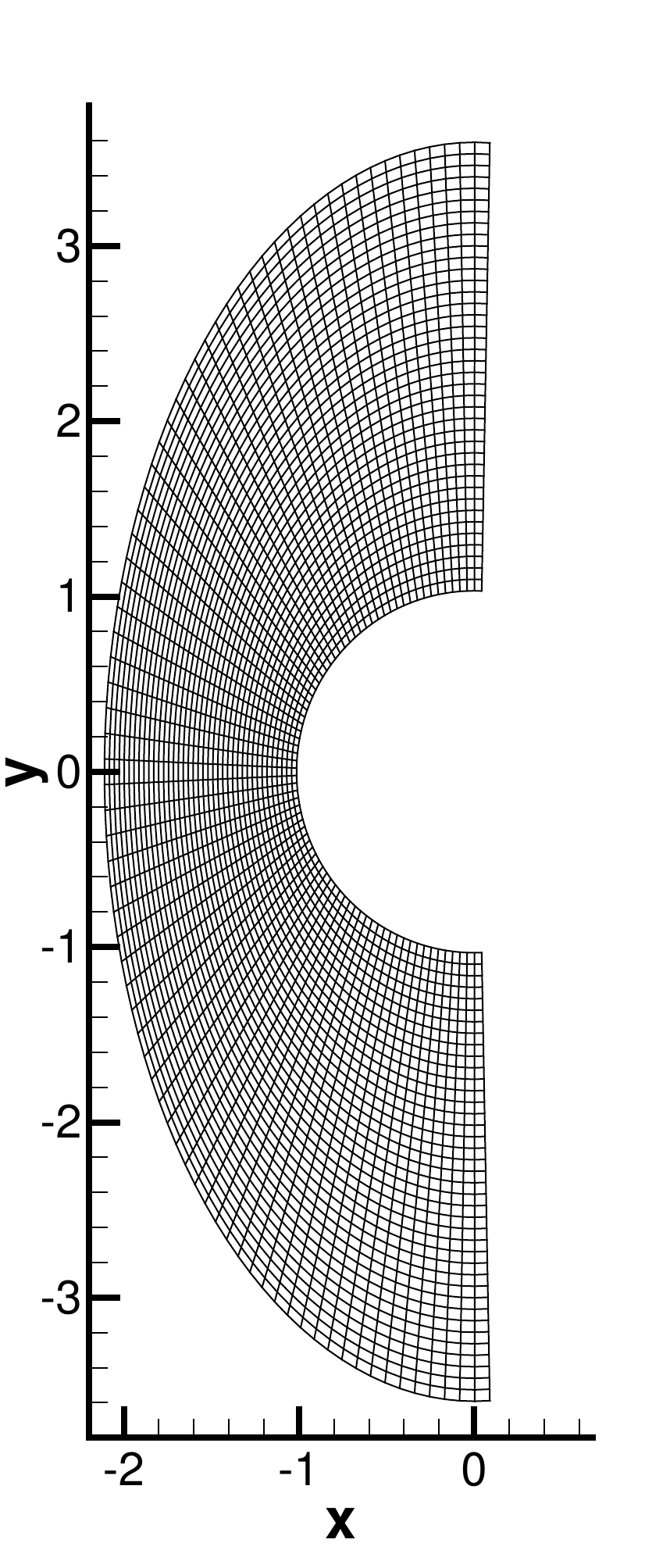}{a}
\includegraphics[height=0.4\textwidth]{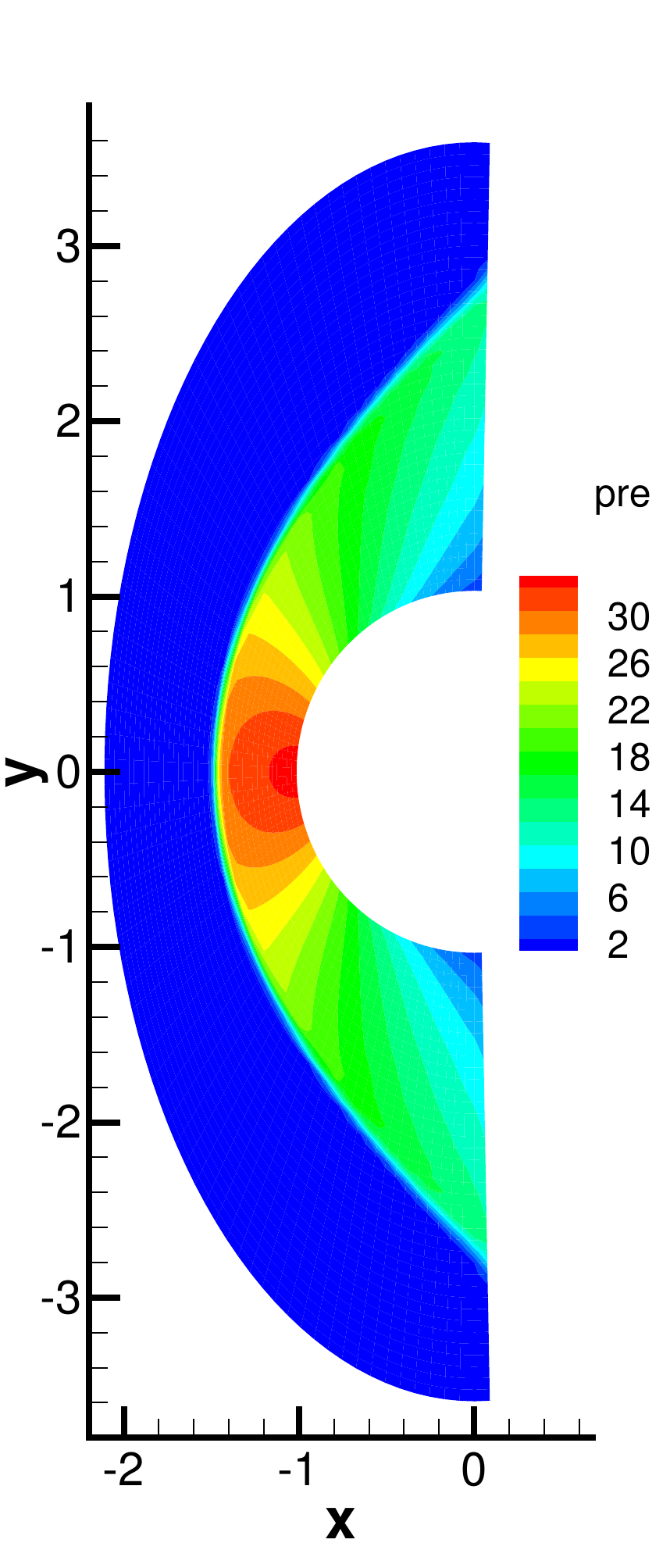}{b}
\includegraphics[height=0.4\textwidth]{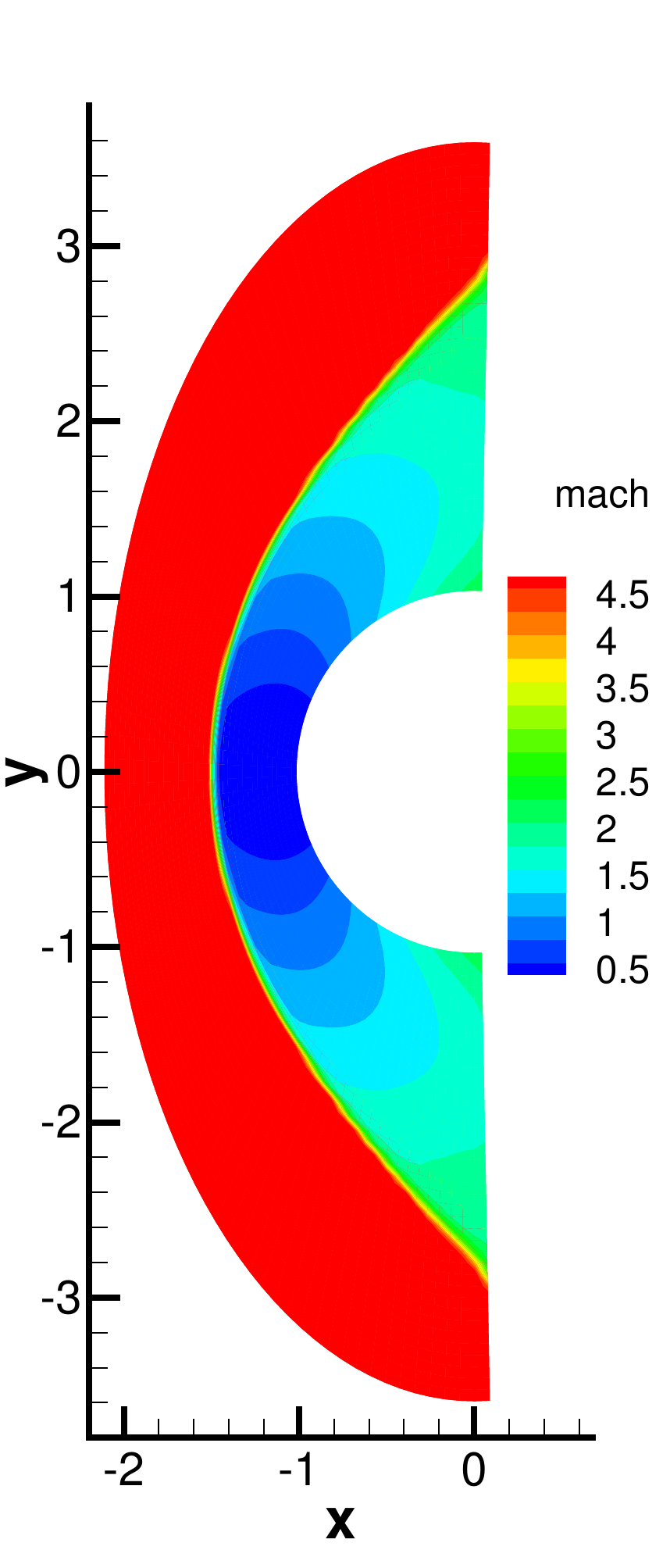}{c}
\includegraphics[height=0.4\textwidth]{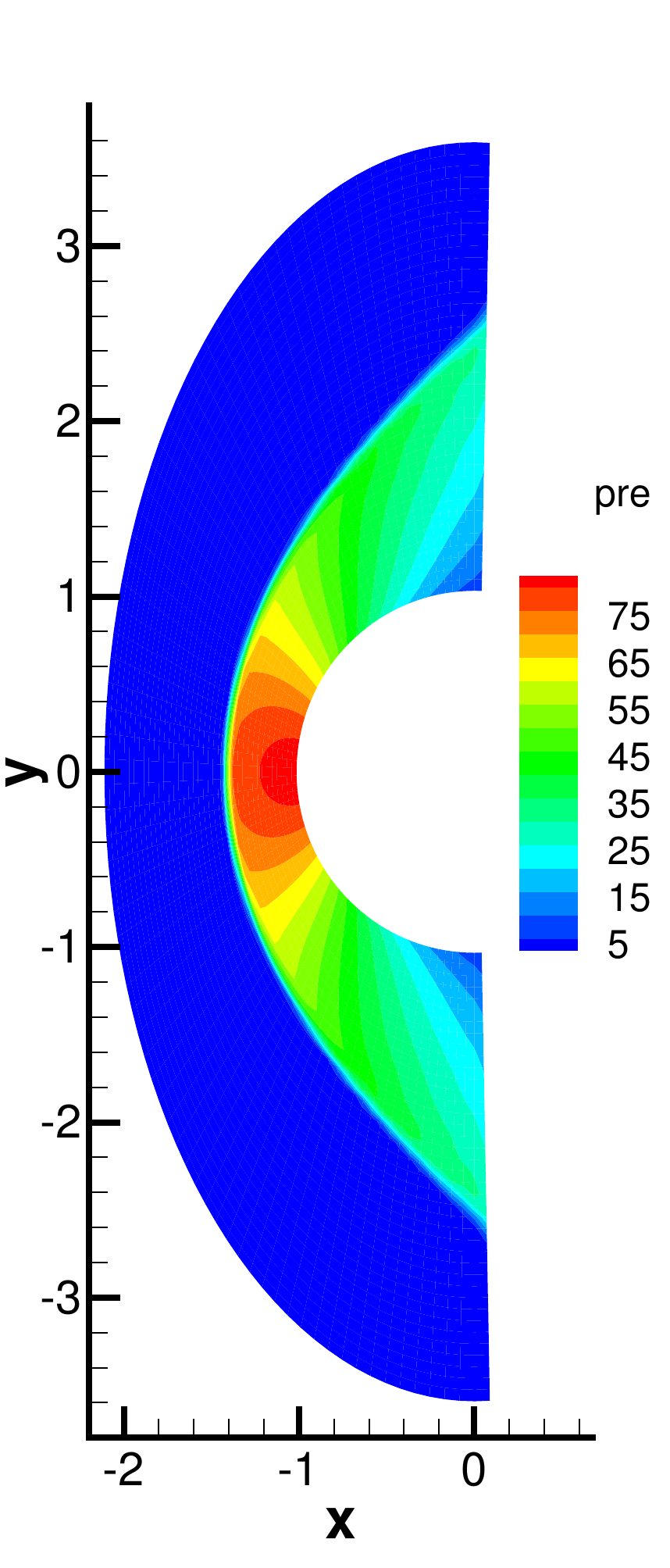}{d}
\includegraphics[height=0.4\textwidth]{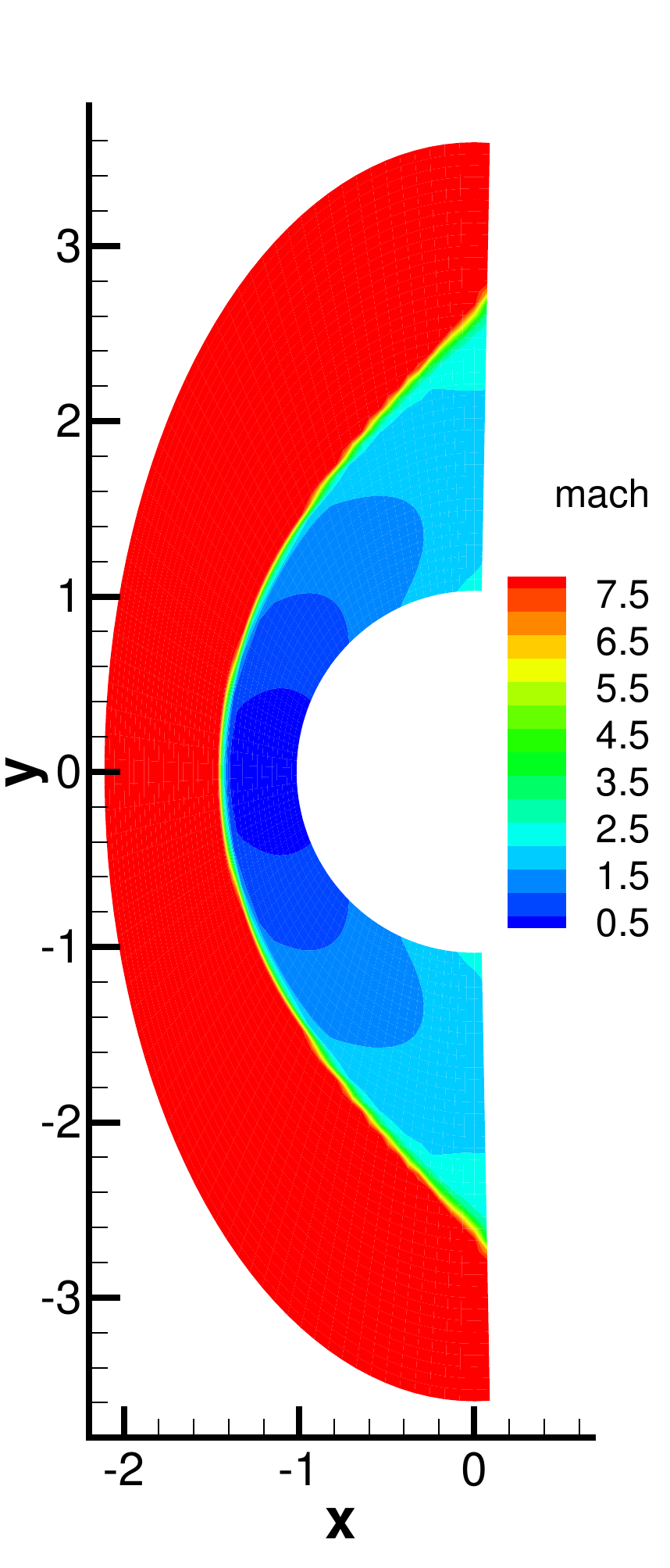}{e}
\caption{\label{sphere-3}Flow impinging on a cylinder: the computational mesh (a),  Mach number (b), pressure (c) distributions for Mach number $Ma=5$, and  Mach number (d), pressure (e) distributions for $Ma=8$.}
\end{figure}

In this case, the inviscid hypersonic flows impinging on a unit sphere are tested to validate robustness of the current scheme with different Mach numbers with $\gamma=1.4$. In the computation, a $40\times40\times40$ mesh shown in Fig.\ref{sphere-0} is used, which representing the domain $[-1.75,-1]\times[-0.4\pi,0.4\pi]\times[0.1\pi,0.9\pi]$ in the spherical coordinate $(r,\phi,\theta)$.  The density distributions in the whole domain, Mach number and pressure distributions at the plane with $\phi=0$, and pressure profiles along the surface of sphere with $\phi=0$ for Mach number $Ma=5$ and $8$ are shown in Fig.\ref{sphere-0}, Fig.\ref{sphere-1} and Fig.\ref{sphere-2}, where the shock is well captured by the current scheme and the carbuncle phenomenon does not appear \cite{Case-Pandolfi}. As comparison, the inviscid hypersonic flows impinging on a two-dimensional unit cylinder is also tested with $40\times80$ mesh. The Mach number and  pressure distributions for Mach number $Ma=5, 8$ are shown in Fig.\ref{sphere-3}.

\begin{figure}[!h]
\centering
\includegraphics[width=0.14\textwidth]{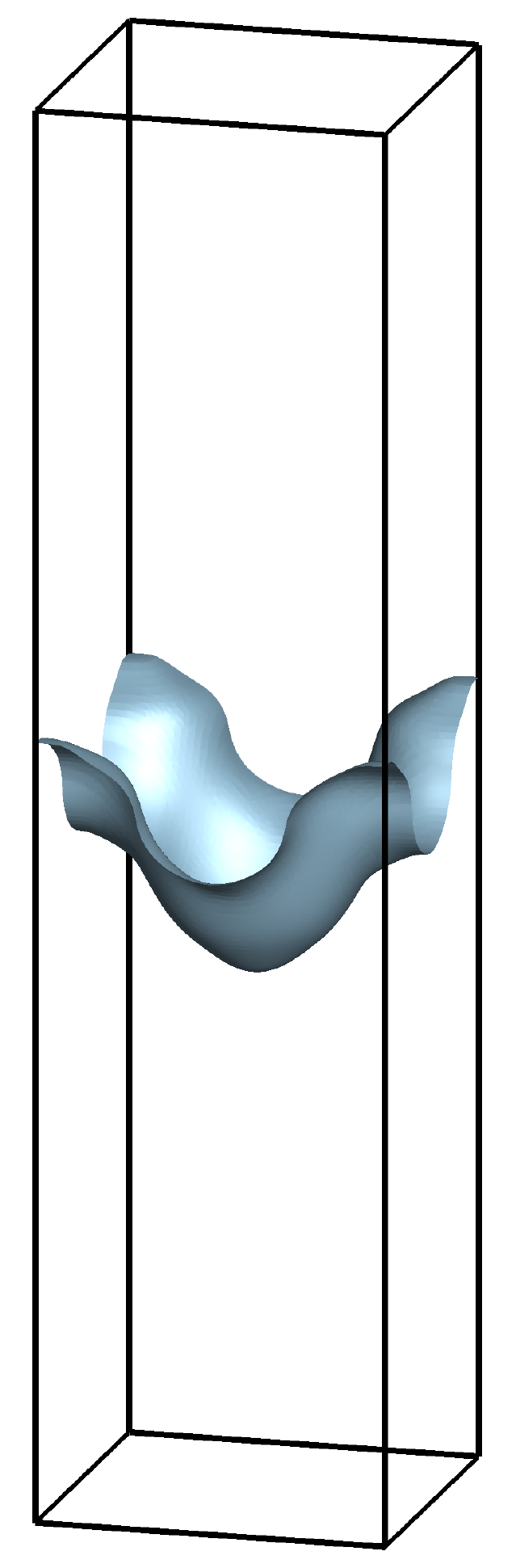}
\includegraphics[width=0.14\textwidth]{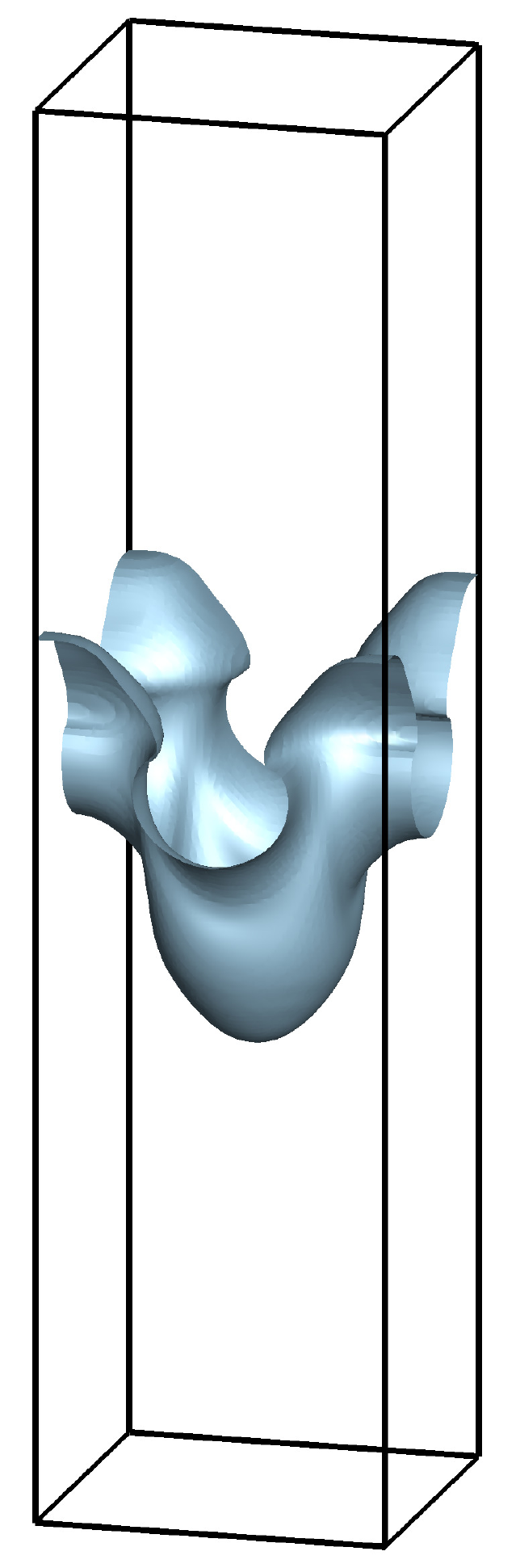}
\includegraphics[width=0.14\textwidth]{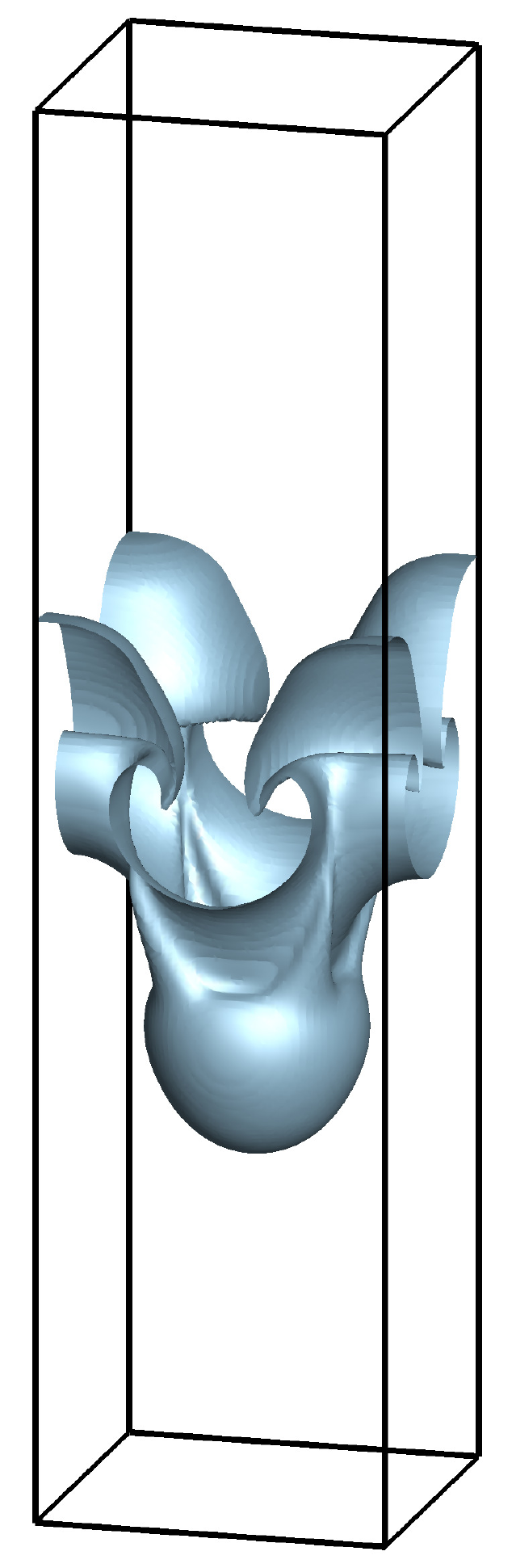}
\includegraphics[width=0.14\textwidth]{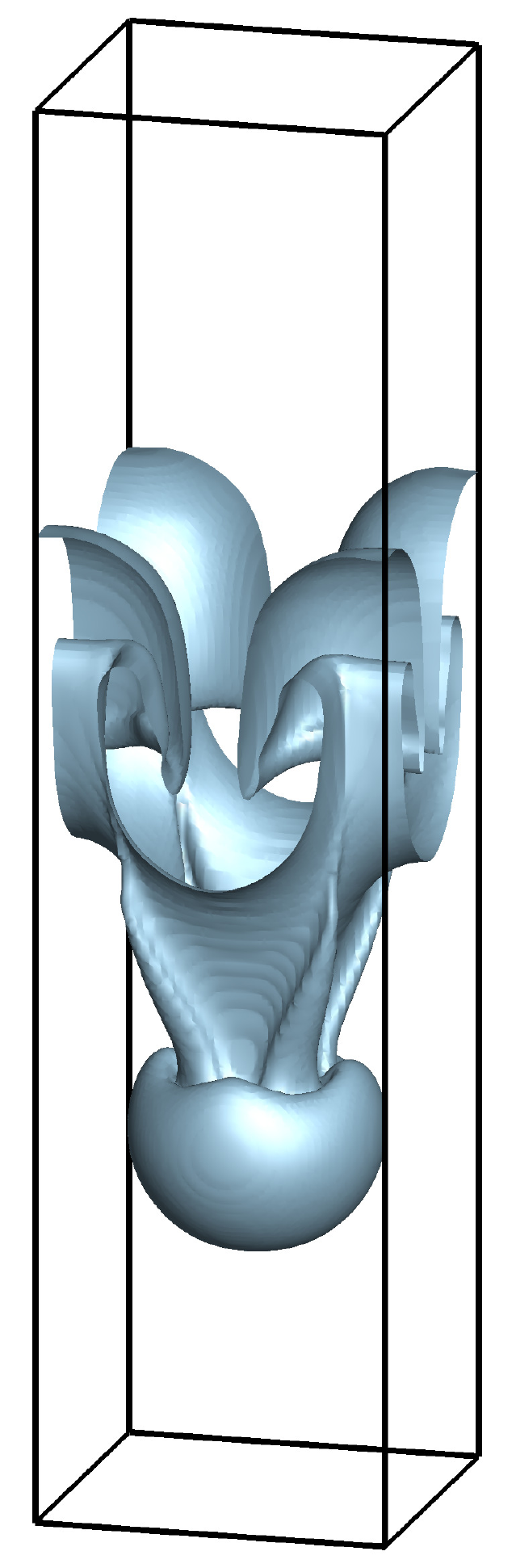}
\includegraphics[width=0.14\textwidth]{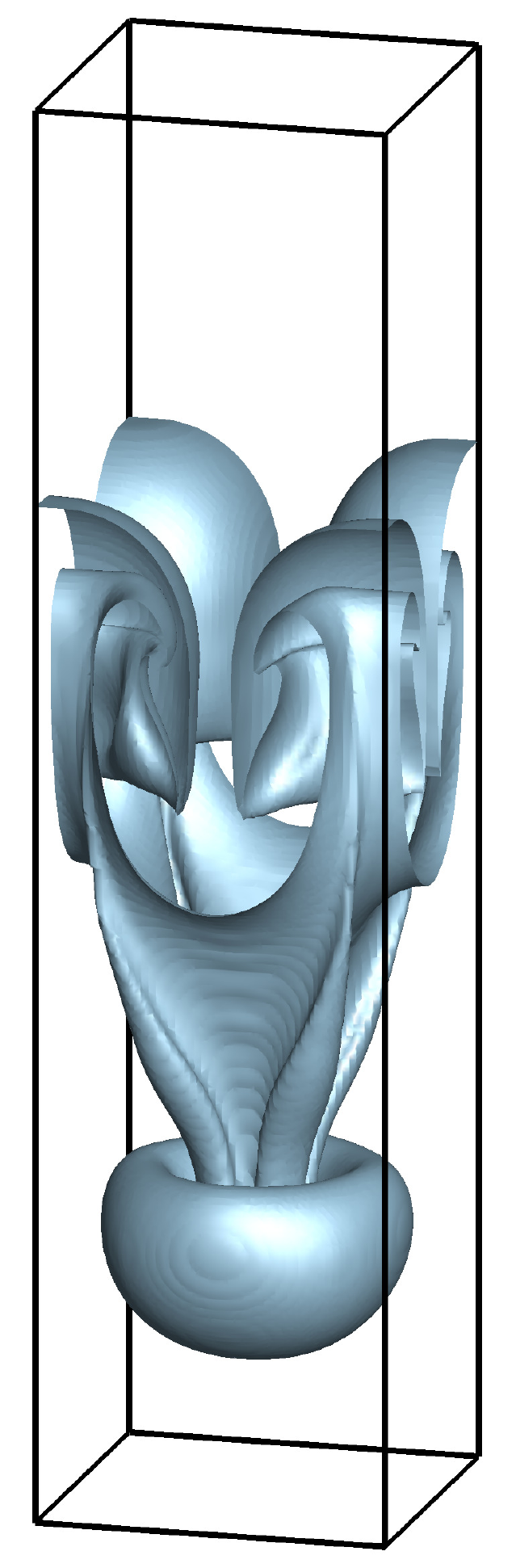}
\includegraphics[width=0.14\textwidth]{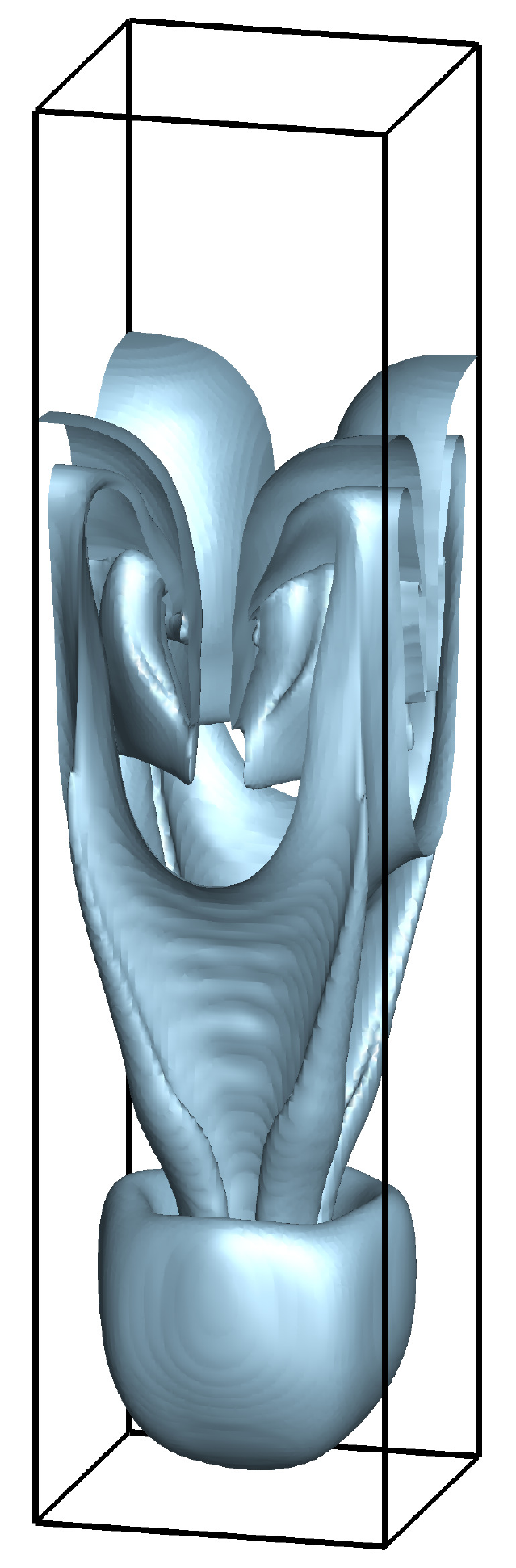}
\caption{\label{rt-interface}Rayleigh-Taylor instability: the fluid interface from a single-mode perturbation at $t=3, 4, 5, 6, 7, 8$.}
\end{figure}

\subsection{Rayleigh-Taylor instability}
\begin{figure}[!h]
\centering
\includegraphics[width=0.3\textwidth]{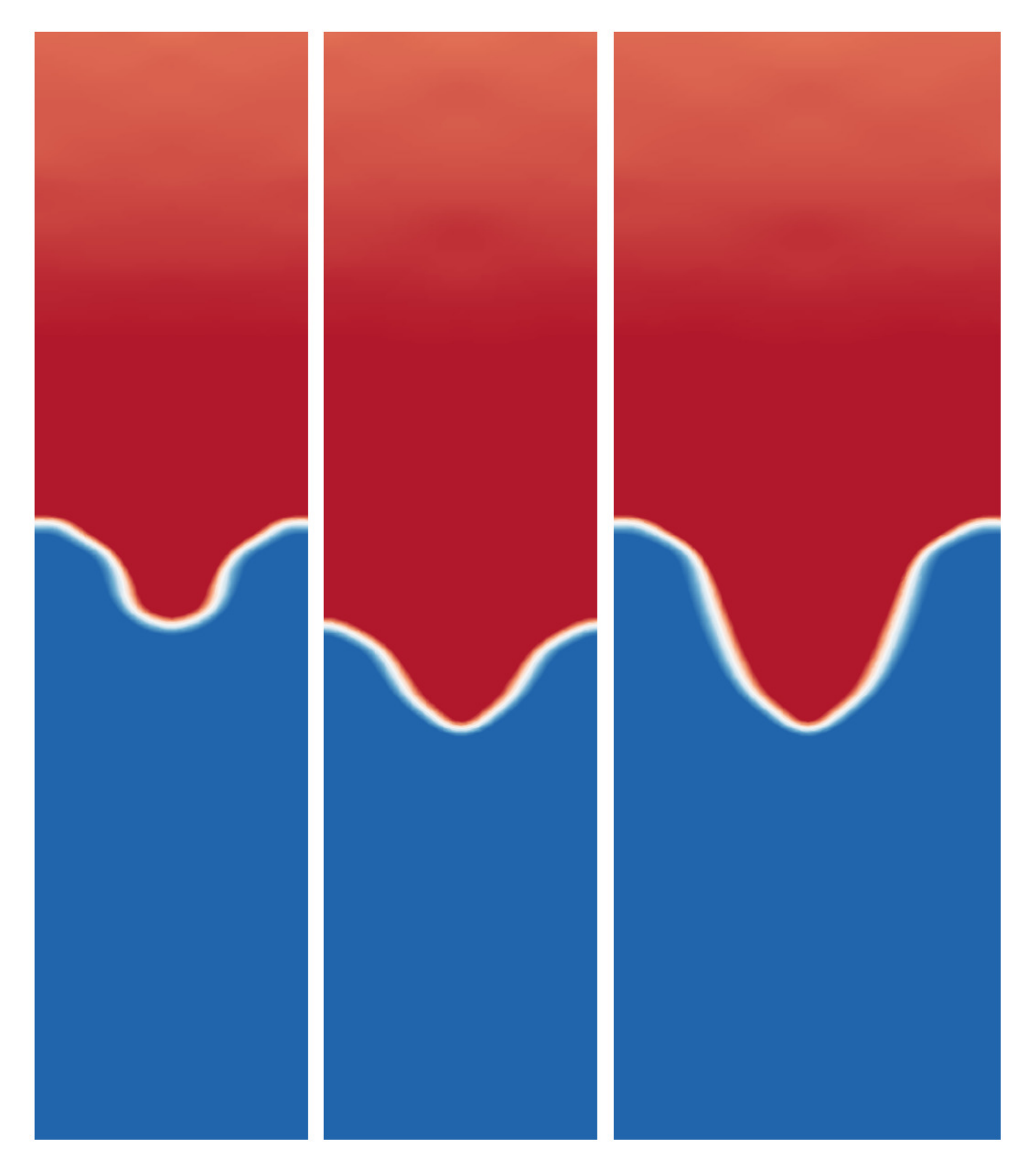}
\includegraphics[width=0.3\textwidth]{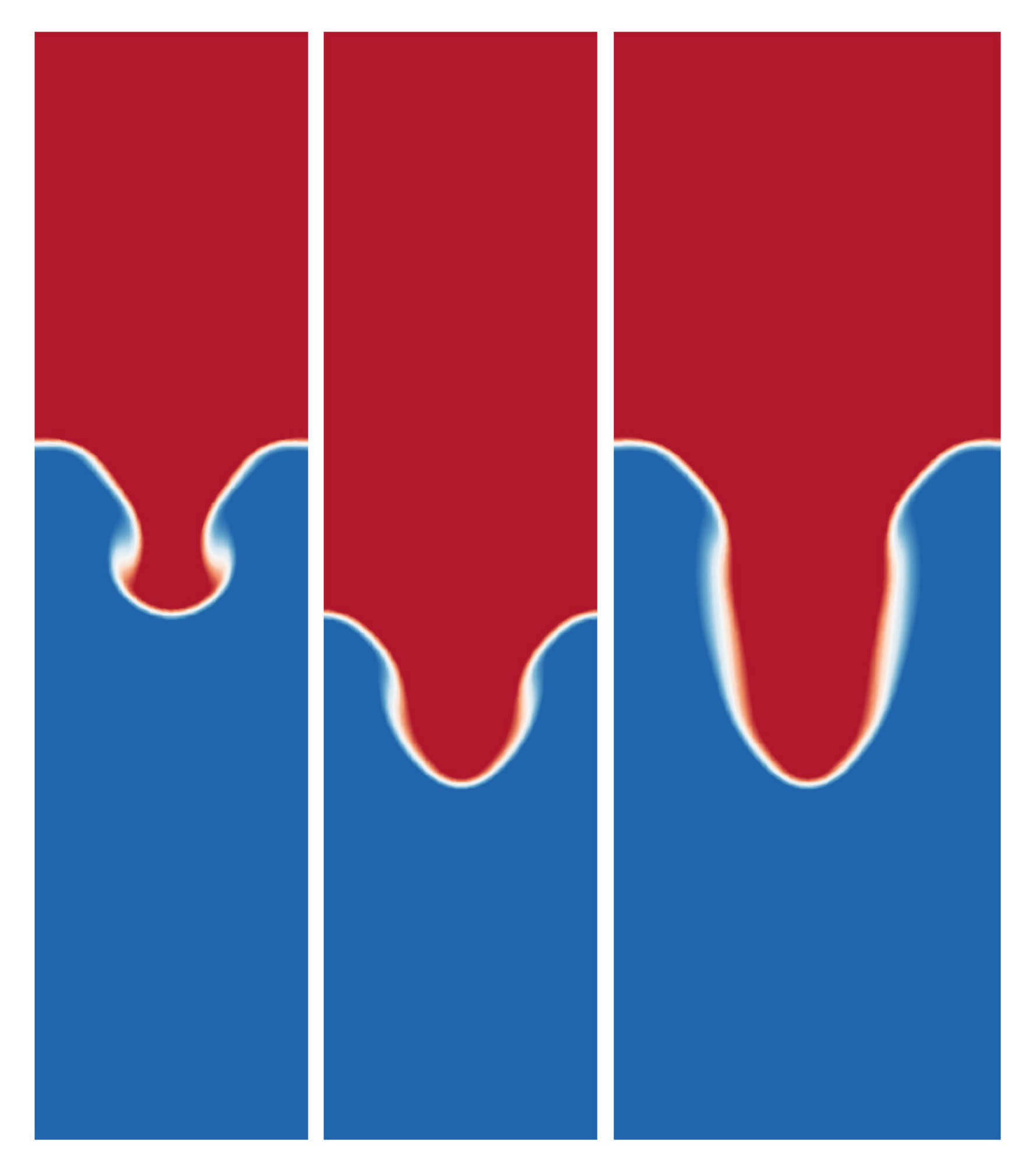}
\includegraphics[width=0.3\textwidth]{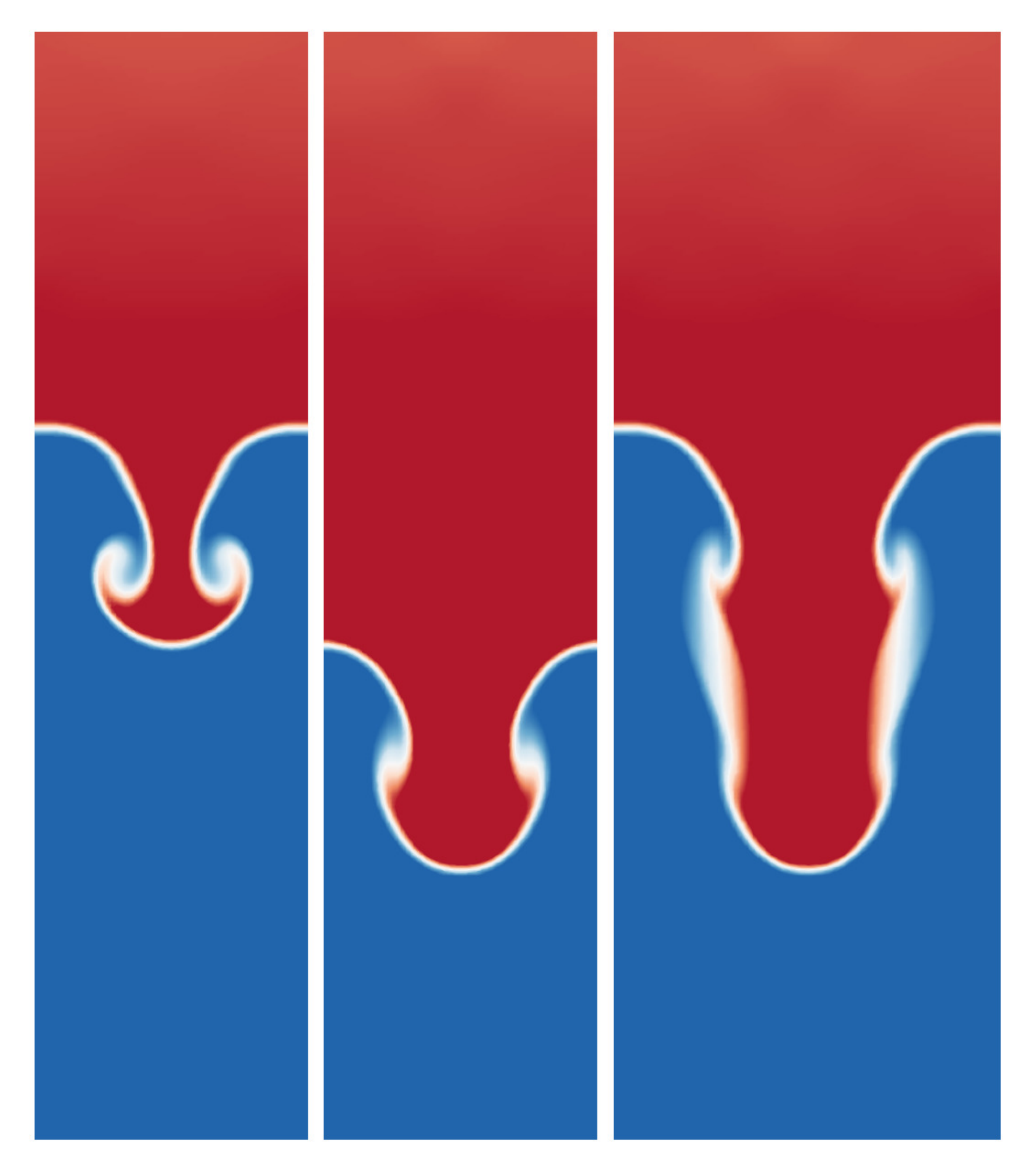}\\
\includegraphics[width=0.3\textwidth]{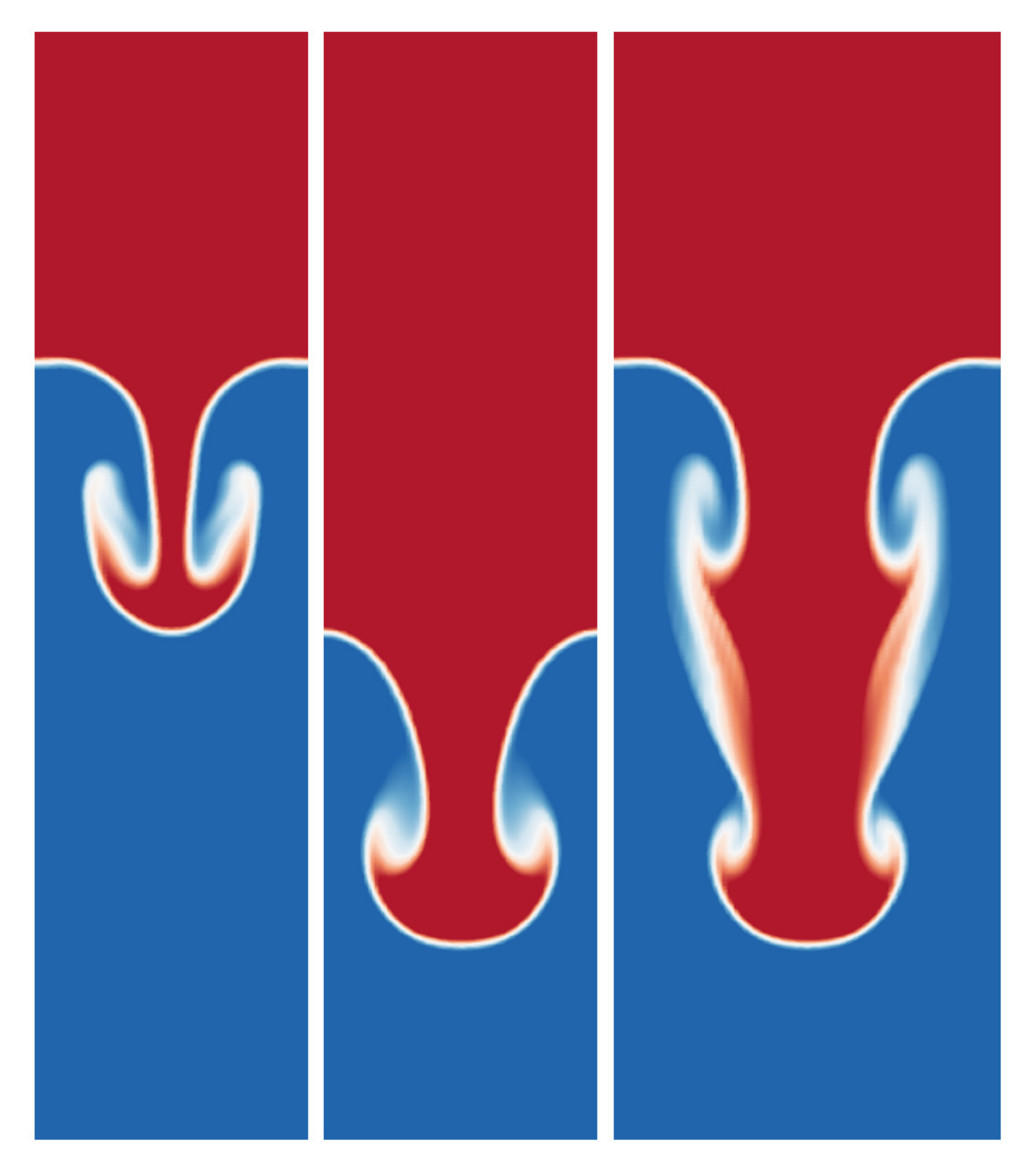}
\includegraphics[width=0.3\textwidth]{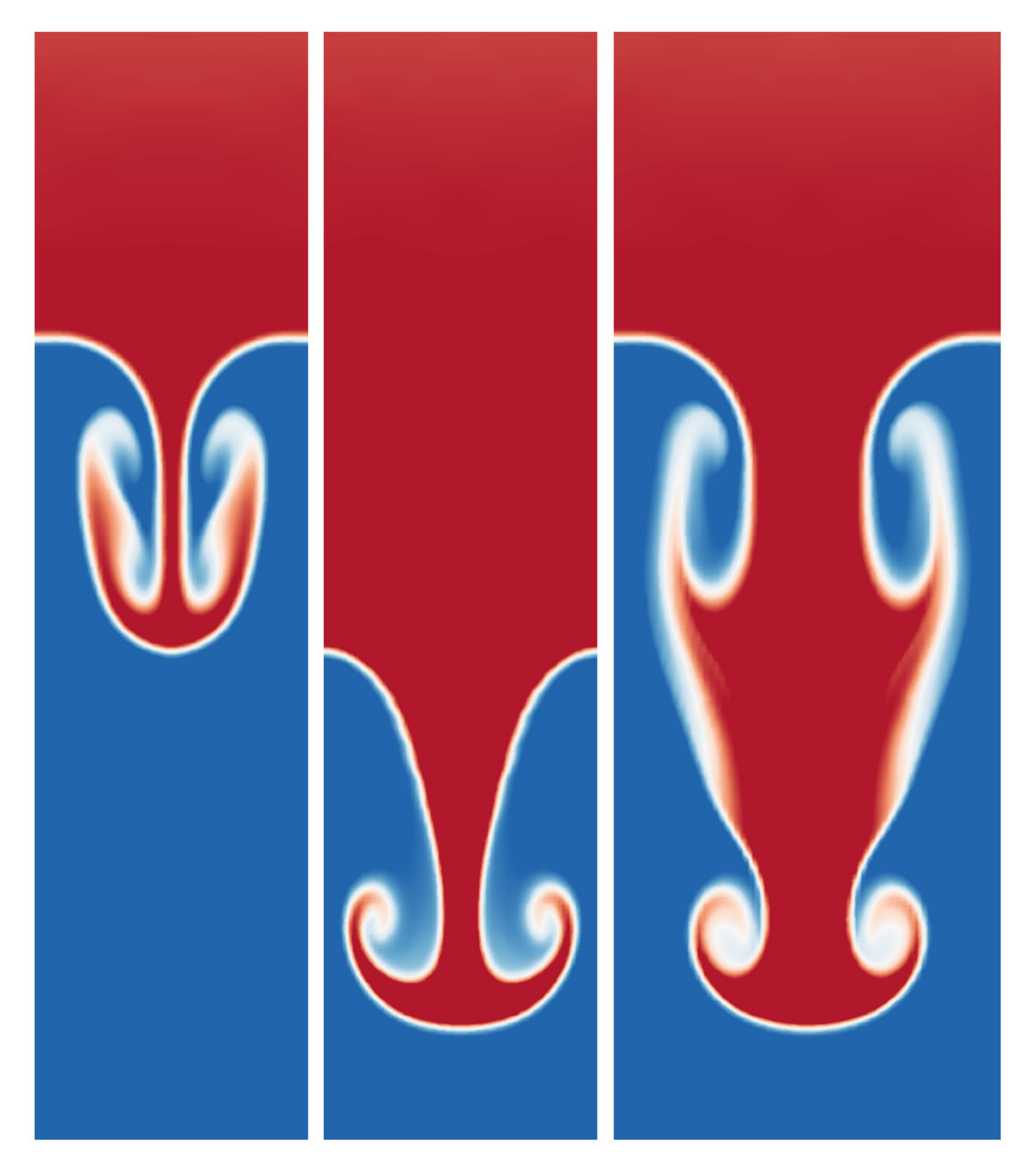}
\includegraphics[width=0.3\textwidth]{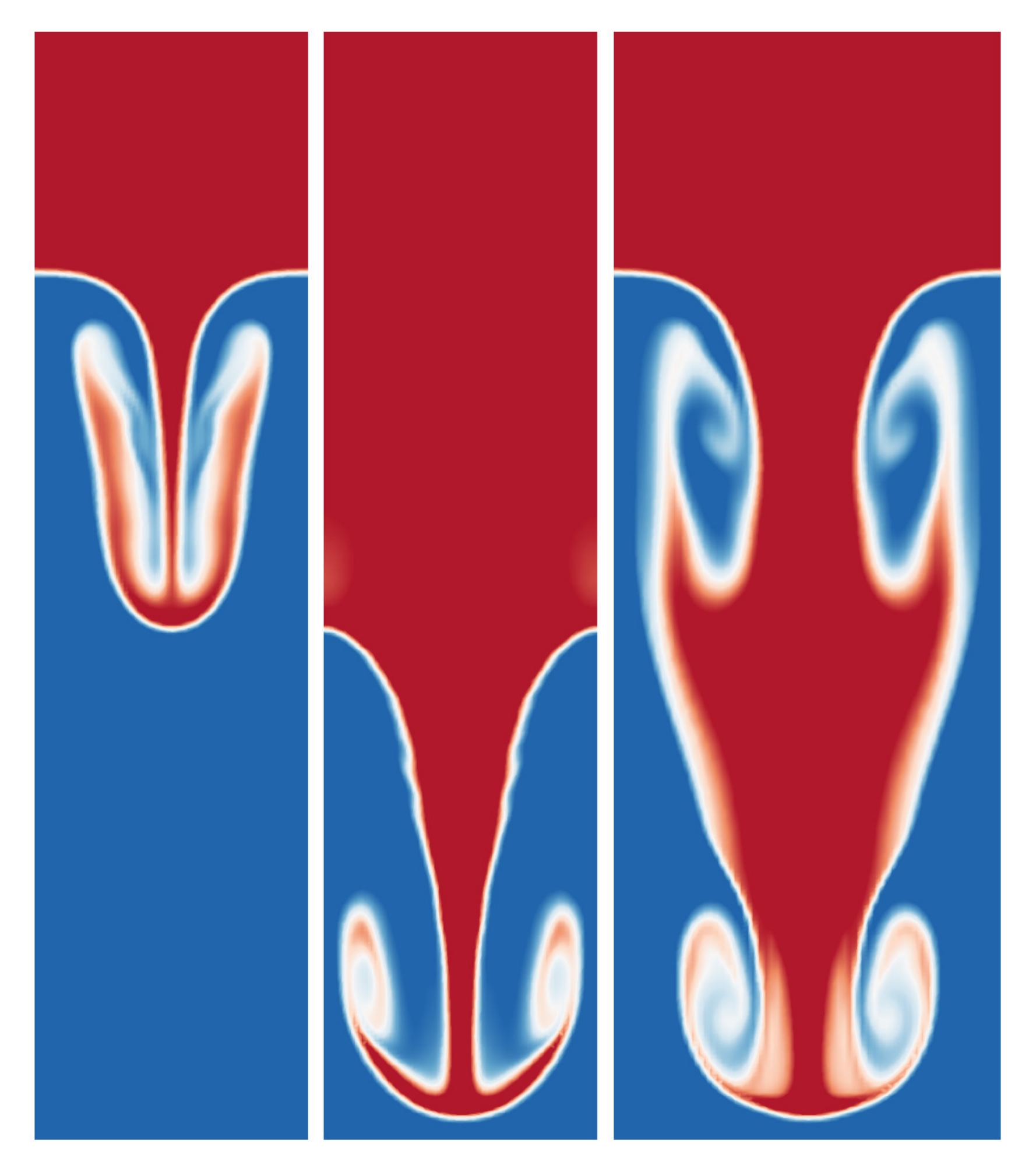}
\caption{\label{rt-interface2} Rayleigh-Taylor instability: the density distribution at $x=0$, $x=W/2$ and $x=y$ with $t=3,4,5,6,7,8$.}
\end{figure}

The Rayleigh-Taylor instability occurs when an interface between two different fluids with different densities experiences a pressure gradient opposing the density gradient \cite{Case-rt}. This interface is unstable and any disturbance tends to grow, leading to the penetration of both fluids into each other. If the initial interface displacement is random, the Rayleigh-Taylor instability usually evolves into complicated turbulent mixing. In this case, the three-dimensional Rayleigh-Taylor instability is tested in a rectangular box with domain $[0,0.25]\times[0,0.25]\times[0,1]$ with a gravity pointed downward $\textbf{g}=(0, 0,-0.1)$.
Boundary conditions are  applied at the four sides, while symmetric boundary conditions are applied at the top and bottom walls. The instability develops from the imposed single mode initial perturbation
\begin{equation*}
\frac{h(x,y)}{W}=0.05\big[\cos(\frac{2\pi x}{W})+\cos(\frac{2\pi y}{W})\big],
\end{equation*}
where $h$ is the height of the interface and $W$ is the box width. The dimensionless parameters is the Atwood number
\begin{equation*}
A=\frac{\rho_h-\rho_l}{\rho_h+\rho_l},
\end{equation*}
where $\rho_h, \rho_l$ are densities of heavy and light fluids, respectively. In the computation, the Atwood number $A=1/3$ is used. The the fluid interface from a single-mode perturbation at $t=3, 4, 5, 6, 7, 8$ are given in Fig.\ref{rt-interface}, and the density distribution at $x=0$, $x=W/2$ and $x=y$ with $t=3,4,5,6,7,8$ are given in Fig.\ref{rt-interface2}. As expected, the heavy and light fluids penetrate into each other as time increases. The light fluid rises to form a bubble and the heavy fluid falls to generate a spike.

\begin{figure}[!h]
\centering
\includegraphics[width=0.35\textwidth]{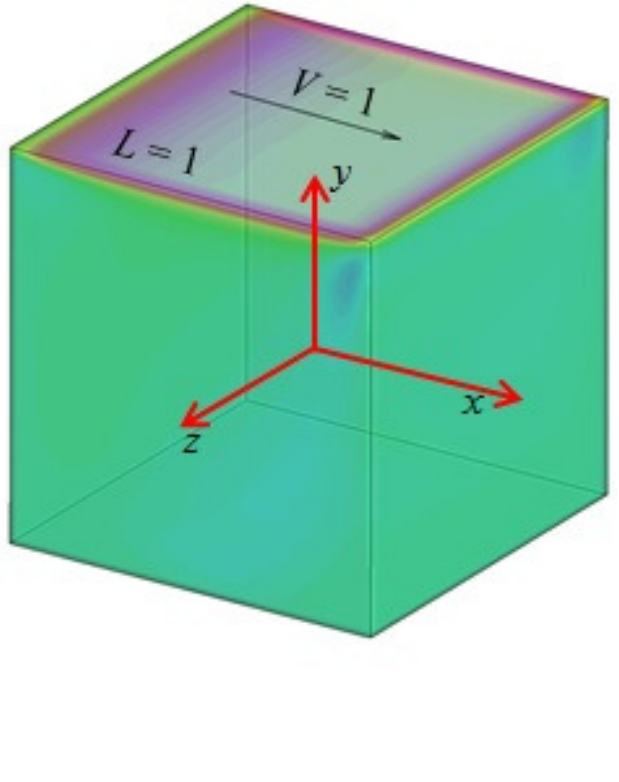}
\vspace{-17mm}
\caption{\label{cavity}3D cavity flow: schematic for 3d cavity flow.}
\end{figure}

\begin{figure}[!h]
\centering
\includegraphics[width=0.415\textwidth]{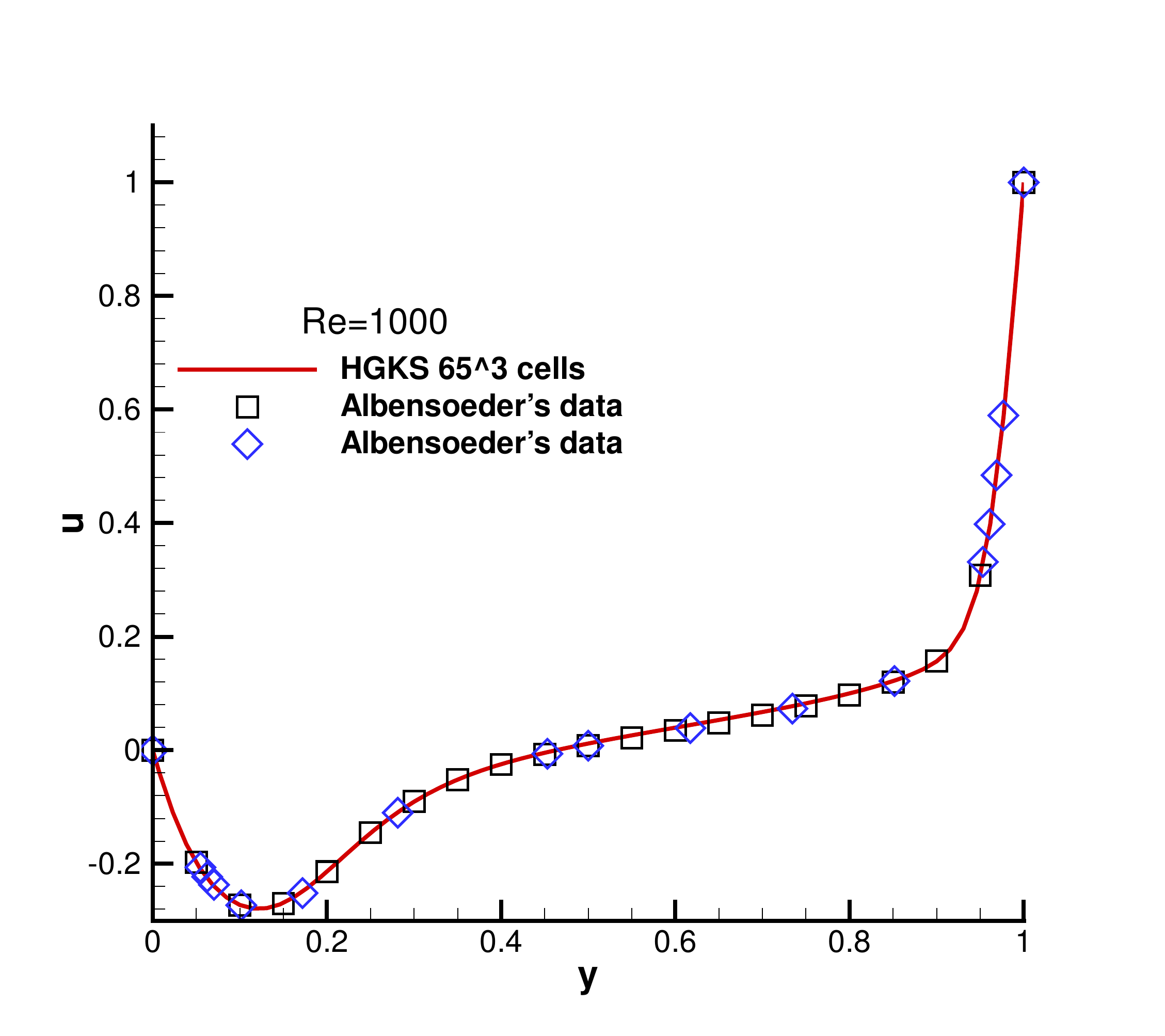}
\includegraphics[width=0.415\textwidth]{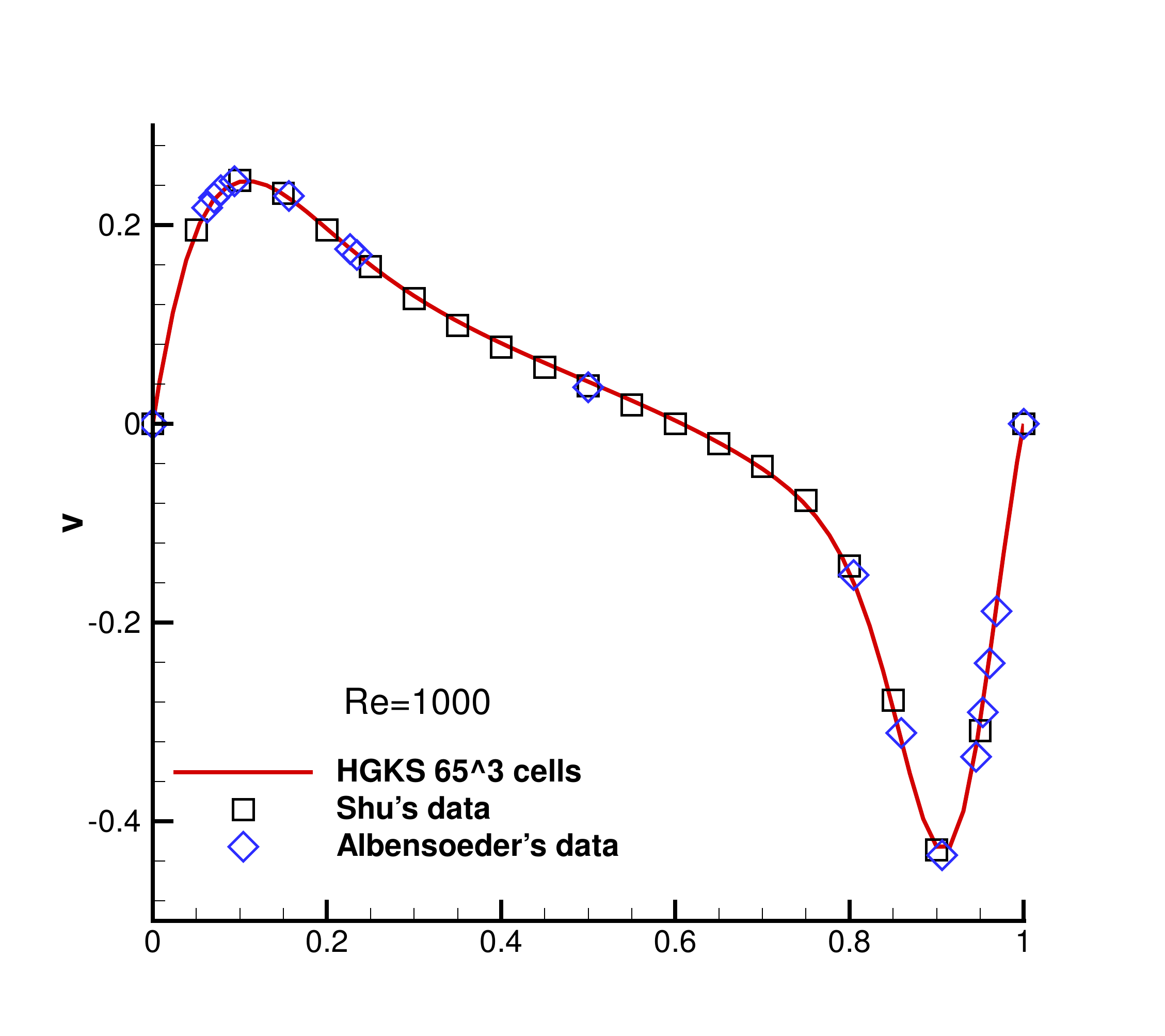}
\includegraphics[width=0.415\textwidth]{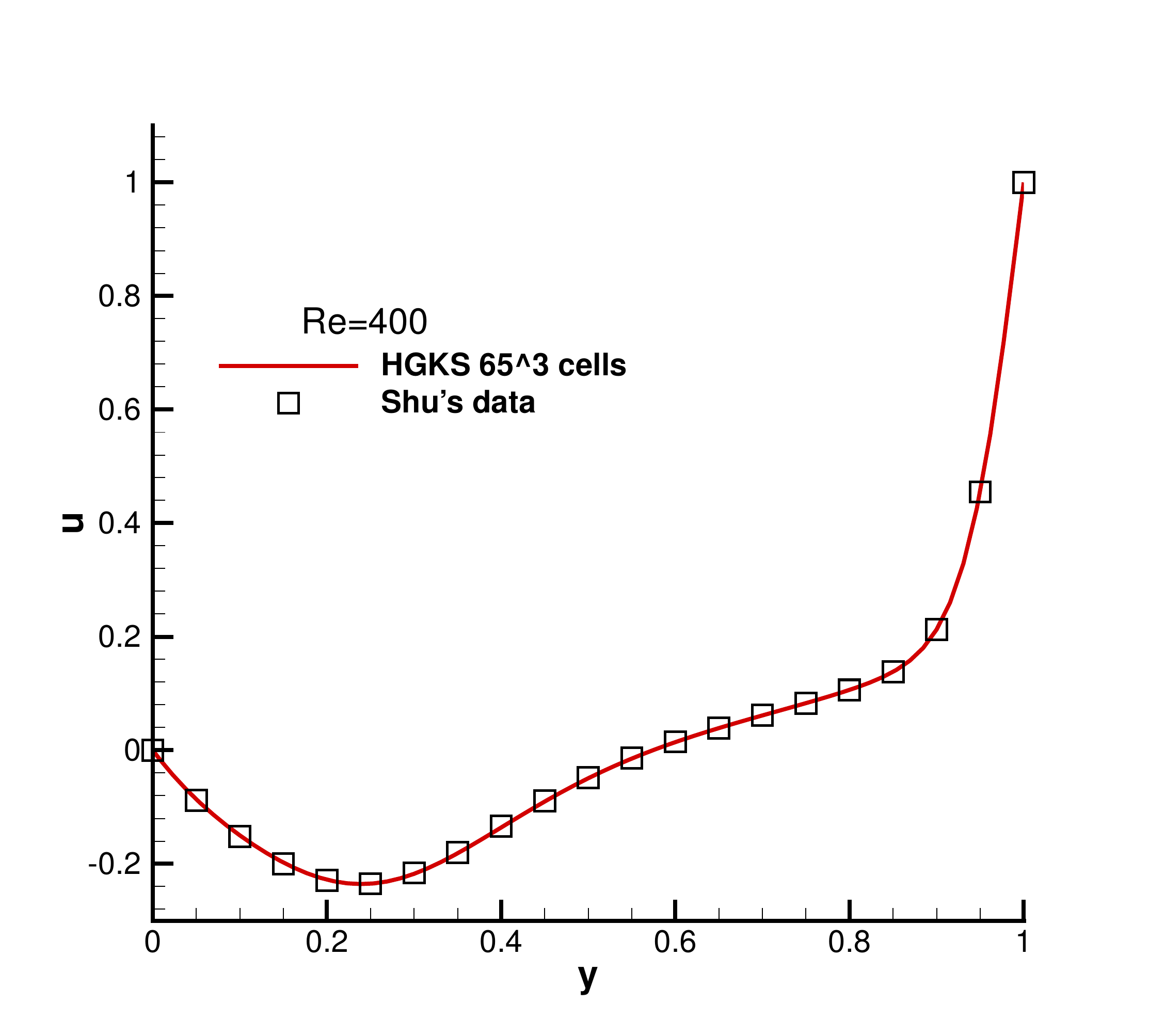}
\includegraphics[width=0.415\textwidth]{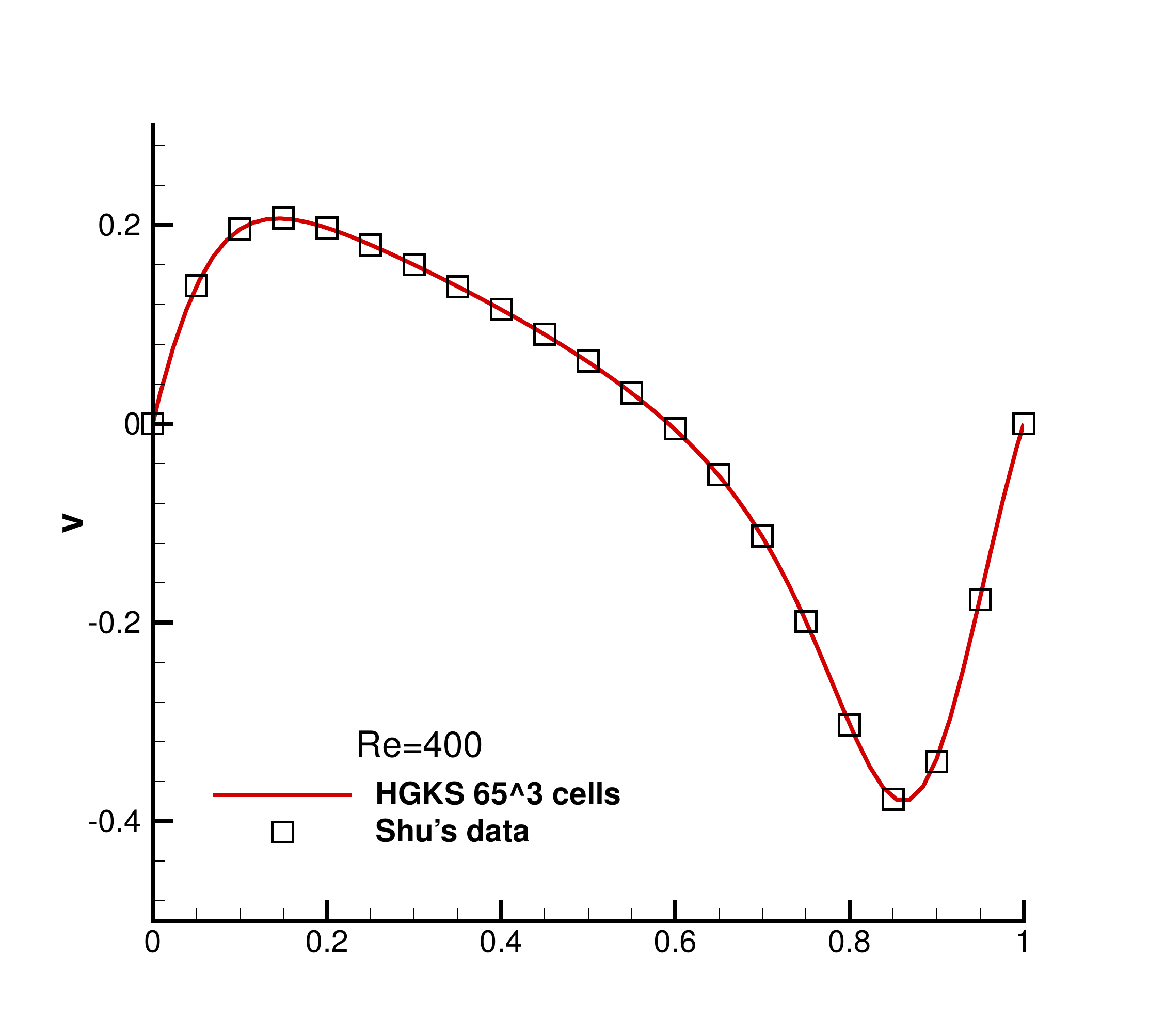}
\includegraphics[width=0.415\textwidth]{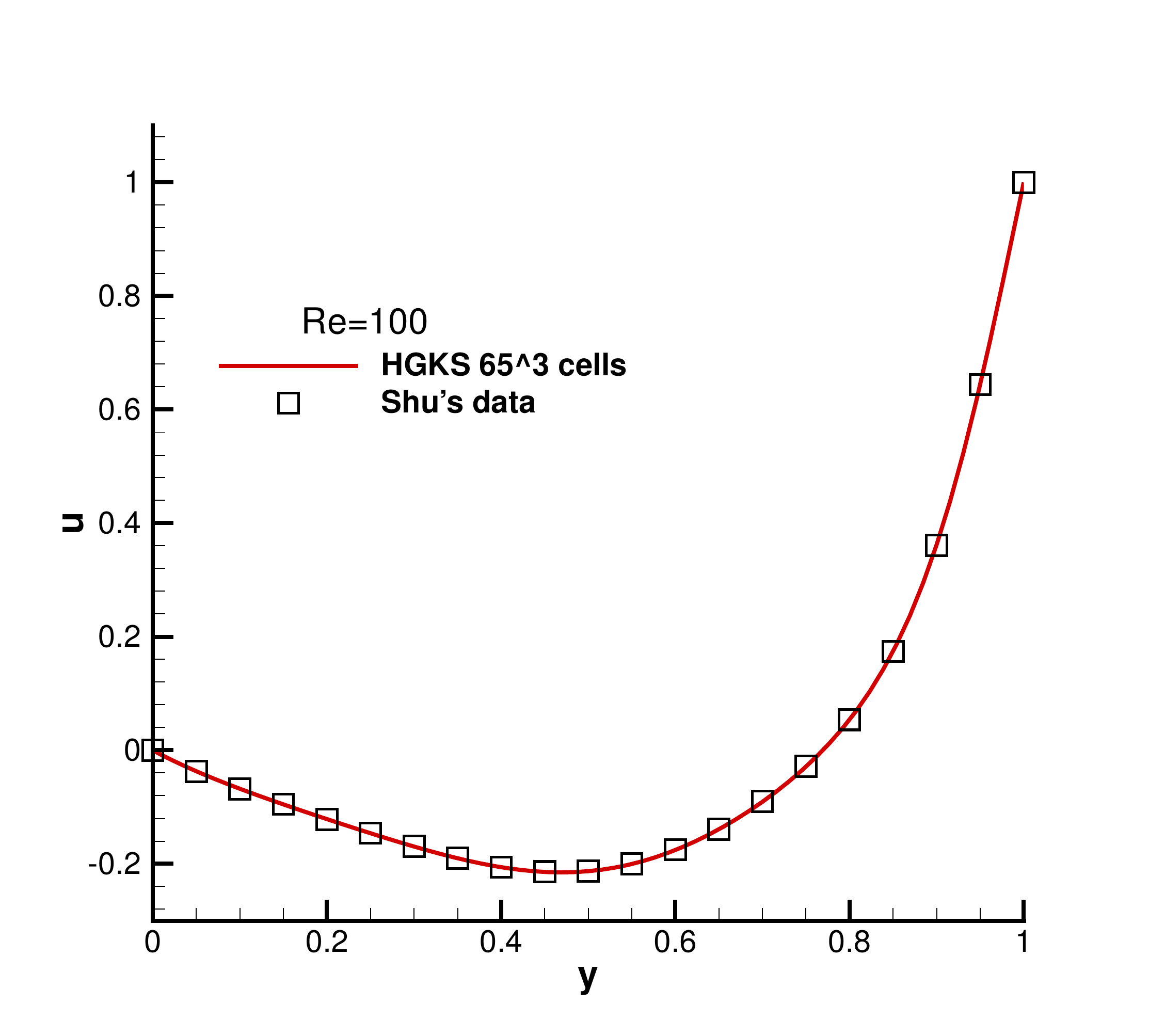}
\includegraphics[width=0.415\textwidth]{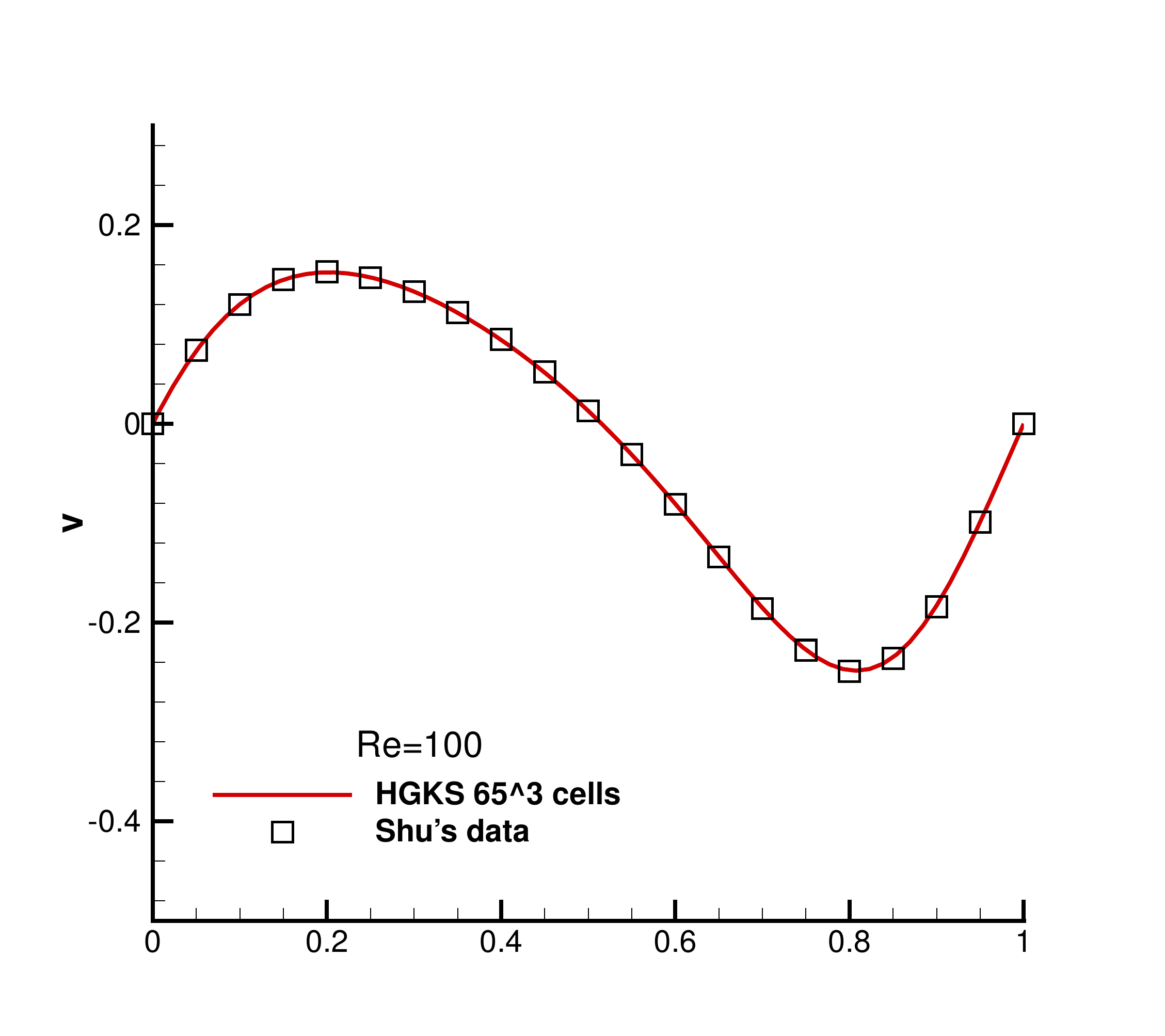}
\caption{\label{cavity-velocity-1} 3D cavity flow: the profiles of $U$-velocity along the vertical centerline line and $V$-velocity along the horizontal centerline at the $x-y$-plane with $z=0$ with $Re=1000, 400$ and $100$.}
\end{figure}

\subsection{Lid-driven cavity flow}
The lid-driven cavity problem is one of the most important benchmarks for numerical Navier-Stokes solvers. As shown in Fig.\ref{cavity}, the fluid is bounded by a unit cubic and driven by a uniform translation of
the top boundary. The monatomic gas with $\gamma = 5/3$ is used, such that there is no bulk viscosity involved. Early three-dimensional cavity-flow calculations were carried out by De Vahl Davis and Mallinson \cite{Case-Davis} and Goda \cite{Case-Goda}. In this case, the flow is simulated with Mach number $Ma=0.15$ and all the boundaries are isothermal and nonslip. Numerical simulations are conducted with three Reynolds numbers of $Re=1000, 400$ and $100$ using $65\times65\times65$ meshes for the domain $[-0.5, 0.5]\times[-0.5, 0.5]\times[-0.5, 0.5]$. The $u$-velocity profiles along the vertical centerline line, $v$-velocity profiles along the horizontal centerline in the $x-y$ plane with $z=0$ and the benchmark data \cite{Case-Albensoeder, Case-Shu} are shown in Fig.\ref{cavity-velocity-1}. The simulation results match well with the benchmark data.

\begin{figure}[!h]
\centering
\includegraphics[width=0.49\textwidth]{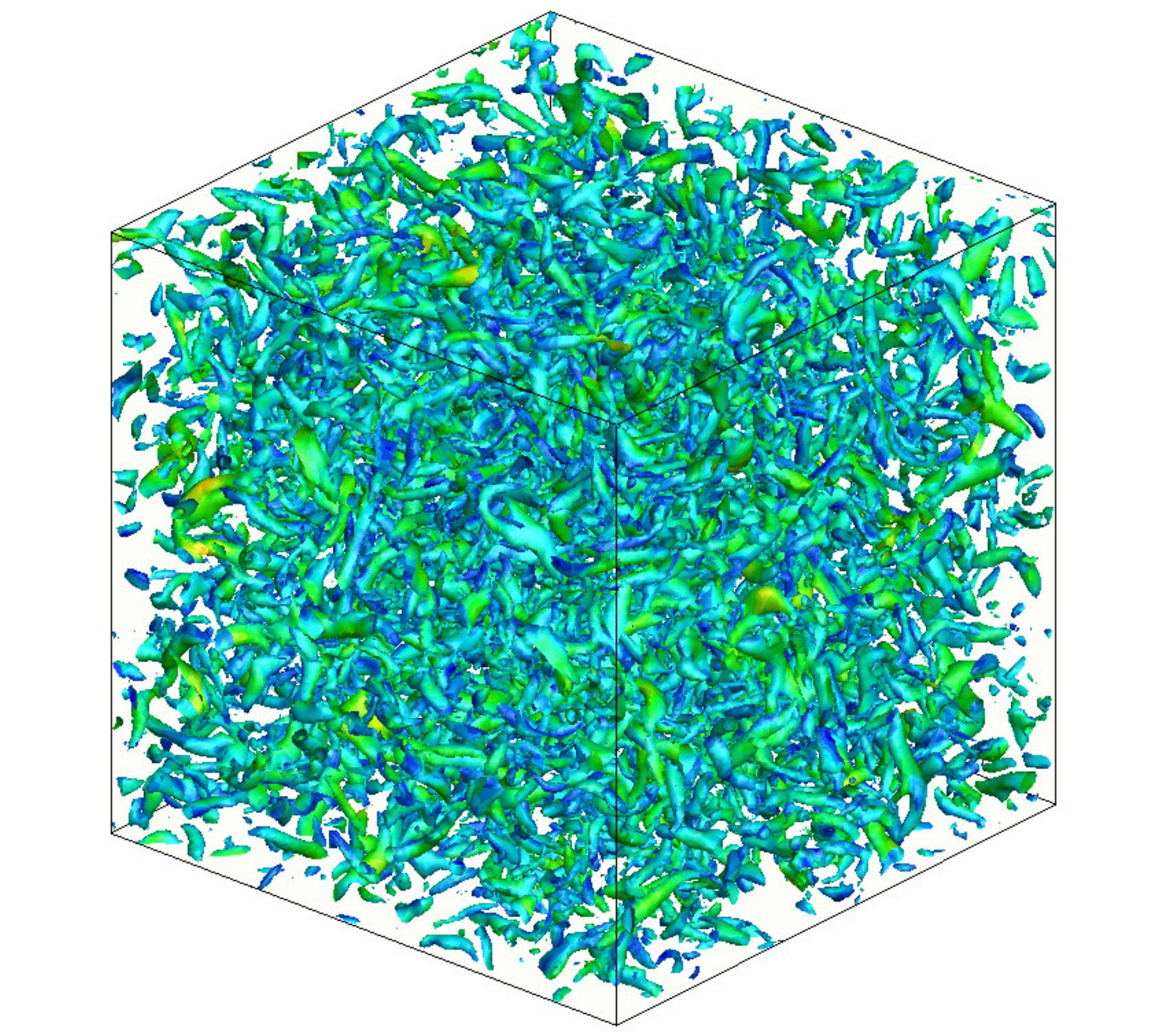}
\includegraphics[width=0.49\textwidth]{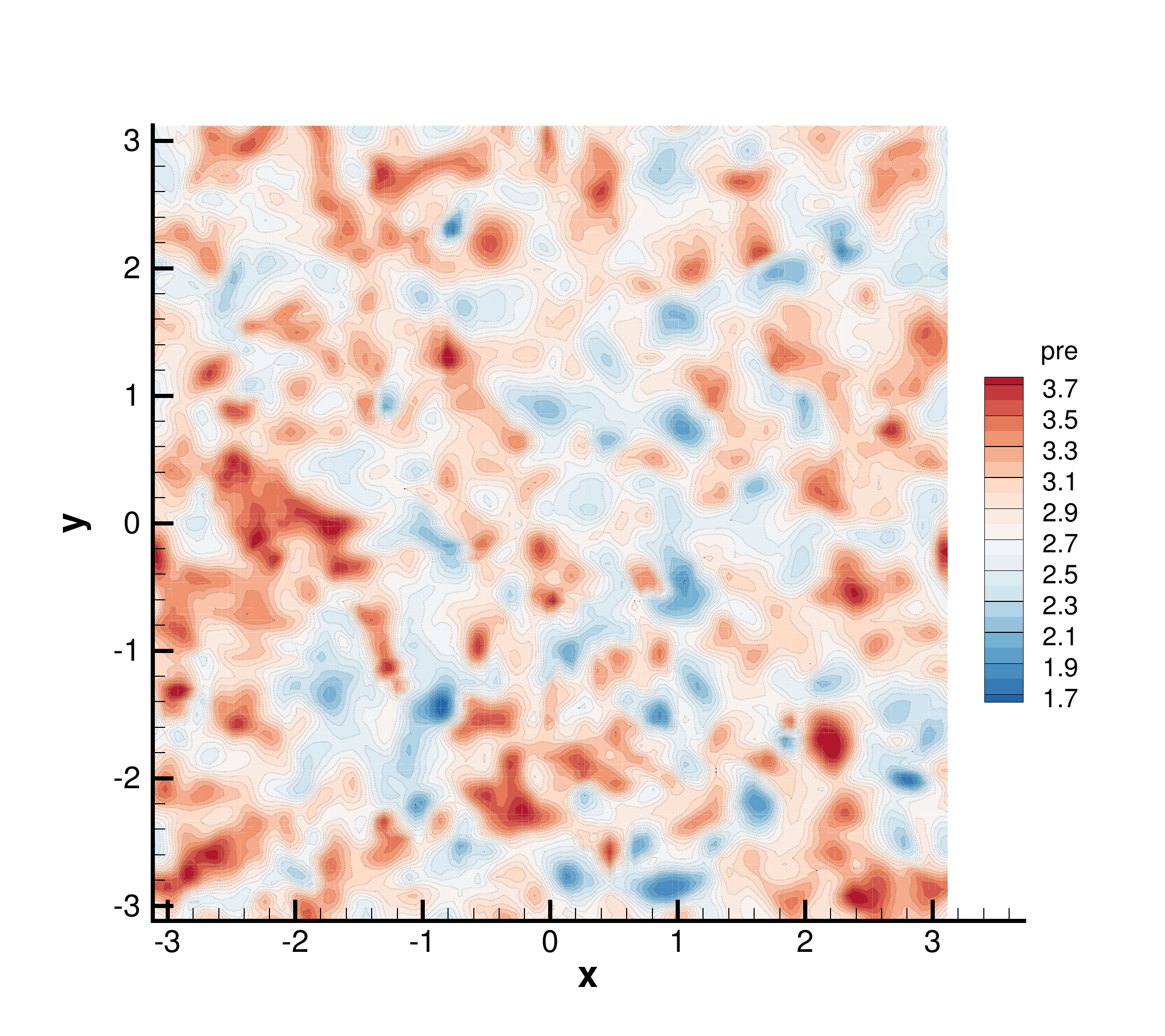}
\caption{\label{homogeneous-1} Compressible homogeneous turbulence. Left: iso-surfaces
of Q criterion colored by velocity magnitude at time $t/\tau=2$ with $128^3$ cells; right: the pressure distribution with $z=-\pi$ at time $t/\tau=1$.}
\end{figure}

\begin{figure}[!h]
\centering
\includegraphics[width=0.45\textwidth]{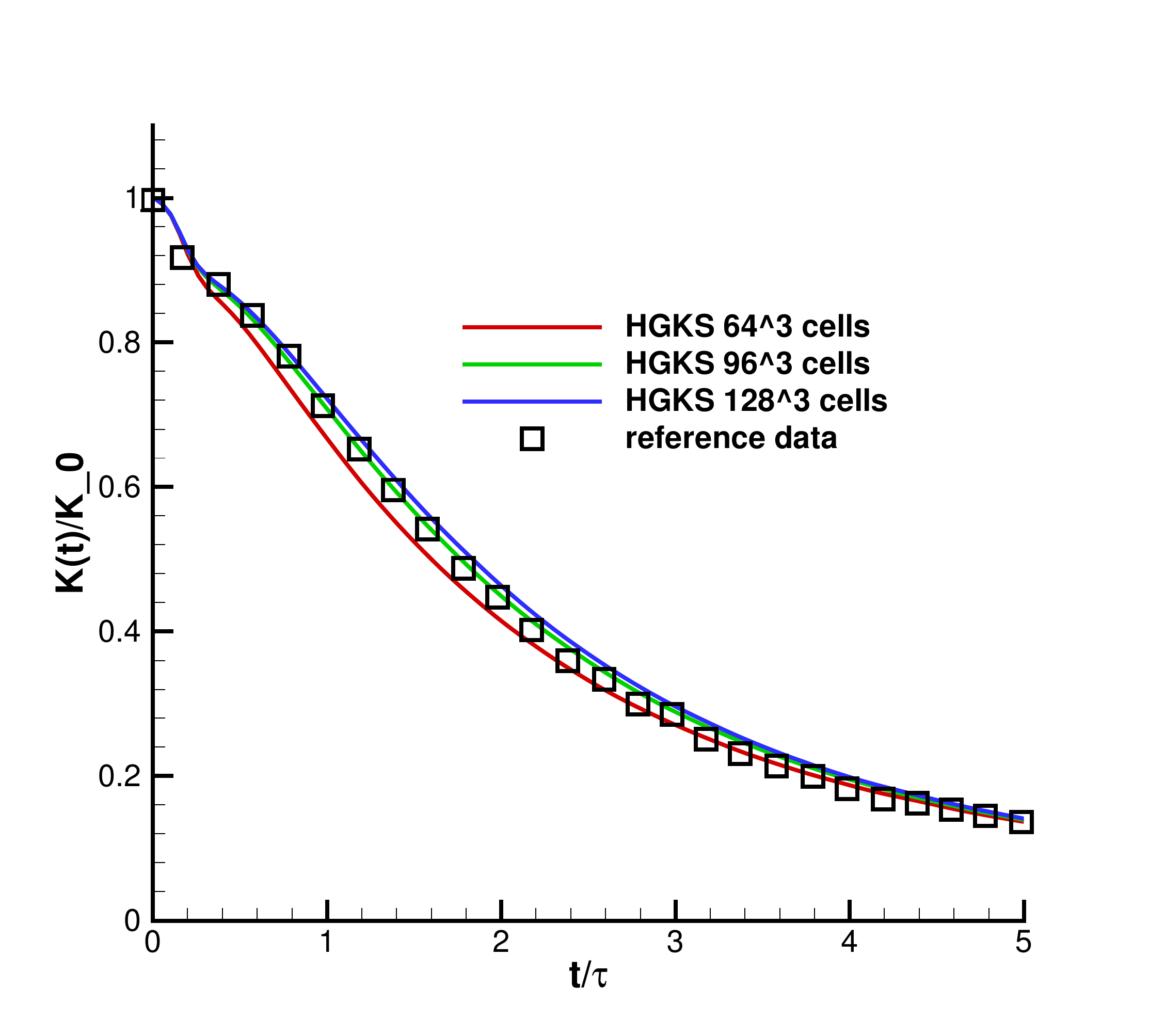}{a}
\includegraphics[width=0.45\textwidth]{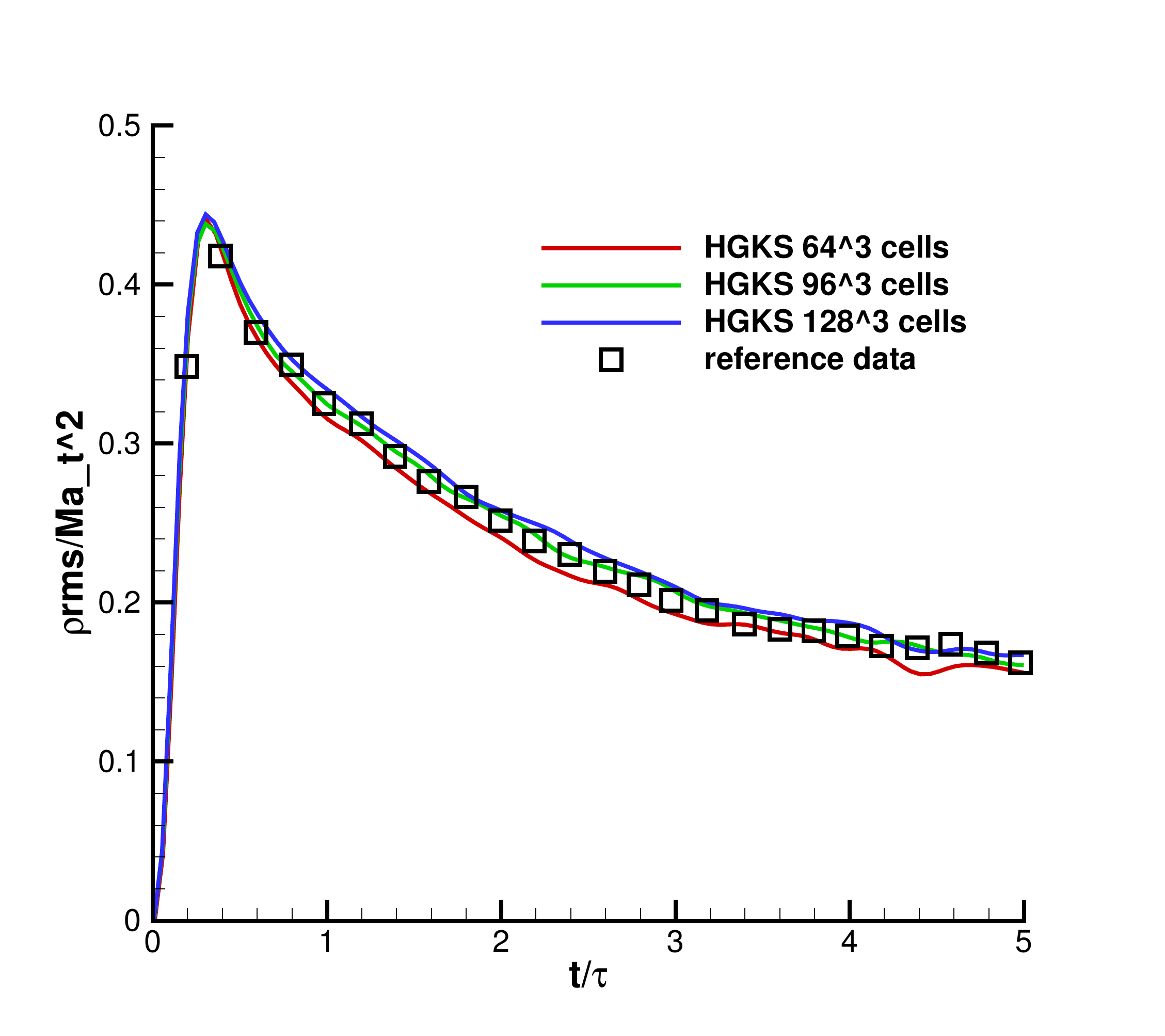}{b}\\
\includegraphics[width=0.45\textwidth]{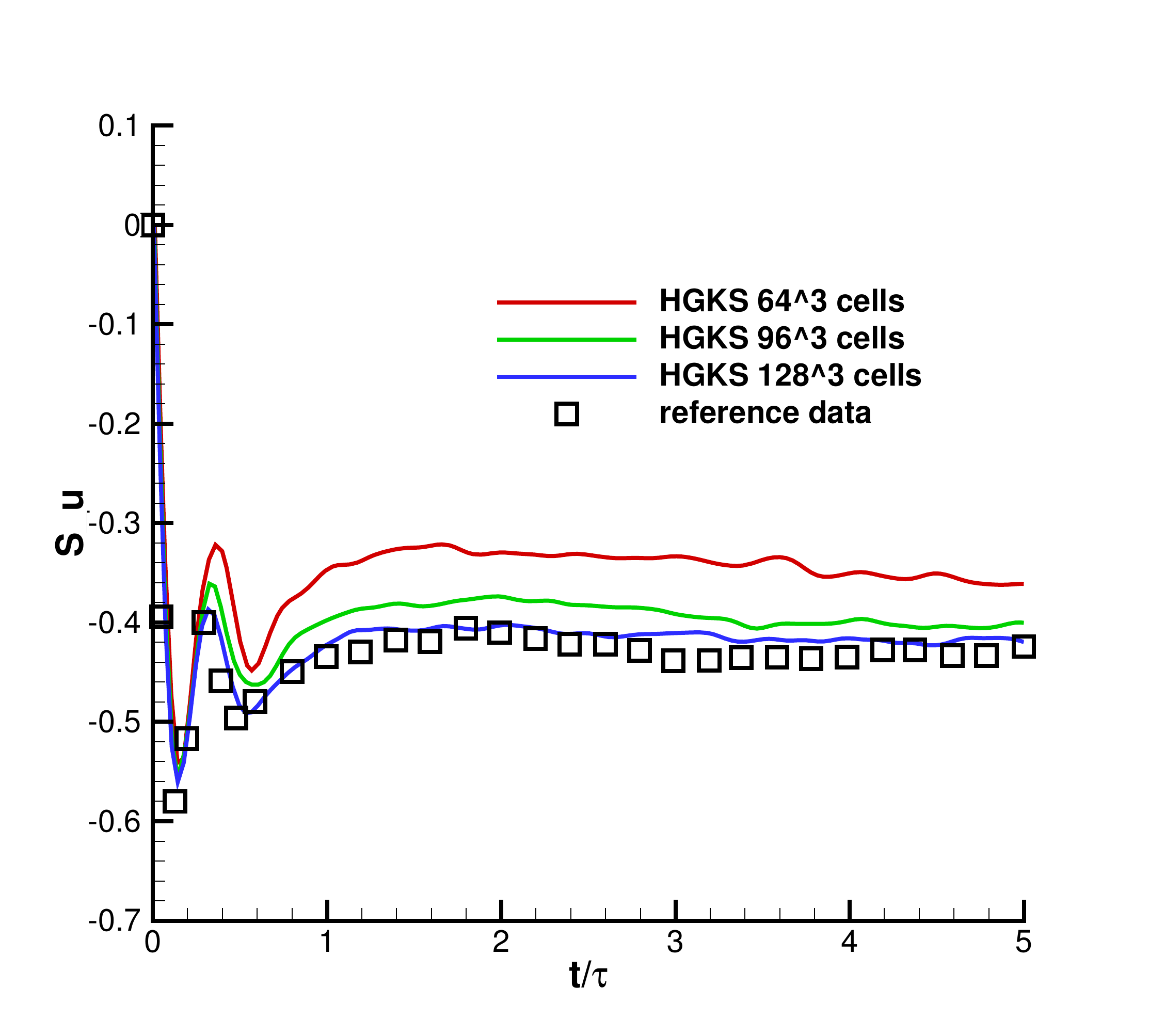}{c}
\includegraphics[width=0.45\textwidth]{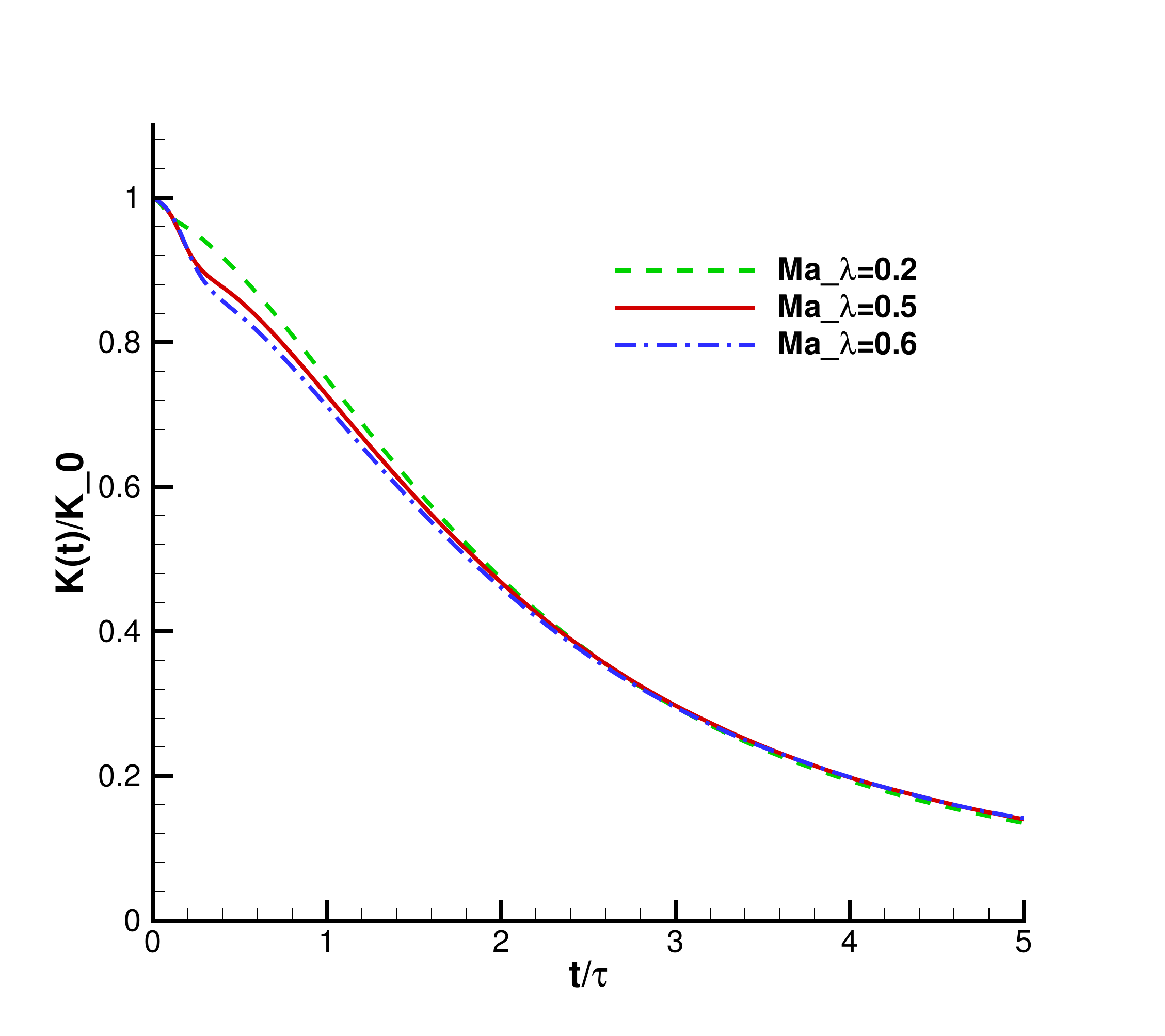}{d}
\caption{\label{homogeneous-2} Compressible homogeneous turbulence: the time history of $K(t)/K_0$ (a), $\rho_{rms}(t)/Ma_t^2$ (b) and $S_u(t)$ and time history of $K(t)/K_0$ (c) with $Ma_\lambda=0.2,0.5$ and $0.6$ with respect to $t/\tau$ (d).}
\end{figure}

\subsection{Compressible homogeneous turbulence}
The high-order gas-kinetic scheme is applied for the direct numerical simulations (DNS) of the compressible decaying homogeneous isotropic turbulence. The flow is computed within a square box defined as $-\pi \leq x, y, z\leq \pi$, and the periodic boundary conditions are used in all directions for all the flow variables \cite{DNS-1,DNS-2,DNS-3}. In the computation, the
domain is discretized with a uniform Cartesian mesh cells $N^3$. A divergence-free random initial velocity field $\textbf{u}_0$ is generated
for a given spectrum with a specified root mean square $u'$
\begin{align*}
u'=<\frac{\textbf{u}\cdot \textbf{u}}{3}>,
\end{align*}
where $<...>$ is a volume average over the whole computational domain. The specified spectrum for velocity is given by
\begin{align*}
E(k)=A_0k^4\exp(-2k^2/k_0^2),
\end{align*}
where $k$ is the wave number, $k_0$ is the wave number at spectrum peaks, $A$ is a constant chosen to get a specified initial kinetic energy.  The initial volume averaged turbulent kinetic energy $K_0$  and the initial large-eddy-turnover time $\tau$ is given by
\begin{align*}
K_0=\frac{3A_0}{64}\sqrt{2\pi}k_0^5,~~\tau=\frac{32}{A_0}(2\pi)^{1/4}k_0^{-7/2}.
\end{align*}
The Taylor microscale Reynolds number $Re_\lambda$ and turbulence Mach number $Ma_t$ are given as
\begin{align*}
Re_\lambda&=\frac{<\rho>u'\lambda}{<\mu>}=\frac{(2\pi)^{1/4}}{4}\frac{\rho_0}{\mu_0}\sqrt{2A_0}k_0^{3/2},\\
Ma_t&=\frac{\sqrt{3}u'}{<c_s>}=\frac{\sqrt{3}u'}{\sqrt{\gamma T_0}},
\end{align*}
where $\lambda$ is Taylor microscale
\begin{align*}
\lambda^2=\frac{(u')^2}{<(\partial_1 u_1)^2>}.
\end{align*}
$\mu_0$ and $T_0$ can be determined from $Re_\lambda$ and $Ma_t$ with initialized $u'$ and $\rho_0=1$, and the dynamic viscosity is determined by
\begin{align*}
\mu=\mu_0\big(\frac{T}{T_0}\big)^{0.76}.
\end{align*}
The time history of the kinetic energy, root-mean-square of density fluctuation and skewness factor for velocity slope are defined as
\begin{align*}
K(t)&=\frac{1}{2}<\rho \textbf{u}\cdot \textbf{u}>,\\
\rho_{rms}(t)&=\sqrt{<(\rho-\overline{\rho})^2>},\\
S_u(t)=&\sum_i\frac{<(\partial_i u_i)^3>}{<(\partial_i u_i)^2>^{3/2}}.
\end{align*}
In the computation, $A_0=1.3\times10^{-4}, k_0=8$, $Re_\lambda=72$ and $Ma_t=0.5$, and the uniform meshes with $64^3$, $96^3$ and $128^3$ cells are used. The iso-surfaces of $Q$ criterion colored by velocity magnitude and the pressure distribution with $z=-\pi$ at time $t/\tau=1$  are given in Fig.\ref{homogeneous-1}. The time history of normalized kinetic energy $K(t)/K_0$, normalized root-mean-square of density fluctuation $\rho_{rms}(t)/Ma_t^2$  and skewness factor $S_u(t)$ with respect to $t/\tau$ are given in Fig.\ref{homogeneous-2}. The numerical results agree well with the reference data. With fixed initial $Re_\lambda=72$, the cases with $Ma_t=0.2, 0.5, 0.6$ are tested, and the time histories of normalized kinetic energy $K(t)/K_0$ are given in Fig.\ref{homogeneous-2} as well. With the increase of $Ma_t$, the dynamic viscosity increases, and the kinetic energy gets dissipated more rapidly.

\begin{figure}[!h]
\centering
\includegraphics[width=0.425\textwidth]{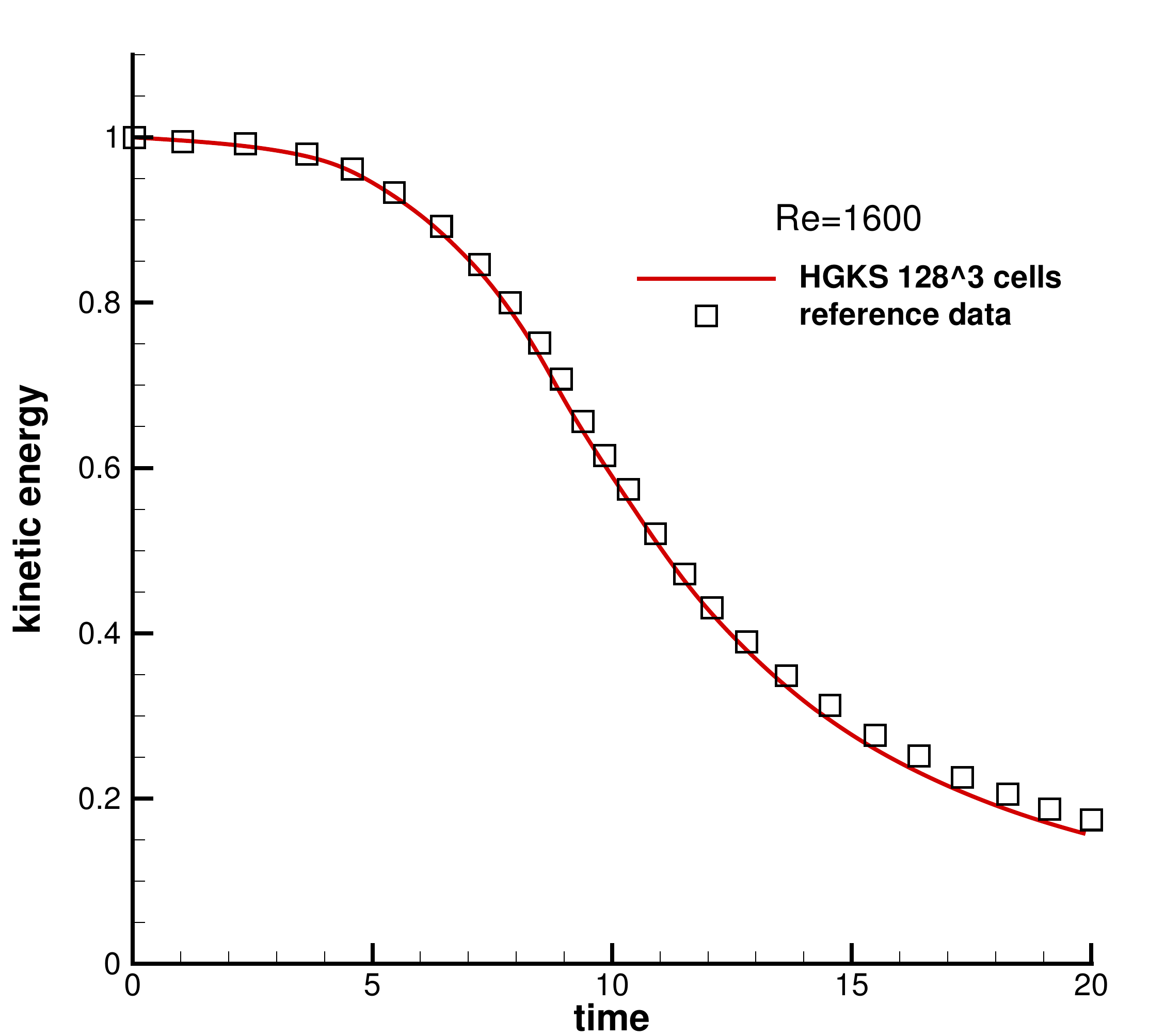}
\includegraphics[width=0.425\textwidth]{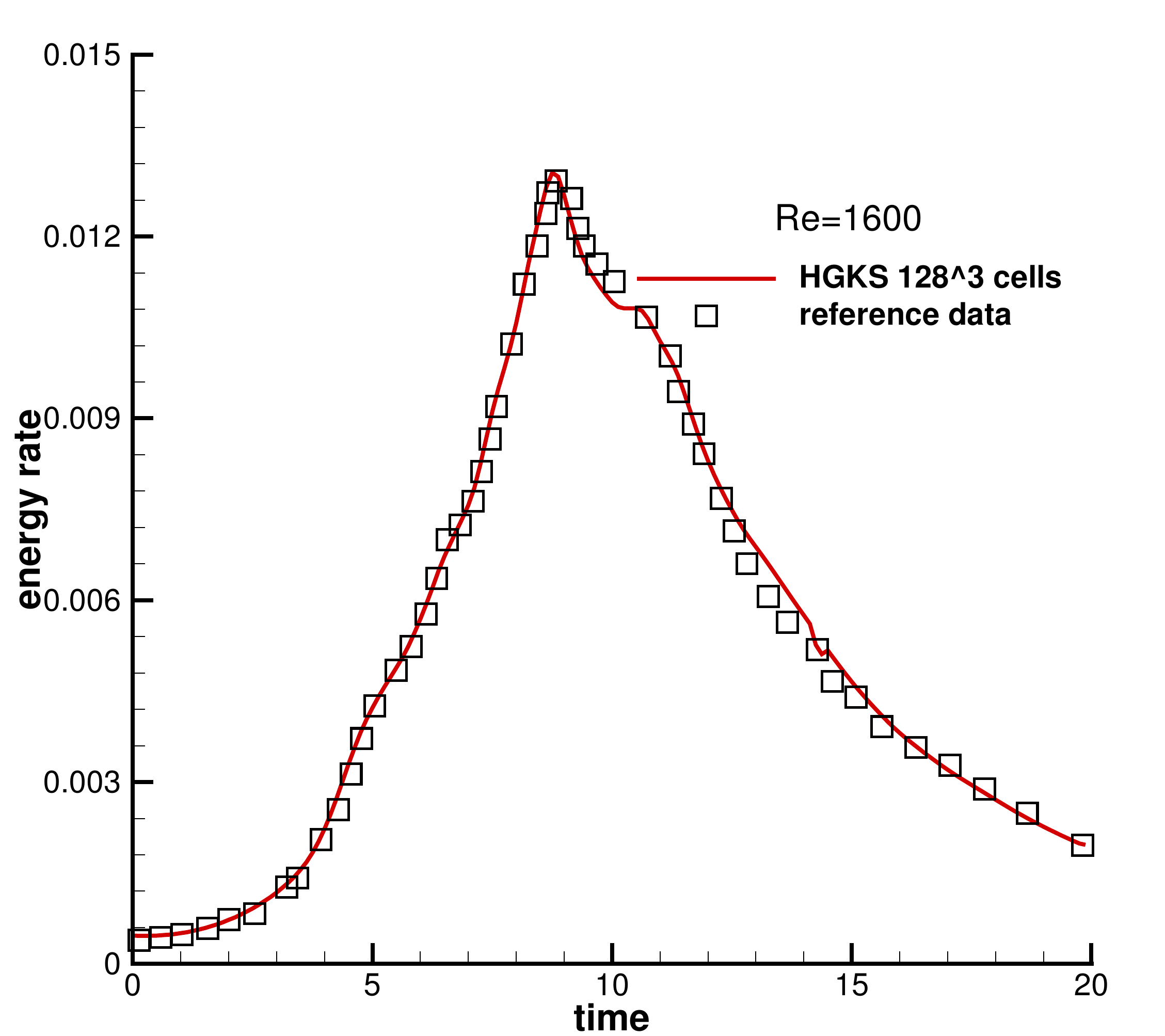}
\caption{\label{tg-vortex1} Taylor-Green Vortex problem: kinetic energy $E_k$ and dissipation rate $-dk/dt$ with fourth-order scheme for $Re=1600$.}
\includegraphics[width=0.425\textwidth]{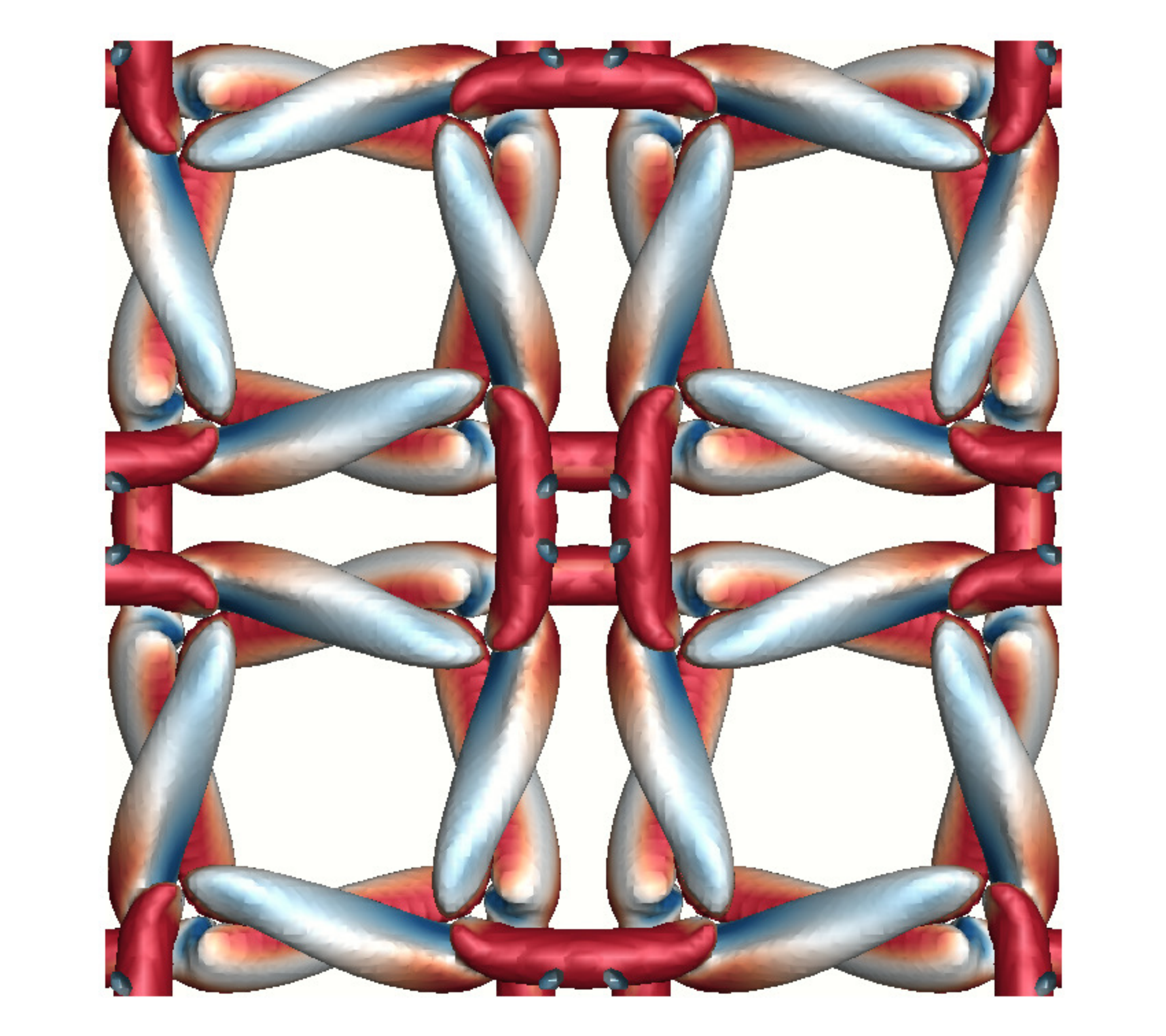}
\includegraphics[width=0.425\textwidth]{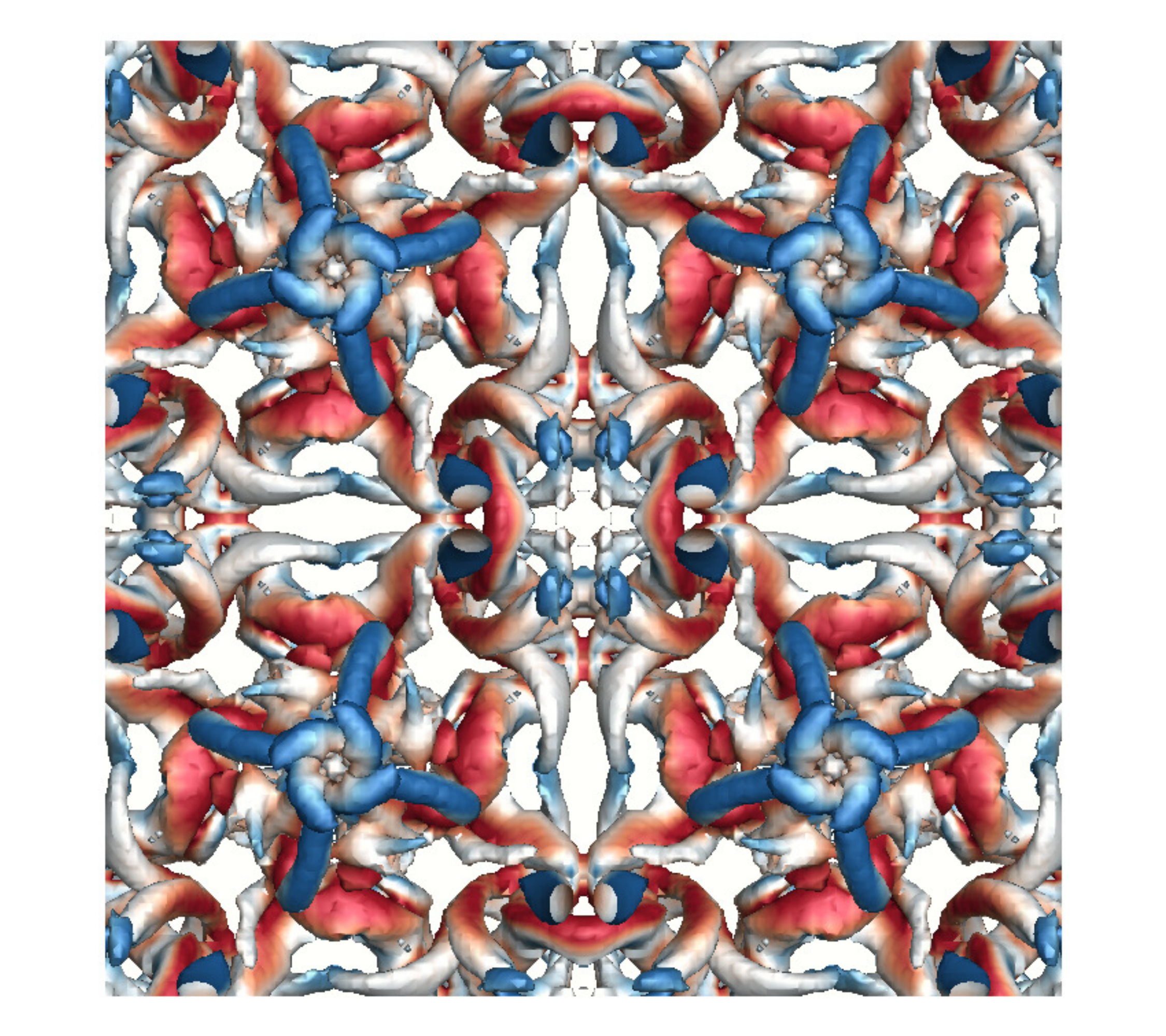}
\caption{\label{tg-vortex2} Taylor-Green Vortex problem: iso-surfaces of $Q$ criterion colored by velocity magnitude at time $t =5, 10$ for $Re=1600$.}
\end{figure}

\begin{figure}[!h]
\centering
\includegraphics[width=0.425\textwidth]{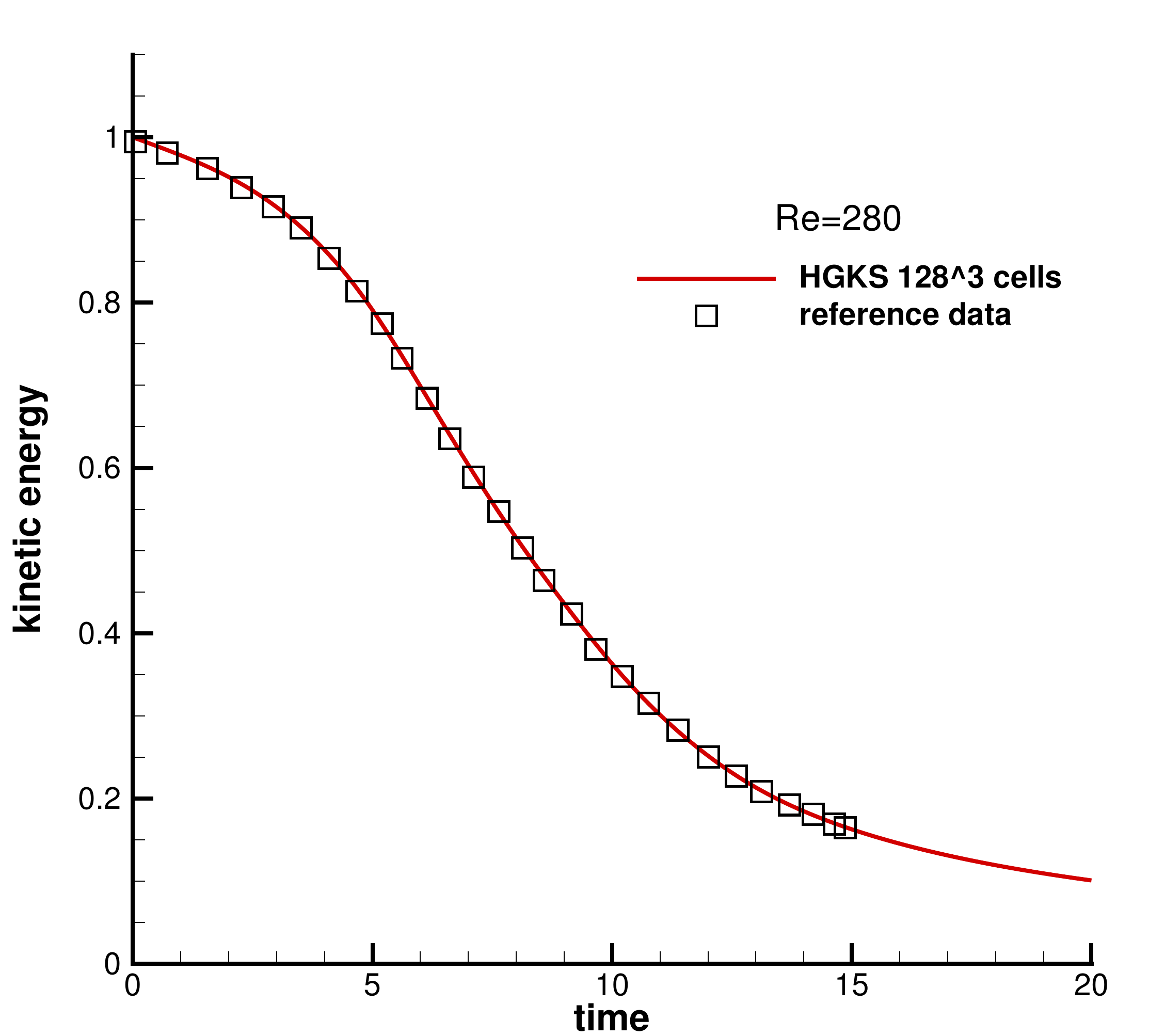}
\includegraphics[width=0.425\textwidth]{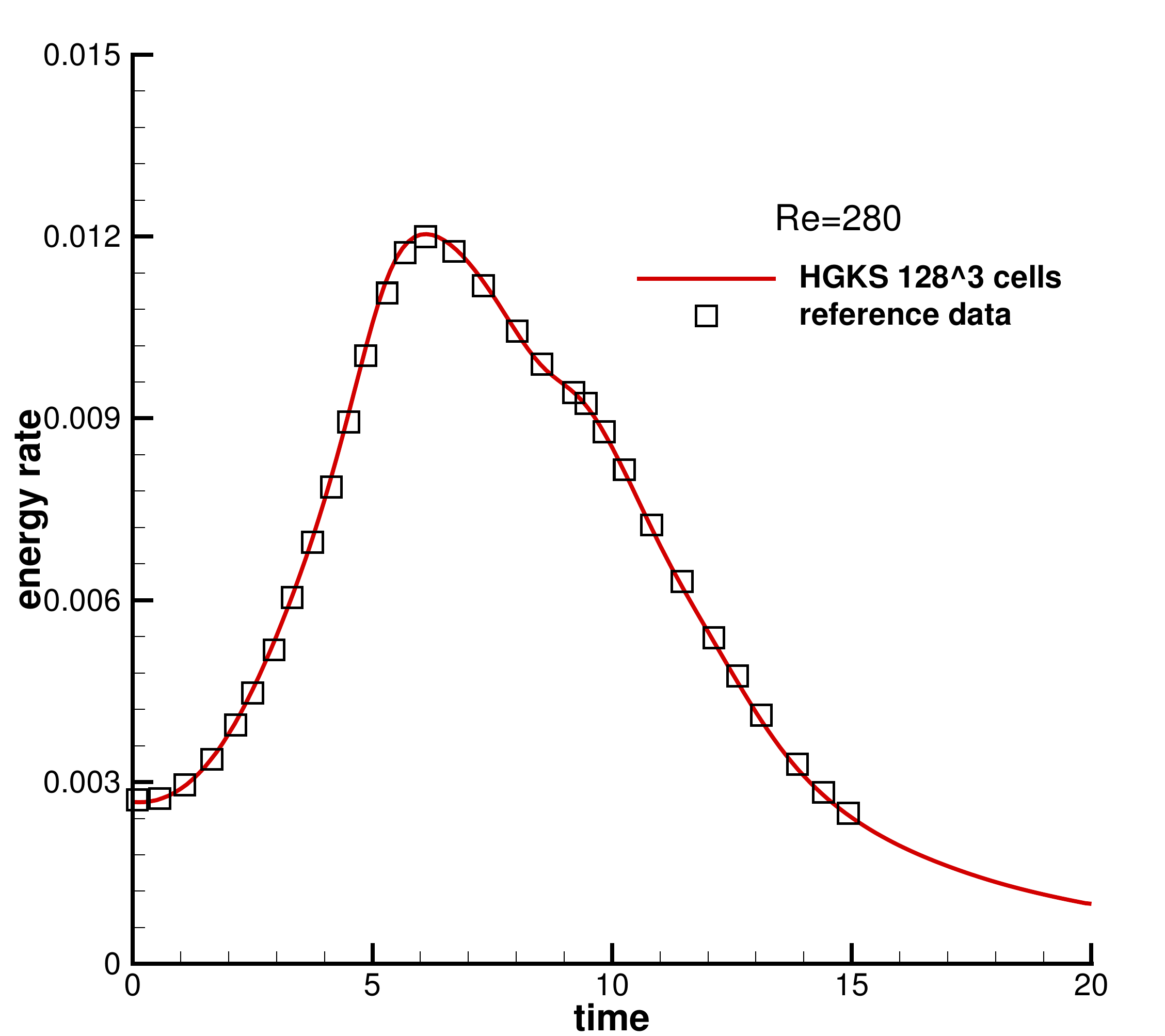}
\caption{\label{tg-vortex3} Taylor-Green Vortex problem: kinetic energy $E_k$ and dissipation rate $-dk/dt$ with fourth-order scheme for $Re=280$.}
\includegraphics[width=0.425\textwidth]{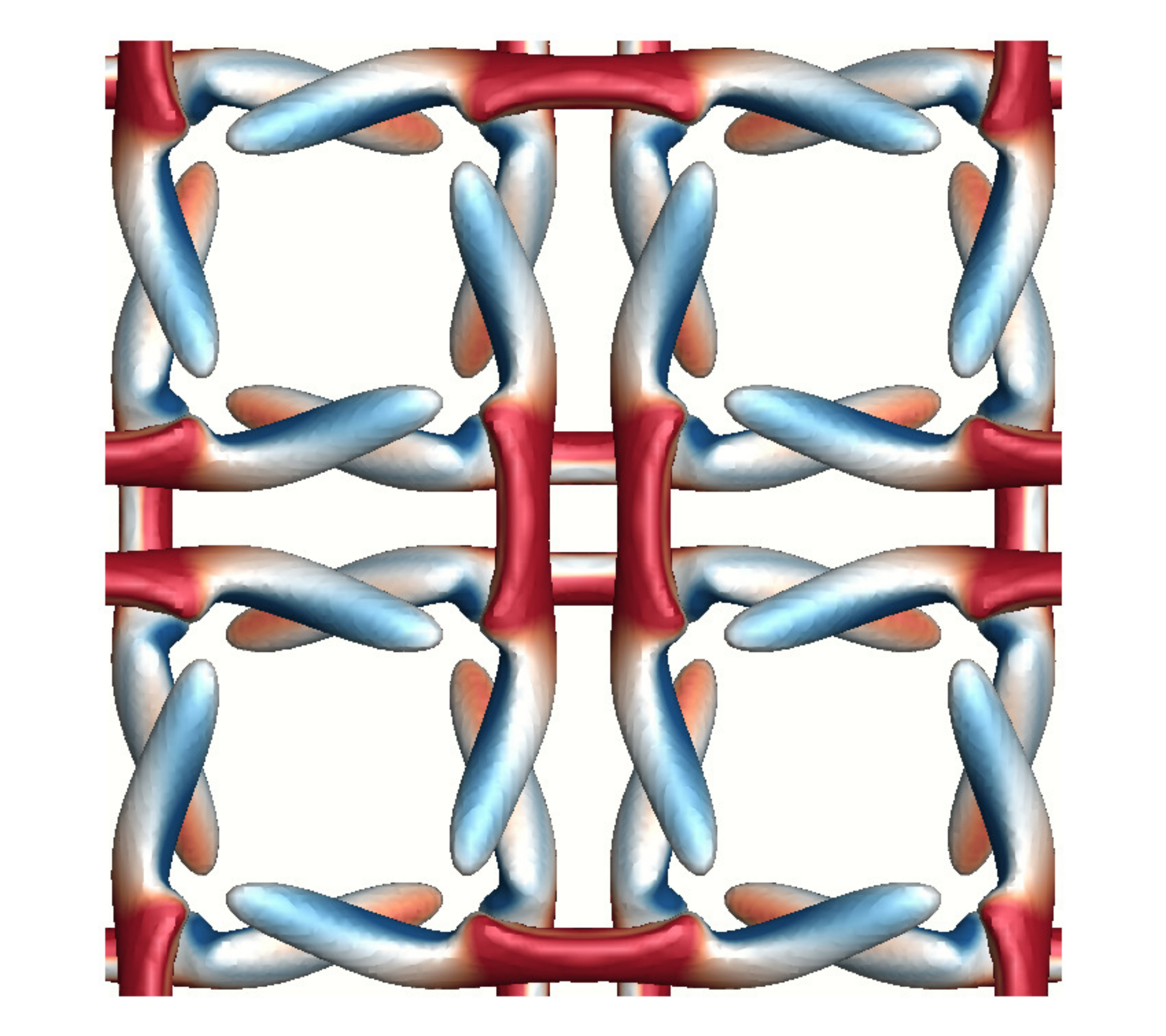}
\includegraphics[width=0.425\textwidth]{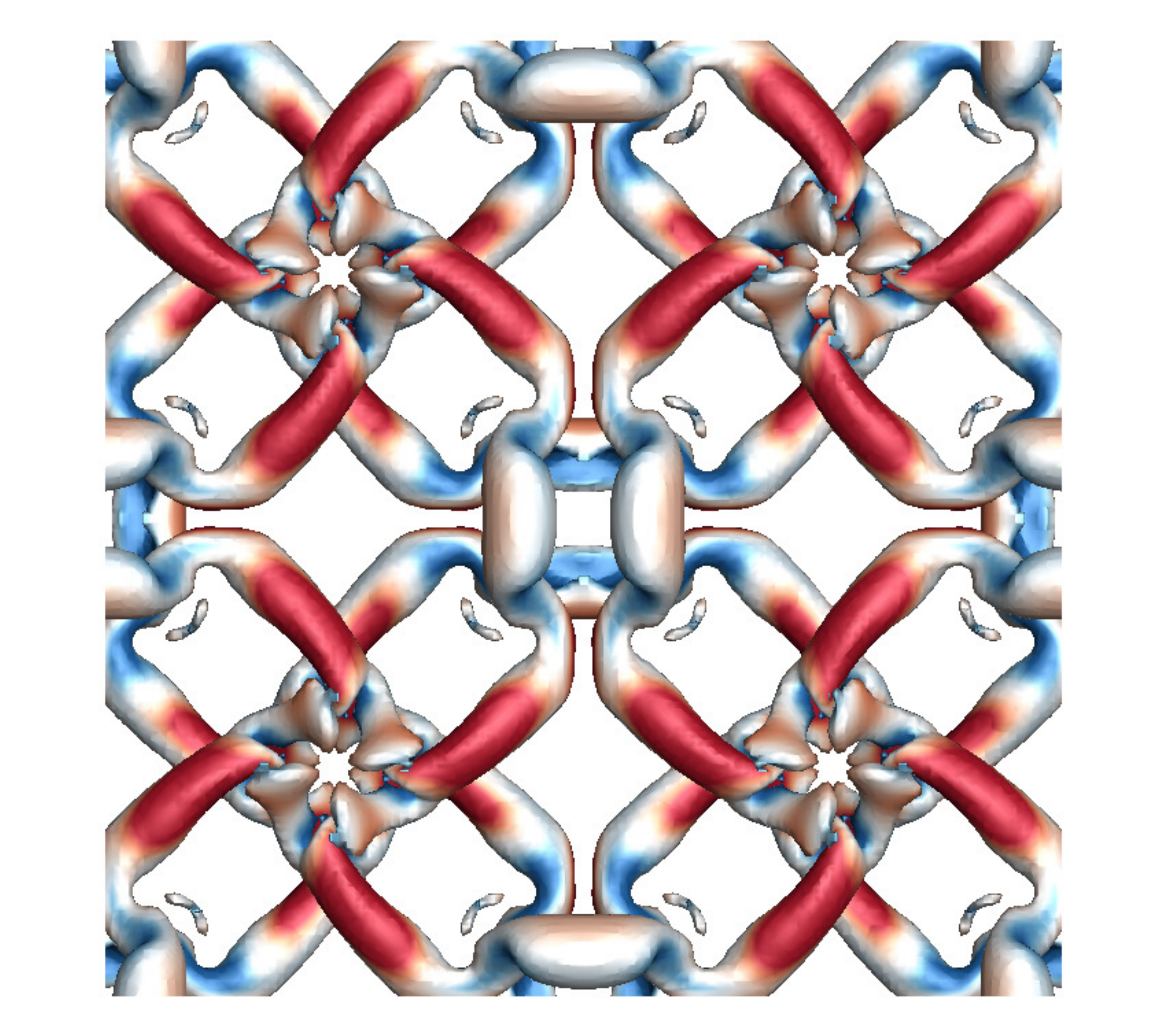}
\caption{\label{tg-vortex4} Taylor-Green Vortex problem: iso-surfaces of $Q$ criterion colored by velocity magnitude at time $t =5, 10$ for $Re=1600$.}
\end{figure}

\subsection{Taylor-Green Vortex}
This problem is aimed at testing the performance of high-order
methods on the direct numerical simulation of a three-dimensional
periodic and transitional flow defined by a simple initial
condition, i.e. the Taylor-Green vortex \cite{Case-Bull,Case-Debonis}. With a uniform temperature field, the
initial flow field is given by
\begin{align*}
u=&V_0\sin(\frac{x}{L})\cos(\frac{y}{L})\cos(\frac{z}{L}),\\
v=&-V_0\cos(\frac{x}{L})\sin(\frac{y}{L})\cos(\frac{z}{L}),\\
w=&0,\\
p=&p_0+\frac{\rho_0V_0^2}{16}(\cos(\frac{2x}{L})+\cos(\frac{2y}{L}))(\cos(\frac{2z}{L})+2).
\end{align*}
The fluid is then a perfect gas with $\gamma=1.4$ and the Prandtl
number is $Pr=0.71$. Numerical simulations are conducted with two
Reynolds numbers $Re=1600$ and $280$. The flow is computed within a
periodic square box defined as $-\pi L\leq x, y, z\leq \pi L$. The
characteristic convective time $t_c = L/V_0$. In the computation,
$L=1, V_0=1, \rho_0=1$, and the Mach number takes $M_0=V_0/c_0=0.1$,
where $c_0$ is the sound speed.

The volume-averaged kinetic energy can be computed from the flow as it evolves in time, which is expressed as
\begin{align*}
E_k=\frac{1}{\rho_0\Omega}\int_\Omega\frac{1}{2}\rho\textbf{u}\cdot\textbf{u}d\Omega,
\end{align*}
where $\Omega$ is the volume of the computational domain, and the dissipation rate of the kinetic energy is given by
\begin{align*}
\varepsilon_k=-\frac{dE_k}{dt}.
\end{align*}
The numerical results of the current scheme with $128\times128\times128$  mesh points for the normalized volume-averaged kinetic energy and dissipation rate with Reynolds numbers $Re=1600$ and $280$ are presented in Fig.\ref{tg-vortex1} and Fig.\ref{tg-vortex3}, which agree well with the data in \cite{Case-Debonis,Case-wang}. The iso-surfaces of $Q$ criterions colored by velocity magnitude at $t=5$ and $10$ are shown in Fig.\ref{tg-vortex2} for $Re=1600$ and in Fig.\ref{tg-vortex4} for $Re=280$. The evolution of flow structure is evident, starting from large vortices and decaying into more complex structures. Different from many other higher-order methods, the current scheme has no internal degrees of freedom to be updated within each cell.

\section{Conclusion}
In this paper, based on the two-stage time stepping method, a fourth-order gas-kinetic scheme is proposed for the three-dimensional inviscid and viscous flow computations.
With the three-dimensional WENO-JS reconstruction, a gas-kinetic scheme with higher-order spatial and temporal accuracy is developed.
In comparison with the classical methods based on the first-order Riemann solver,
for the same fourth-order accuracy in time the current scheme only uses two stages instead of four stages with the Runge-Kutta time-stepping technique. As a result, the two-stage GKS can be more efficient than the four-stage methods with the absence of two time consuming reconstructions.
For the Navier-Stokes solutions, the current scheme doesn't separate inviscid and viscous terms and they are treated uniformly from the same
initial WENO-type reconstruction. The GKS can present very accurate viscous flow solutions due to its multidimensional
flux function at a cell interface, where the gradients in both normal and tangential directions of flow variables participate in the gas evolution.
The fourth-order GKS not only has the expected order of accuracy for the
smooth flow, but also has favorable shock capturing property for the discontinuous solutions.
Most importantly, the numerical tests clearly demonstrate that the current fourth-order scheme has the same robustness as the second-order one.
The scheme has been tested from the smooth flows to the flows with discontinuities, and from the low speed to the hypersonic ones.
For the three dimensional Navier-Stokes solutions, the current scheme is one of the state-of-art methods in the capturing of complicated flow structures.

\section*{Acknowledgements}
The work of Pan is supported by the grants from NSFC (11701038) and China Postdoctoral Science Foundation (2016M600065),
and the work of Xu is supported by Hong Kong research grant council (16206617, 16207715, 16211014)
and the grants from NSFC (11772281, 91530319).

\end{document}